\documentclass[reqno,letter]{amsart}
\usepackage{amsmath,amssymb,epsfig}
\usepackage{float}
\usepackage[pdftex,pdftitle={anexoc},colorlinks=true,breaklinks=true,linkcolor=blue]{hyperref}

\numberwithin{equation}{section}
\newcommand{\diag}{\operatorname*{diag}}
\renewcommand{\aa}{\ensuremath{\mathbf{a}}}

\newcommand{\abs}[1]{\ensuremath{\left|#1\right|}}

\newcommand{\Div}{\nabla\cdot\,}
\newcommand{\ee}{\ensuremath{\mathbf{e}}}
\newcommand{\eps}{\ensuremath{\varepsilon}}
\newcommand\escala{\tilde{\varphi}}

\newcommand{\Grad}{\ensuremath{\mathrm{\nabla}}}

\newcommand{\M}{\ensuremath{\mathcal{M}}}
\newcommand{\N}{\ensuremath{\mathbb{N}}}
\newcommand{\n}{\ensuremath{\mathbf{n}}}
\newcommand{\Om}{\ensuremath{\Omega}}
\newcommand{\pt}{\ensuremath{\partial_t}}

\renewcommand{\r}{\ensuremath{\mathbf{r}}}
\newcommand{\RR}{\mathbb{R}}

\newcommand\wavelet{\tilde{\psi}}

\renewcommand{\u}{\mathbf{u}}

\newcommand{\Iap}{I_{\mathrm{app}}}
\newcommand{\Ion}{I_{\mathrm{ion}}}

\newcommand{\IapK}{I_{\mathrm{app},K}}

\newcommand{\IonK}{I_{\mathrm{ion},K}}

\newcommand{\dx}{\ensuremath{\, dx}}

\newcommand{\dt}{\ensuremath{\, dt}}

\newcommand{\EE}{\mathcal{E}}
\newcommand{\TT}{\mathcal{T}}

\newcommand{\bM}{\ensuremath{\mathbf{M}}}
\newtheorem{alg}{Algorithm}

\newtheorem{defi}{Definition}[section]

\allowdisplaybreaks

\title[Multiresolution scheme for the bidomain model]{A multiresolution space-time
adaptive scheme for the bidomain
model in electrocardiology}

\author[Bendahmane]{Mostafa Bendahmane$^{\mathrm{a}}$}
\author[B\"urger]{Raimund B\"urger$^{\mathrm{a}}$}
\author[Ruiz]{Ricardo Ruiz Baier$^{\mathrm{a}}$}

\thanks{$^{\mathrm{a}}$Departamento de Ingenier\'{\i}a Matem\'{a}tica,
 Facultad de Ciencias F\'{\i}sicas y Matem\'{a}ticas,
 Universidad de Concepci\'{o}n, Casilla 160-C,  Concepci\'{o}n,
 Chile.
E-mail: {\tt mostafab@ing-mat.udec.cl}, {\tt rburger@ing-mat.udec.cl},
{\tt rruiz@ing-mat.udec.cl}}

\date{\today}
\begin{document}

\begin{abstract}
This work deals with the numerical solution of
the  monodomain and bidomain  models of  electrical activity of
myocardial tissue. The bidomain model is a system consisting of
a  possibly degenerate  parabolic PDE coupled with  an elliptic
PDE for the transmembrane  and extracellular potentials,
respectively. This system of two scalar PDEs
  is supplemented by   a time-dependent ODE
  modeling the evolution of the so-called gating variable. In
the simpler sub-case of the monodomain model, the elliptic PDE
reduces to an algebraic equation.  Two simple models for the membrane
and ionic currents are considered, the Mitchell-Schaeffer
 model and  the  simpler FitzHugh-Nagumo model. Since typical solutions
 of the bidomain and monodomain  models  exhibit  wavefronts with steep
gradients, we propose  a finite volume scheme enriched by a fully
adaptive multiresolution method,  whose basic purpose is to
concentrate computational effort on zones of strong variation
 of the solution.  Time adaptivity
is achieved by two alternative devices, namely  locally varying time
stepping and a Runge-Kutta-Fehlberg-type adaptive time  integration.
 A  series of numerical examples
 demonstrates   that
these methods are  efficient and sufficiently accurate
to simulate the electrical activity in myocardial tissue with
affordable effort.  In addition, an optimal
threshold for discarding   non-significant  information
 in the multiresolution representation of the solution   is derived,  and the
numerical efficiency and accuracy of the  method is measured in terms of CPU time
speed-up, memory compression, and errors in different norms.
\end{abstract}



\maketitle

\section{Introduction}

\subsection{Scope}
The obvious  difficulty of performing
direct measurements in
electrocardiology has motivated  wide interest in the   numerical
simulation
 of  cardiac models. In  1952, Hodgkin and Huxley \cite{HH} introduced the
first mathematical  model of wave propagation  in squid nerve,
 which was modified  later on
 to describe  several phenomena in biology. This
  led  to  the first physiological  model of cardiac tissue \cite{Noble} and
many others. Among these models, the \emph{bidomain model}, firstly introduced
 by Tung \cite{tung}, is  one of
the most accurate and complete models for the theoretical and numerical study
of the electric activity in cardiac tissue.  The bidomain equations
result from the principle of conservation of current between the intra- and extracellular
domains,  followed by a  homogenization process (see  e.g.
\cite{Bend-Karl:cardiac,Colli5,keenersneyd}) derived from a scaled version of a
cellular model on a periodic structure of cardiac tissue.  Mathematically,
 the  bidomain model  is   a  coupled system consisting of a
 scalar, possibly degenerate
 parabolic PDE coupled with  a scalar  elliptic PDE for the
transmembrane  potential and the
 extracellular potential, respectively. These equations
 are supplemented by   a time-dependent ODE for
 the so-called gating variable, which is defined at every point of the
 spatial computational domain. Here, the term ``bidomain''  reflects that in
general,  the intra- and extracellular tissues have different longitudinal and
transversal (with respect to the fiber)   conductivities; if these
are equal, then the model is termed {\em monodomain model},
 and  the elliptic PDE reduces to an algebraic equation. The degenerate structure
 of the mathematical formulation of the bidomain model
is essentially due to the differences between the intra- and extracellular
anisotropy of the cardiac tissue \cite{Bend-Karl:cardiac,Colli4}.

The bidomain model represents a computational challenge
 since    the  width  of an excitation front is roughly  two orders of
magnitude smaller than the long axis of a human-size right ventricle. This local
 feature, along with  strongly  varying time scales in the reaction
 terms,  produces  solutions with   sharp propagating wave fronts in the
 potential field, which almost precludes  simulations with  uniform grids.  Clearly,
  cardiac simulations should be based on  space- (and
also time-) adaptive methods.

It is the purpose of this paper  to develop a fully adaptive multiresolution
(MR) scheme with locally varying space-time stepping (LTS) and adaptive time
step control by means of  a Runge-Kutta-Fehlberg (RKF) method. These strategies
are of different nature, but  do not exclude  each
other; rather, they may be combined to obtain a potentially  more powerful
 method (as is suggested e.g.\ in  \cite{DRS062}; however, herein  we do not
 pursue that approach). We furthermore   address the deduction
of an optimal threshold value for discarding non-significant data,
which permits to achieve significant data compression.
 Previous experience  with
degenerate parabolic equations and reaction-diffusion systems
\cite{bbrs,bks,brss,brss2,RS05,RSTB03}  suggests  that the MR device
should provide an efficient tool  for solving the bidomain equations,
and in this same spirit, we construct the corresponding extension of
the MR method with the novel application to the bidomain and
 monodomain models  in mind.

The efficiency of the MR method is a consequence of the fact
 that at each time step,
the solution is encoded  with respect to a MR basis corresponding to a
hierarchy of nested grids. The size of the details  determines the level of
refinement needed to obtain an accurate local representation of the solution.
Therefore,  an adaptive mesh is evolved in time by
refining and coarsening in a suitable way, by means of a strategy based on the
prediction of the displacement and creation of singularities in the solution.

We apply the MR approach to an explicit finite volume (FV) method   in
each time step. Since  the computational effort required  for integrating a
system of equations for one time step is usually  substantially higher
for an implicit scheme when compared to explicit schemes,
implicit schemes  may be less efficient than explicit ones, especially when
the overall number of time steps is large (see e.g. \cite{cgh:03}).

\subsection{Related work}
To further put this work  into the  proper perspective, we first mention that
standard theory for coupled parabolic-elliptic systems (see e.g.
\cite{chen}) does not apply naturally to  the bidomain
equations, since  the anisotropies of the intra- and extracellular media differ
 and the resulting  system  is of degenerate parabolic type.  Colli  Franzone and
Savar\'e~\cite{Colli4} present a weak formulation for the bidomain model and show
that it has a structure suitable for applying the theory of  evolution variational
inequalities in Hilbert spaces. Bendahmane and Karlsen~\cite{Bend-Karl:cardiac}
prove existence and uniqueness for the bidomain equations using
the Faedo-Galerkin method and compactness theory for the existence
part, and  Bourgault et al.\    \cite{coudiere1}  prove existence and uniqueness
for the bidomain equations by  first reformulating the problem as a single
 parabolic PDE,  and then applying a semigroup approach.

From a computational point of view,  substantial contributions have been
made  in adaptivity for cardiac models. However, the approach
presented herein differs to the best of our knowledge from
other adaptive approaches in the literature. These alternative techniques
include  adaptive mesh refinement (AMR)
(e.g., \cite{cgh:03,tk:04}), adaptive finite element methods using a
posteriori error techniques (see, e.g., \cite{Colli5}) or multigrid methods
applied to finite elements. Furthermore, Quan et al.\ \cite{qehh}
 present a domain decomposition approach
using an alternating direction implicit (ADI) method.
With respect to time adaptivity, Skouibine et al.\   \cite{stm:00}  present a
predictor-corrector time stepping strategy to  accelerate a given finite
differences scheme for the bidomain equations using active membrane kinetics
(Luo-Rudy phase II).  Cherry et al.\   \cite{cgh:03}  use  local time stepping,
similar to the method introduced in the germinal work of Berger and Oliger
  \cite{BO}, to  accelerate a reference scheme.
Parallelized versions of part of the methods
mentioned  above are presented, for example,
 by Colli Franzone and Pavarino \cite{Colli-Pava} and Saleheen and Ng
 \cite{Saleheen}.

MR schemes for hyperbolic partial differential equations
 were first proposed
by  Harten \cite{harten:1995}. We refer to the work of M\"uller
\cite{Muller} for a survey on MR methods, see also Chiavassa  et al.\ \cite{chiav}.
As stated above, the  idea behind the MR method is to accelerate a reference
discretization scheme while controlling the error.
In the context of fully adaptive MR methods \cite{Cohen}, the mathematical analysis
is complete only in the case of a scalar conservation law, but in practice, these techniques
have  been used by several groups (see  e.g.\ \cite{bbrs,DRS062,Muller,MS,RSTB03}) to successfully
solve a wide class of problems, including applications to multidimensional systems.
For more details on the framework of classical MR methods for
hyperbolic partial differential equations, we also refer  to
 Cohen et al.\  \cite{Cohen} and  Dahmen et al.\ \cite{DGM}.

\subsection{Outline  of the paper} The remainder of this
paper is organized as follows. In Section~\ref{sec:model}, the
 bidomain and  monodomain  models of cardiac tissue are  introduced.
  The  general bidomain model  can be expressed as a
 coupled system of a parabolic PDE and an elliptic PDE plus an ODE
 for the evolution of the local gating variable,
  while the monodomain model, which arises as a particular
sub-case of the bidomain model, is defined by a reaction-diffusion
 equation, which is again supplemented with an ODE for the gating
 variable. Section~\ref{sec:FV} deals with the construction of an
appropriate FV method for the solution of both
the parabolic--elliptic system  and the reaction--diffusion equation arising
from the  bidomain and monodomain models, respectively. Next, in
 Section~\ref{sec:MR} we develop the MR
analysis used to endow the  reference finite FV schemes
with space adaptivity.  More precisely, in Section~\ref{subsec4.1} we
introduce the wavelet basis underlying the multiresolution
representation with the pertinent projection operator. In
Section~\ref{subsec4.2}, the  prediction operator and the detail
coefficients are introduced. Small detail coefficients
 on fine  levels of resolution may be discarded (this operation is
called thresholding), which allows for substantial data compression.
 In Section~\ref{subsec4.3}, we  recall   the graded tree data structure
used for storage of the numerical solution,   and which is
introduced for ease of navigation. In Section~\ref{sec:err-analysis}
 we outline an error analysis, similar to that conducted  in
 \cite{brss,brss2,RSTB03} and motivated by the rigorous analysis
 of Cohen et al.\  \cite{Cohen},  which justifies the  choice of a
 reference tolerance~$\varepsilon_{\mathrm{R}}$. In turn, this
quantity determines the comparison values $\varepsilon_l$
 used for the thresholding operation  at each level~$l$ of
 multiresolution. The basic goal is to choose the threshold
values in such a way that the resulting multiresolution scheme
has the same order of accuracy as the usual finite volume scheme.

In Section~\ref{sec:time} we address two strategies for the adaptive
 evolution in time of the space-adaptive MR scheme,  namely
the locally varying time stepping (LTS, Section~\ref{sec:LTS}) and
a  variant of  the well-known Runge-Kutta-Fehlberg
(RKF, Section~\ref{sec:RKF}) method.
Finally, in Section~\ref{sec:num_results} we present
  numerical examples putting into evidence the efficiency of the underlying methods.
  Some  conclusions that can be drawn from the paper about the effectiveness of our
  methods and statement of possible further extensions to our research  are given in
  Section \ref{sec:concl}, and in the Appendix we present a brief description of the
 LTS and general MR algorithms.

\section{The macroscopic bidomain and monodomain models}
\label{sec:model}
The spatial domain
for  our  models is a bounded open subset $\Omega \subset \RR^2$ with
a piecewise smooth boundary $\partial \Om$. This represents  a two-dimensional slice of
the cardiac muscle regarded as two interpenetrating and superimposed
(anisotropic) continuous media, namely
 the intracellular ($\mathrm{i}$) and extracellular
($\mathrm{e}$) tissues. These tissues  occupy
  the same two-dimensional area, and   are separated
from each other (and connected at each point) by the cardiac cellular
membrane. The quantities of interest are  \textit{intracellular} and
\textit{extracellular} electric potentials,
$u_\mathrm{i}=u_\mathrm{i}(x,t)$ and $u_\mathrm{e}=u_\mathrm{e}(x,t)$,
at $(x,t)\in \Om_T:=\Om\times(0,T)$. Their  difference
$v=v(x,t):=u_\mathrm{i}-u_\mathrm{e}$   is known as the
\textit{transmembrane  potential}. The conductivity of the tissue is represented by
scaled tensors $\bM_\mathrm{i}(x)$ and $\bM_\mathrm{e}(x)$ given by
$$\bM_j(x)=\sigma_j^{\mathrm{t}}\mathbf{I} +
(\sigma_j^{\mathrm{l}}-\sigma_j^{\mathrm{t}})\aa_l(x)\aa_l^{\mathrm{T}}(x),$$
where $\smash{\sigma_j^{\mathrm{l}} =\sigma_j^{\mathrm{l}}(x)\in C^1(\RR^2)}$ and
$\smash{\sigma_j^{\mathrm{t}}=\sigma_j^{\mathrm{t}}(x)\in
  C^1(\RR^2)}$,
$j \in \{ \mathrm{e}, \mathrm{i} \}$,  are the intra- and extracellular
conductivities  along and transversal to   the direction of the fiber
(parallel to $\aa_l(x)$),  respectively.

For fibers aligned with the axis, $\bM_\mathrm{i}(x)$ and $\bM_\mathrm{e}(x)$ are
diagonal matrices:
$\smash{\bM_\mathrm{i}(x)=\diag
(\sigma_\mathrm{i}^{\mathrm{l}} ,  \sigma_\mathrm{i}^{\mathrm{t}})}$
 and $\smash{\bM_\mathrm{e}(x)= \diag
(\sigma_\mathrm{e}^{\mathrm{l}} ,  \sigma_\mathrm{e}^{\mathrm{t}})}$.
 When  the so-called
 \emph{anisotropy ratios} $\sigma_\mathrm{i}^{\mathrm{l}}/\sigma_\mathrm{i}^{\mathrm{t}}$
and $\sigma_\mathrm{e}^{\mathrm{l}}/\sigma_\mathrm{e}^{\mathrm{t}}$  are
equal, we are in the case of \emph{equal anisotropy}, but generally the conductivities
in the longitudinal direction~$\mathrm{l}$ are higher than those across the fiber
(direction~$\mathrm{t}$);
such a case is called \emph{strong anisotropy} of electrical conductivity. When the
fibers rotate from bottom to top, this type of anisotropy is often
referred to as
\emph{rotational anisotropy}.

The bidomain model is given by the  following coupled
reaction-diffusion system   \cite{sundes,ying}:
\begin{align}\label{S1}
\begin{split}
\beta c_\mathrm{m} \pt v-\Div\bigl(\bM_\mathrm{i}(x)\Grad u_\mathrm{i} \bigr)+\beta
\Ion(v,w)& = 0, \\
\beta c_\mathrm{m} \pt v+\Div\bigl(\bM_\mathrm{e}(x)\Grad u_\mathrm{e}\bigr) +\beta
\Ion(v,w)& =  \Iap,\\
\pt w -H(v,w)& =0, \qquad \qquad  (x, t) \in \Om_T.
\end{split}
\end{align}
Here, $c_\mathrm{m}>0$ is the so-called {\em surface capacitance}
  of the membrane, $\beta$ is the
surface-to-volume ratio, and $w(x,t)$ is the
so-called {\em gating} or {\em recovery variable},
 which also takes into account the concentration
variables. The stimulation currents applied to the extracellular space are
represented by the function  $\Iap=\Iap(x,t)$.
The functions $H(v, w)$ and~$\Ion(v, w)$ correspond to
 the fairly simple Mitchell-Schaeffer membrane model \cite{MS:Ion}
  for the membrane and ionic currents:
\begin{align}
H(v,w)&=\frac{w_\infty(v/v_p)-w}{R_{\mathrm{m}}c_{\mathrm{m}}\eta_\infty(v/v_p)},
\quad
\Ion(v,w) =\frac{v_p}{R_{\mathrm{m}}}\biggl(\frac{v}{v_p\eta_2}-
\frac{v^2(1-v/v_p)w}{v_p^2\eta_1}\biggr),\label{Ion}
\end{align}
where the dimensionless functions~$\eta_{\infty}(s)$  and~$w_{\infty}(s)$
 are  given by
$\eta_\infty(s)= \eta_3 + ( \eta_4- \eta_3) \mathcal{H} (s- \eta_5)$
 and $w_{\infty}(s) = \mathcal{H} (s- \eta_5)$, where
 $\mathcal{H}$
 denotes the Heaviside function,
  $R_{\mathrm{m}}$ is the surface resistivity of the membrane, and
$v_p$ and  $\eta_1,\dots,\eta_5$ are given parameters. A simpler
choice for the membrane kinetics is  the widely known FitzHugh-Nagumo
 model \cite{FH,Nagumo}, which is often used to avoid computational
difficulties arising from  a large number of coupling variables.
This model is specified by
\begin{align} \label{fhn}
H(v,w) =av-bw, \quad
\Ion(v,w)=-\lambda\bigl(w-v(1-v)(v-\theta)\bigr),
\end{align}
where $a$, $b$, $\lambda$ and $\theta$ are given parameters.

We
rewrite  \eqref{S1}   equivalently in terms of~$v$
 and~$u_\mathrm{e}$ as the  strongly coupled parabolic-elliptic PDE-ODE
system   (see e.g. \cite{sundes,ying}):
\begin{subequations}\label{S4}
\begin{align}
\beta c_{\mathrm{m}}\pt v+\Div\bigl(\bM_\mathrm{e}(x)\Grad
u_\mathrm{e}\bigr)+\beta \Ion(v,w)& =
\Iap, \\
\Div \bigl((\bM_\mathrm{i}(x)+\bM_\mathrm{e}(x))\Grad u_\mathrm{e}\bigr)+
\Div\bigl(\bM_\mathrm{i}(x)\Grad v\bigr)&=\Iap, \label{eliptica}\\
\pt w -H(v,w)& =0, \qquad \qquad \qquad  (x, t) \in \Om_T.
\end{align}
\end{subequations}

We utilize zero flux boundary conditions, representing an
isolated piece of cardiac tissue:
\begin{equation}\label{S2}
\bigl(\bM_j(x)\Grad u_j)\cdot\n =0\text{ on }
\Sigma_T:=\partial\Om\times(0,T),\quad j \in \{ \mathrm{e},
\mathrm{i} \},
\end{equation}
and impose   initial conditions (which are degenerate for the
transmembrane potential $v$):
\begin{equation}\label{S3}
v(0,x)=v_0(x),\quad w(0,x)=w_0(x), \quad x \in \Omega.
\end{equation}

For the solution $v$ of the bidomain model, we
require the initial datum $v_0$ to be compatible with \eqref{S2}. Therefore,
if we fix both $u_j(0,x)$, $j\in\{ \mathrm{e},\mathrm{i}\}
$ as initial data, the problem may become unsolvable, since the time
derivative involves only $v=u_\mathrm{i}-u_\mathrm{e}$ (this is also referred as
\emph{degeneracy in time}). Thus,  we impose the compatibility condition
\begin{equation}\label{compat_ue}
\int_\Om u_\mathrm{e}(x,t)\dx=0 \quad \text{for a.e. $t\in(0,T)$.}
\end{equation}

In the case
that $\bM_\mathrm{i}\equiv \lambda \bM_\mathrm{e}$  for some constant
$\lambda \in \RR$,  the system \eqref{S1} is equivalent to a scalar
parabolic equation for the transmembrane potential $v$, coupled to an ODE
for the gating variable~$w$. This parabolic equation is obtained by
multiplying the first equation
in \eqref{S1} by $1/(1+ \lambda)$, the second by
 $\lambda/(1+ \lambda)$ and adding the resulting equations.
 The final
\emph{monodomain model} can be stated as follows:
\begin{align}\label{monodomain}
 \begin{split}
 \displaystyle \beta c_{\mathrm{m}}\pt v-\Div
 \biggl(\frac{\bM_\mathrm{i}}{1+\lambda}\Grad v \biggr)
 +\beta \Ion(v,w)& = \frac{\lambda}{1+\lambda}\Iap , \\
 \displaystyle \pt w -H(v,w)&=0,  \qquad  (x, t) \in  \Om_T.
 \end{split}
 \end{align}
This model is, of course, significantly  less involved and requires substantially less
computational effort than the bidomain model, and even though the assumption
of equal anisotropy ratios is very strong and generally unrealistic, the
monodomain model is adequate for a qualitative investigation of
repolarization sequences and the distribution of
 patterns of  durations of the action potential
  \cite{Colli3}.


We  assume that the functions $\bM_j$, $j\in\{\mathrm{e,i}\}
$,  $\Ion$,  and $H$ are sufficiently smooth so that the following
definitions of weak solutions make sense. Furthermore, we assume that
$\Iap \in L^2(\Om_T)$ and
$\smash{\bM_j\in L^\infty(\Om)}$ and
$\smash{\bM_j\boldsymbol{\xi}\cdot\boldsymbol{\xi}\geqslant
C_M|\boldsymbol{\xi}|^2}$
for a.e. $x\in\Om$, for all  $\boldsymbol{\xi}\in \RR^2$, $j\in
\{\mathrm{e,i}\}$,
and a constant $C_M>0$. For later reference, we now state the
definitions
 of a weak solution for the bidomain and the monodomain model, respectively.

\begin{defi} \label{def1} A triple $\u=(v,u_\mathrm{e},w)$
 of functions is a {\em weak solution of the bidomain model}
 \eqref{S4}--\eqref{S3} if  $v,u_\mathrm{e}\in L^2(0,T;H^1(\Om))$,
  $w\in C([0,T],L^2(\Om))$,  \eqref{compat_ue} is satisfied, and the
  following identities hold for all test functions
  $\varphi,\psi,\xi\in \mathcal{D}([0,T)\times\bar{\Om})$:
\begin{align*}
\beta c_{\mathrm{m}}\int_\Om v_0(x)\varphi(0,x)\dx+
 \iint_{\Om_T} \Bigl\{
\beta c_{\mathrm{m}} v\pt
\varphi -  \bM_\mathrm{e}(x)\Grad u_\mathrm{e}\cdot \Grad \varphi
 + \beta \Ion \varphi \Bigr\}
\dx\dt
&=\iint_{\Om_T}\Iap
\varphi \dx\dt, \\
\iint_{\Om_T}\Bigl\{ - (\bM_\mathrm{i}(x)+\bM_\mathrm{e}(x))\Grad u_\mathrm{e}\cdot
\Grad \psi
- \bM_\mathrm{i}(x)\Grad v\cdot \Grad \psi\Bigr\}  \dx\dt&=\iint_{\Om_T}\Iap
\varphi \dx\dt,\\
-\int_\Om w_0(x)\xi(0,x)\dx-\iint_{\Om_T}w\pt\xi \dx \dt & =\iint_{\Om_T}H\xi \dx \dt.
\end{align*}
\end{defi}

\begin{defi} A pair $\u=(v,w)$  of functions is a {\em weak solution
    of the monodomain model} \eqref{monodomain} if $v\in L^2(0,T;H^1(\Om))$,
$w\in C([0,T],L^2(\Om))$, and the
  following identities hold for all test functions
  $\varphi,\xi\in \mathcal{D}([0,T)\times\bar{\Om})$:
\begin{align*}
 \beta c_{\mathrm{m}}\int_\Om v_0(x)\varphi(0,x)\dx & \\
 +\iint_{\Om_T} \biggl\{ \beta c_{\mathrm{m}}
 v\pt \varphi +\beta \Ion \varphi
  - \frac{1}{1+\lambda}\bM_\mathrm{i}\Grad v\cdot \Grad \varphi \biggr\} \dx\dt
 & =\frac{\lambda}{1+\lambda}\iint_{\Om_T}\Iap \varphi \dx\dt,
   \\
 -\int_\Om w_0(x)\xi(0,x)\dx-\iint_{\Om_T}w\pt\xi \dx \dt
 & =\iint_{\Om_T}H\xi \dx \dt.
 \end{align*}
\end{defi}

\section{The reference finite volume scheme}\label{sec:FV}
To define the  FV scheme for approximating
solutions to the bidomain equations \eqref{S4}, we follow
the framework of \cite{Ey-Gal-Her:book}. An admissible mesh for $\Om$ is formed
by a family $\TT$ of control volumes (open and convex
polygons) of maximum diameter $h$. (From the next section on, we will use only
Cartesian meshes, however the following description is in a more general
setting.)  For all $K \in \TT$, $x_K$ denotes the
center of~$K$, $N(K)$~is the set of neighbors of $K$, $\EE_{\text{int}}(K)$ is
the set of edges of $K$ in the interior of $\TT$,  and $\EE_{\text{ext}}(K)$
the set of edges of $K$ on the boundary $\partial \Om$. For all $L \in
N(K)$,
$d(K,L)$ denotes the distance between~$x_K$ and~$x_{L}$, $\sigma_{K,L}$ is
the interface between~$K$ and~$L$, and $\eta_{K,L}$
($\eta_{K,\sigma}$,
respectively) is the unit normal vector to $\sigma_{K,L}$
($\sigma \in \EE_{\text{ext}}(K)$,  respectively) oriented from $K$ to $L$
(from $K$ to $\partial \Om$,  respectively). For all $K \in \TT$, $|K|$
stands for the measure of the cell $K$. The admissibility of $\TT$ implies
that $\overline{\Om}=\cup_{K\in \TT} \overline{K}$, $K\cap L=\varnothing$
if $K,L\in \TT$ and $K \ne L$,  and there exist a finite sequence
$(x_{K})_{K\in \TT}$ for which  $\overline{x_{K}x_{L}}$ is orthogonal to
$\sigma_{K,L}$.

Now, consider $K \in \TT$ and $L \in N(K)$ with common vertices
$(a_{\ell,K,L})_{1\le \ell\le I}$ with $I \in \N \backslash \{0\}$ and
let $T_{K,L}$ ($T^{\text{ext}}_{K,\sigma}$ for
$\sigma\in \EE_{\text{ext}}(K)$, respectively) be the open and convex polygon
with vertices $(x_K,x_L)$ ($x_K$,  respectively) and
$(a_{\ell,K,L})_{1\le \ell\le I}$. Notice that $\Om$ can be decomposed into
$\smash{\overline{\Om}=\cup_{K\in \TT}
((\cup_{L\in N(K)}\overline{T}_{K,L})\cup
 (\cup_{\sigma \in
\EE_{\text{ext}}(K)}\overline{T}^{\text{ext}}_{K,\sigma}))}$.

For all $K \in \TT$, the approximation $\Grad_h u_{h}$ of $\Grad u$
is defined by
$$
\Grad_h u_{h}(x)
:=
\begin{cases}
	{\displaystyle \frac{\displaystyle
	|\sigma_{K,L}|}{\displaystyle |T_{K,L}|}}(u_{L}-u_{K})\eta_{K,L}
	& \text{if $x \in T_{K,L}$}, \\
	0
	& \text{if $x \in T^{\text{ext}}_{K,\sigma}$}.
\end{cases}
$$

To discretize \eqref{S4}--\eqref{S3}, we choose an
  admissible discretization  of $Q_T$, consisting of
an admissible mesh of $\Om$ and  a time step size $\Delta t>0$.
 For example, we could choose  $N>0$
as the smallest integer such that $N\Delta t\ge T$, and  set
$t^n:=n \Delta t$ for $n\in  \{0,\ldots,N\}$. However,
  for reasons stated further below,   we assume that the
time step size $\Delta t$ is determined anew for each iteration,
and present the scheme as it applies to advance the
solution from $t^n$ to $t^{n+1}:= t^n + \Delta t$.

On each cell $K \in \TT$, (positive definite) conductivity tensors are defined
by
$$
\bM_{j,K}=\frac{1}{|K|}\int_\Om \bM_{j}(x)\dx,\quad
j \in \{ \mathrm{e}, \mathrm{i} \}.
$$
Let $F_{j,K,L}$ be an approximation of
\begin{align*}
\int_{\sigma_{K,L}}\bM_j(x)\Grad u_j \cdot \eta_{K,L}\,
  d\gamma
\end{align*}
for $j \in \{ \mathrm{e}, \mathrm{i} \}
$, and for $K \in \TT$ and $L \in N(K)$, let
$$
M_{j,K,L}=\abs{\frac{1}{|K|} \int_K \bM_j(x)\dx \,\eta_{K,L}}
:=\abs{\bM_{j,K}\,\eta_{K,L}},\quad j \in \{ \mathrm{e}, \mathrm{i} \}
.
$$
Here $\abs{\cdot}$ stands for the Euclidean norm. The diffusive fluxes
$\bM_j(x)\Grad u_j \cdot \eta_{K,L}$ on $\sigma_{K,L}$
are  approximated by
\begin{align*}
& \int_{\sigma_{K,L}} \bigl(\bM_j(x)\Grad u_j\bigr) \cdot \eta_{K,L} \, d\gamma
\\ & \approx |\sigma_{K,L}|\Grad u_{j}(y_\sigma)
\cdot (\bM_{j,K}\, \eta_{K,L})
= |\sigma_{K,L}|M_{j,K,L} \Grad u_{j}(y_\sigma)
\cdot \frac{y_\sigma-x_K}{d(K,\sigma_{K,L})}
\approx |\sigma_{K,L}|M_{j,K,L}
\frac{u_{j,\sigma}-u_{j,K}}{d(K,\sigma_{K,L})},
\end{align*}
where $y_\sigma$ is the center of $\sigma_{K,L}$ and $u_{j,\sigma}$ is an
approximation of $u_j(y_\sigma)$, $j \in \{ \mathrm{e}, \mathrm{i} \}
$. The resulting
approximation of fluxes is consistent \cite{Ey-Gal-Her:book}. In
addition,  the scheme should be conservative. This property enables
us to determine the additional unknowns $u_{j,\sigma}$, and to compute the
numerical fluxes on internal edges:
\begin{align}
F_{j,K,L}=d_{j,K,L}^*
\frac{|\sigma_{K,L}|}{d(K,L)}(u_{j,L}-u_{j,K}) \quad \text{if $L\in N(K)$},
\end{align}
where we define
$$
d_{j,K,L}^*:=\frac{M_{j,K,L}M_{j,L,K}}
{d(K,\sigma_{K,L})M_{j,K,L}+d(L,\sigma_{K,L})M_{j,L,K}} d(K,L),
$$
while we discretize the zero-flux boundary condition by setting
\begin{equation}\label{bound-edges}
F_{j,K,\sigma}=0
\end{equation}
on boundary edges. Here, $d(K,\sigma_{K,L})$ and $d(L,\sigma_{K,L})$ are
the distances from $x_K$ and $x_{L}$ to $\sigma_{K,L}$, respectively.
We define cell averages of the unknowns $H(v, w)$ and $\Ion(v, w)$:
\begin{align*}
H_{K}^{n+1}:=\frac{1}{\Delta t |K|}\int_{t^n}^{t^{n+1}}
\int_KH \bigl(v(x,t),w(x,t) \bigr)\dx\dt, \quad
\IonK^{n+1}:=\frac{1}{\Delta t |K|}\int_{t^n}^{t^{n+1}}
\int_K\Ion \bigl(v(x,t),w(x,t) \bigr)\dx\dt,
\end{align*}
and of the given function $\Iap$:
\begin{align*}
\IapK^{n+1}& :=\frac{1}{\Delta t |K|}\int_{t^n}^{t^{n+1}}
\int_K\Iap(x,t)\dx\dt.
\end{align*}
The computation starts from the initial cell averages
\begin{align}\label{prob:init}
v_K^0=\frac{1}{|K|} \int_{K} v_0(x) \dx,\quad
w_K^0=\frac{1}{|K|} \int_{K} w_0(x) \dx.
\end{align}

We now describe the finite  volume scheme
employed to advance the numerical solution from $t^n$ to $t^{n+1}$,
 which is based on  a simple
explicit Euler time discretization. Assuming that
 at $t=t^n$,
 the quantities  $\smash{u_{j,K}^{n}}$,
   $j\in \{ \mathrm{e}, \mathrm{i} \}$,
 $\smash{v_{K}^{n}=(u_{\mathrm{i},K}^{n}-u_{\mathrm{e},K}^{n})}$,
  and $w_{K}^{n}$ are known for all
 $K \in \TT$,  we compute the values of these cell averages
  $\smash{u_{j,K}^{n+1}}$,
   $j\in \{ \mathrm{e}, \mathrm{i} \}$,
 $\smash{v_{K}^{n+1}=(u_{\mathrm{i},K}^{n+1}-u_{\mathrm{e},K}^{n+1})}$
  and $w_{K}^{n+1}$ at $t=t^{n+1}$ from
\begin{align}
\label{S1-discr}
\beta c_{\mathrm{m}} |K|\frac{v^{n+1}_K-v^{n}_K}{\Delta t}+\sum_{L \in N(K) }
     d_{\mathrm{e},K,L}^*\frac{|\sigma_{K,L}|}{d(K,L)}
(u^{n}_{\mathrm{e},L}-u^{n}_{\mathrm{e},K})
     +\beta |K|\IonK^n
& =|K|\IapK^{n},
\\
\label{S2-discr}
\sum_{L \in N(K) } \frac{|\sigma_{K,L}|}{d(K,L)}    \left\{
\bigl(d_{\mathrm{i},K,L}^*+d_{\mathrm{e},K,L}^*\bigr)
\bigl(u^{n+1}_{\mathrm{e},L}-u^{n+1}_{\mathrm{e},K}\bigr)
+ 
 d_{i,K,L}^*\bigl(v^{n+1}_L-v^{n+1}_K\bigr)\right\} &=|K|\IapK^{n},
\\ \label{S3-discr}
     |K|\frac{w^{n+1}_K-w^{n}_K}{\Delta t}-|K|H^{n}_K&=0.
\end{align}
The order in which these equations are used to advance
is explicitly stated
in Algorithm~\ref{alg:general} in Section~\ref{sec:num_results}.

The boundary condition \eqref{S2} is taken into account by imposing
zero fluxes on the external edges as in \eqref{bound-edges}:
\begin{equation}\label{moyenne-discr}
d_{j,K,\sigma}^*\frac{|\sigma_{K,L}|}{d(K,L)}(u^{n}_{j,L}-u^{n}_{j,K})=0
\quad \text{ for }\sigma \in \EE_{\text{ext}}(K),
\quad j \in \{ \mathrm{e}, \mathrm{i} \}, \quad n=0,1,2,\dots.
\end{equation}
and the compatibility condition \eqref{compat_ue} is discretized via
\begin{equation*}
\sum_{K\in \TT}|K|u^n_{\mathrm{e},K}=0, \quad  n =0,1,2, \dots.
\end{equation*}

Analogously, a FV method for the monodomain model \eqref{monodomain}
is given by determining vectors $v_{K}^{n+1}$
and $w_{K}^{n+1}$ for $K \in \mathcal{T}$ and $n =0,1,2, \dots$  such that for
all $K \in \TT$,  we start from the
 initial data given by  \eqref{prob:init},  and use the following
 formulas to advance the solution over one time step:
\begin{align*}
      \beta c_{\mathrm{m}} |K|\frac{v^{n+1}_K-v^{n}_K}{\Delta t}+\sum_{L \in N(K)
 }
 \frac{1}{1+\lambda}d_{\mathrm{i},K,L}^*\frac{|\sigma_{K,L}|}{d(K,L)}(v^{n}_{L}-v^{n}_{K})
       +\beta |K|\IonK^n &
 =\frac{\lambda}{1+\lambda}|K|\IapK^{n}, \\
   |K|\frac{w^{n+1}_K-w^{n}_K}{\Delta t}-|K|H^{n+1}_K&=0.
 \end{align*}

A FV method for a slightly different version of the bidomain equations is
analyzed in  \cite{Ben-Karl:FV-bidomain}. In that paper,
 the authors prove existence and
uniqueness of solutions to an implicit FV scheme, and provide convergence
results. On the other hand, following
\cite{Ben-Karl:FV-bidomain,coudiere2,Ey-Gal-Her:book}, we prove
in \cite{bbrconvergence}
existence and  uniqueness of approximate solutions
 (that is, well-definedness of the scheme)
 for an  implicit version of the scheme considered herein,
 and show that it converges  to a weak solution
in the sense of Definition~\ref{def1}, under a mild condition
 that limits the time step size~$\Delta t$ (a CFL-type condition
 is, however, not imposed in    \cite{bbrconvergence}).
%
Moreover, as in~\cite{bbrs}, we may  deduce that
in the case of Cartesian meshes,
the
explicit version of the FV method  utilized herein,
\eqref{prob:init}--\eqref{moyenne-discr} ,
is stable under the CFL condition
\begin{align}\label{cfl}
\Delta t\leqslant h \Bigl(
2\displaystyle{\max_{K\in\TT}\bigl(|\IonK^n|+ |\IapK^n|\bigr)}+
4h^{-1}\max_{K\in\TT}\bigl(|M_{\mathrm{i},K}|+|M_{\mathrm{e},K}|\bigr) \Bigr)^{-1}.
\end{align}
Notice that the values of   $\smash{\IonK^n}$ and $\smash{ \IapK^n}$
 depend on time. However, while $\Iap$  is a given control function
 for our model and therefore $\smash{\max_{K\in\TT} |\IapK^n|}$ can
 assumed to be bounded, the quantity
 $\IonK^n$ is not bounded a priori for
   arbitrarily large  times.
Consequently, in our computations,
 we evaluate the right-hand side of   \eqref{cfl} after each
iteration at $t=t^n$, and use
 \eqref{cfl} to {\em define} the time step size $\Delta t$
 to advance the solution
 from $t^n$ to   $t^{n+1}=t^n + \Delta t$.

\section{Multiresolution and Wavelets}\label{sec:MR}
\subsection{Wavelet basis} \label{subsec4.1}
Consider a rectangle which after a change of variables can be regarded as
$\Om=[0,1]^2$. We determine a nested mesh hierarchy
\mbox{$\Lambda_0\subset\cdots\subset\Lambda_L$}, using an uniform dyadic
partition of $\Om$. Here each grid $\smash{\Lambda_l:=\{V_{(i,j),l}\}_{(i,j)}}$,
with $(i,j)$ to be defined, is formed by the control volumes at  each level
$V_{(i,j),l}:=2^{-l}[i,i+1]\times[j,j+1]$, $i,j\in I_l=\{0,\ldots,2^l-1\}$,
$l=0,\ldots,L$. Here,  $l=0$ corresponds to the coarsest and $l=L$ to the
finest level. The nestedness of the grid hierarchy is made precise  by the
refinement sets $\M_{(i,j),l}=\{2(i,j)+\ee\}$, $\ee\in E:=\{0,1\}^2$,
 which satisfy
$\#\M_{(i,j),l}=4$.  For $x=(x_1,x_2)\in V_{(i,j),l}$ the \emph{scale box
function} is defined as
\begin{equation*}
\escala_{(i,j),l}(x): =\frac{1}{|V_{(i,j),l}|}\chi_{V_{(i,j),l}}(x)=
2^{2l}\chi_{[0,1]^2}(2^lx_1-i,2^lx_2-j),
\end{equation*}
and the averages of any function $u(\cdot,t)\in L^1(\Om)$ for the cell
$V_{(i,j),l}$ may be expressed equivalently as the inner product
$\smash{u_{(i,j),l}:=\langle
  u,\escala_{(i,j),l}\rangle_{L^1(\Om)}}$.
 We are now ready to define the following  two-scale relation for cell
averages and box functions:
\begin{align}\label{wav-2scale-u}
\begin{split}
\escala_{(i,j),l}&=\sum_{\r\in\M_{(i,j),l}}\frac{|V_{\r,l+1}|}{|V_{(i,j),l}|}
\escala_{\r,l+1} = \frac{1}{4} \sum_{(p,q) \in E} \escala_{(2i+p,
  2j+q), l+1},\\
\bar{u}_{(i,j),l}&=\sum_{\r\in\M_{(i,j),l}}\frac{|V_{\r,l+1}|}{|V_{(i,j),l}|}
u_{\r,l+1} =  \frac{1}{4} \sum_{(p,q) \in E} u_{(2i+p,2j+q), l+1},
\end{split}
\end{align}
which   defines a \emph{projection} operator, which
allows us to move from finer to coarser levels. For $x\in V_{2(i,j)+\aa,l+1}$
with $\aa\in E$, we define the \emph{wavelet function} depending on the box
functions on a finer level
\begin{align*}
\wavelet_{(i,j),\ee,l}&=\sum_{\aa\in E}2^{-2}(-1)^{\aa\cdot\ee}
\escala_{2(i,j)+\aa,l+1}
=\sum_{\r\in\M_{(i,j),l}}
\frac{|V_{\r=2(i,j)+\aa,l+1}|}{|V_{(i,j),l}|}(-1)^{\aa\cdot\ee}
\escala_{\r,l+1}.
\end{align*}

The number of related wavelets is $\#\mathcal{M}_{(i,j),l}-1=3$.
Since $\r\cdot\ee\in\{0,1,2\}$  for $\r,\ee\in E$, we have for instance that
$$\wavelet_{(i,j),(1,0),l}=\frac{1}{4}\left(\escala_{2(i,j),l+1}
+\escala_{2(i,j)+(0,1),l+1}-\escala_{2(i,j)+(1,0),l+1}
-\escala_{2(i,j)+(1,1),l+1}\right).$$
Doing this for all $\ee\in E^*:=E\setminus\{(0,0)\}$ yields  an
inverse two-scale relation (see \cite{Muller}), namely
\begin{equation*}
\escala_{2(i,j)+\aa,l+1} =\sum_{\ee\in E}(-1)^{\aa\cdot\ee}
\wavelet_{(i,j),\ee,l} ,\quad \aa\in E.
\end{equation*}
This equation is related to the concept of \emph{stable completions}
 \cite{Muller}. Roughly speaking, the $L^\infty$-counterparts of the
wavelet functions $\smash{\{\wavelet_{(i,j),l}\}_{i,j\in I_l}}$ form a
completion of the $L^\infty$-counterpart of the basis system
$\smash{\{\escala_{(i,j),l}\}_{i,j\in I_l}}$,  and this determines
the existence of a biorthogonal system.
\subsection{Detail coefficients} \label{subsec4.2}
For $\ee\in E^*$, we introduce the \emph{details},  which  will be crucial
to detect zones with steep gradients:
$\smash{d_{(i,j),\ee,l}:=\langle u,\wavelet_{(i,j),\ee,l}\rangle}$.
 These detail coefficients also satisfy a two-scale relation, namely
\begin{equation}\label{wav-2scale-det}
d_{(i,j),\ee,l}=\frac{1}{4}\sum_{2(i,j)
+\aa\in\M_{(i,j),l}}(-1)^{\aa\cdot\ee}u_{2(i,j)+\aa,l+1}.
\end{equation}
An appealing feature is that we can determine a transformation between
the cell averages on level $L$ and the cell averages on level 0 plus
a series of details. This can be achieved by applying recursively the
two-scale relations \eqref{wav-2scale-u} and \eqref{wav-2scale-det};
but we also require  this transformation to be reversible:
\begin{align}\label{tm_a}
\tilde{u}_{(i,j),l+1}=\sum_{\r\in \bar{S}^l_{(i,j)}}g^l_{(i,j),\r}u_{\r,l},
\quad \bar{S}^l_{(i,j)}:=\bigl\{V_{([i/2]+r_1,[j/2]+
r_2),l}\bigr\}_{r_1,r_2\in\{-s,\ldots,0,\ldots,s\}},
\end{align}
where
$\smash{\bar{S}^l_{(i,j)}}$ is
the stencil of interpolation or \emph{coarsening set},
$\smash{g_{(i,j),\r}^l}$ are
coefficients,  and the tilde over $u$ in the left-hand side
 of \eqref{tm_a} means that this corresponds to a predicted value.

Relation \eqref{tm_a} defines
the so-called \emph{prediction} operator, which allows us to move from coarser
to finer resolution levels.  In contrast to  the projection, the prediction
operator is not unique, but we will impose two constraints: to be consistent
with the projection, in the sense that the prediction operator
  is the \emph{right inverse}
of the projection operator, and to be local, in the sense that the predicted value
depends only on $\smash{\bar{S}^l_{(i,j)}}$. For sake of notation, in our case we may
write \eqref{tm_a} as
 \begin{equation*}
 \tilde{u}_{(2i+e_1,2j+e_2),l+1}=u_{(i,j),l}-(-1)^{e_1}Q_x-(-1)^{e_2}Q_y
 +(-1)^{e_1e_2}Q_{xy},
 \end{equation*}
 where $e_1,e_2\in \{ 0,1\}$ and
 \begin{align*}
 Q_x&:=\sum_{n=1}^s\tilde{\gamma}_n \bigl(u_{(i+n,j),l}-u_{(i-n,j),l}\bigr),\quad
 Q_y:=\sum_{p=1}^s\tilde{\gamma}_p \bigl(u_{(i,j+p),l}-u_{(i,j-p),l}\bigr), \\
 Q_{xy}&:=
 \sum_{n=1}^s\tilde{\gamma}_n\sum_{p=1}^s\tilde{\gamma}_p \bigl(u_{(i+n,j+p),l}-
      u_{(i+n,j-p),l}-u_{(i-n,j+p),l}+u_{(i-n,j-p),l}\bigr).
 \end{align*}
Here the corresponding coefficients are
$\tilde{\gamma}_1=\smash{-\frac{22}{128}}$ and
$\tilde{\gamma}_2=\smash{\frac{3}{128}}$ (see \cite{RSTB03}).

From \cite{DGM} we know that details are related to the regularity of a given
function: if $u$ is sufficiently smooth, then its detail coefficients decrease
when going from coarser to finer levels:
\begin{equation*}
\bigl|d^u_{(i,j),l}\bigr|\leqslant C2^{-2lr}\|\Grad^{(r)} u\|_{L^\infty(V_{(i,j),l})},
\end{equation*}
where $r=2s+1$ is the number of vanishing moments of the wavelets. This means that the
more regular $u$ is over $V_{(i,j),l}$, the smaller is the
corresponding detail coefficient. In view of this property, it is
natural to
 attain data compression by
 discarding  the information corresponding to small details. This is called
\emph{thresholding}. Basically, we discard all the elements
corresponding to details  that are smaller in absolute
 value than  a level-dependent
tolerance~$\varepsilon_l$,
$$\bigl|d^u_{(i,j),l}\bigr|<\eps_l.$$
Given a reference tolerance~$\eps_{\mathrm{R}}$, which is determined
by means
 of an error analysis (see Section~\ref{sec:err-analysis}), we
determine $\eps_l$ by
\begin{align} \label{eq4.4}
 \eps_l=2^{2(l-L)}\eps_{\mathrm{R}}.
\end{align}

\subsection{Graded tree data structure} \label{subsec4.3}
We  organize the cell averages and corresponding details at
different levels in a \emph{dynamic graded tree}: whenever
an element is included in the tree, all other elements corresponding to
the same spatial region in coarser resolution levels are also
included,  and
neighboring cells will differ by  at most one refinement level. This
choice guarantees the stability of the multiscale operations
 \cite{Cohen}. We  denote by  {\em root}  the basis  and by \emph{node}
an element of the tree. In two space dimensions, a parent node has four
  sons, and the sons of the same parent are called
 \emph{brothers}. A node without sons is called a \emph{leaf}.
 A given node has $s'=2$ nearest neighbors in each spatial direction, called
 \emph{nearest cousins}, needed for the computation of the fluxes of leaves;
 if these nearest cousins do not exist,  we create them as \emph{virtual
leaves}.  The leaves of the tree are the control volumes forming the adaptive
mesh.  We denote by~$\Lambda$ the set of all nodes of the tree
 and by $\mathcal{L}(\Lambda)$ the restriction of $\Lambda$ to
the leaves. We  apply this MR representation  to the spatial part of the
function $\u=(v,u_\mathrm{e},w)$,
which corresponds to the numerical solution of the underlying problem for each
time step, so we need to update the tree structure  for the
proper representation of the solution during the evolution. To this end, we
 apply  a thresholding strategy, but always keep the graded tree structure
of the data. Once the  thresholding is performed, we add to the tree a
\emph{safety zone}, so the new tree may contain the adaptive mesh
for the next time step. The  safety zone is generated  by adding one finer level
to the tree in all possible positions without violating the graded tree data
structure. This device, first  proposed by Harten \cite{harten:1995},  ensures that
the graded tree  adequately represents  the solution in the next time step.
 Its effectiveness depends strongly on the assumption of finite propagation
 speed of the singularities.

Note that the fluxes are only computed at level $l+1$ and we set the
ingoing flux on the leaf  at level $l$ equal to the sum of the outgoing
fluxes on the leaves of  level $l+1$ sharing the same edge
 \begin{equation}\label{flux_conserv}
 F_{(i+1,j),l\to(i,j),l}=F_{(2i+1,2j),l+1\to(2i+2,2j),l+1}+
 F_{(2i+1,2j+1),l+1\to(2i+2,2j+1),l+1}.
 \end{equation}
It is known that this choice decreases the number of costly flux evaluations without
 loosing the conservativity in the flux computation, and this represents a real
 advantage when using a graded tree structure, see
 e.g. \cite{RSTB03} for more details. This advantage is lost for a non-graded
tree structure, for which fluxes for leaves on an immediately finer level are
not always available.

The \emph{data compression rate}  \cite{brss,brss2}
$\smash{\eta:= \mathcal{N}/(2^{-(L+1)}\mathcal{N}+\#\mathcal{L}(\Lambda))}$
is used to  measure the improvement in data compression.
Here,  $\mathcal{N}$ is the number of elements in the full finest grid at
level $L$, and $\#\mathcal{L}(\Lambda)$ is the size of the set of leaves.
We also measure the \emph{speed-up} $\mathcal{V}$ between the CPU time of
the numerical solution obtained by the FV method and the CPU time of
the numerical solution obtained by the MR method:
$\smash{\mathcal{V}:=\mathrm{CPU\, time}_{\mathrm{FV}} /
 \mathrm{CPU\, time}_{\mathrm{MR}}}$.

\subsection{Error analysis of the multiresolution scheme}
\label{sec:err-analysis}
Using the main properties of the reference FV scheme,
such as the CFL stability condition and
order of approximation in space, we can derive the optimal choice
for the threshold parameter $\eps_{\mathrm{R}}$ for the adaptive
MR scheme. We can decompose the global error between the cell
average values of the exact solution vector at the level $L$, denoted by
$\u^L_{\mathrm{ex}}=(v^L_{\mathrm{ex}},u^L_{\mathrm{e,ex}},
w^L_{\mathrm{ex}})$,
and those of the MR computation with
a maximal level $L$, denoted by $\u^L_{\mathrm{MR}}$, into two errors
\begin{align*}
\bigl\|\u^L_{\mathrm{ex}} - \u^L_{\mathrm{MR}} \bigr\| \leq
\bigl\|\u^L_{\mathrm{ex}} - \u^L_{\mathrm{FV}} \bigr\| +
\bigl\|\u^L_{\mathrm{FV}} - \u^L_{\mathrm{MR}} \bigr\|.
\end{align*}
The first error on the right-hand side is called \emph{discretization
 error} and the second \emph{perturbation error}. Using this, the CFL condition
\eqref{cfl}, and the fact that $h=|\Om|^{1/2}2^{-L}$;
we obtain  that if the so-called \emph{reference tolerance} for the numerical
computations in Section~\ref{sec:num_results} is set  to
\begin{align}
\varepsilon_{\mathrm{R}} &= C \frac{2^{-(\alpha+2)L}}{\displaystyle{
|\Om|\max_{K\in\TT}\Bigl(|\IonK|+2|\IapK|\Bigr)+|\Omega|^{3/2}2^{2+L}\,
\max_{K\in\TT}\Bigl(|M_{i,K}|+|M_{e,K}|\Bigr)}}\label{equ:epsref1},
\end{align}
then we can expect the discretization error and the perturbation error
 to be  kept at
the same order (see \cite{brss,Cohen,RSTB03} for more details).

To measure errors between a reference solution $u_{\mathrm{ex}}$ and an approximate
solution $u_{\mathrm{MR}}$, we will use $L^p$-errors:
$e_p=\smash{ \|u^n_{\mathrm{ex}}-u^n_{\mathrm{MR}}\|_p}$, $p=1,2,\infty$, where
\begin{gather*}
e_\infty=\max_{(i,j,l)\in \mathcal{L}(\Lambda)}\bigl|{u^n_{\mathrm{ex}}}_{i,j,l}-
{u^n_{\mathrm{MR}}}_{i,j,l} \bigr|; \quad
e_p=\Biggl(\frac{1}{|\mathcal{L}(\Lambda)|}\sum_{(i,j,l)\in \mathcal{L}(\Lambda)}
\bigl|{u^n_{\mathrm{ex}}}_{i,j,l}-{u^n_{\mathrm{MR}}}_{i,j,l} \bigr|^p\Biggr)^{1/p},\quad p=1,2.
\end{gather*}
Here $\smash{{u^n_{\mathrm{ex}}}_{i,j,l}}$ stands for the projection of the reference
solution onto  the leaf $(i,j,l)$.
\section{Time--step accelerating methods}\label{sec:time}

\subsection{Local time stepping}\label{sec:LTS}
We employ a version of the locally varying time
stepping strategy introduced by M\"{u}ller and Stiriba
\cite{MS}, and summarize here its principles. The basic
idea is to enforce  a local CFL condition by using the same CFL number for
all scales, and  evolving all leaves on level
$l$ using the local time step size
\begin{equation*}
\Delta t_l=2^{L-l}\Delta t,\quad l=L-1,\ldots,0,
\end{equation*}
where $\Delta t=\Delta t_L$ corresponds to the time step size on the finest
 level $L$. This strategy allows to increase the time step for the major
part of the adaptive mesh without violating  the CFL stability condition.
The synchronization of the time stepping for the portions of the solution
lying on different resolution levels will be automatically achieved
after $2^l$ time steps using $\Delta t_l$. To additionally save
computational effort, the  the tree  is updated  only
each odd intermediate time step $1,3,\ldots,2^L-1$, and furthermore,  the
projection and prediction operators are performed only on scales occupied
by the leaves of the current tree. For the rest of the intermediate time
steps, we  use the current (old) tree structure. For the sake of
synchronization and conservativity of the flux computation, for coarse levels
(scales without leaves), we use \emph{the same} diffusive fluxes and sources
computed in the previous intermediate time step, because the cell average
on these scales are the same that in the previous intermediate time step.
Only for scales containing leaves, we compute fluxes in the following way: if
there is a leaf at the corresponding edge and at the same resolution
level~$l$, we simply perform a flux computation using the brother leaves,
and the virtual leaves at the same level if necessary; and if there is a
leaf at the corresponding cell edge but on a finer resolution level $l+1$
(\emph{interface edge}), the flux will be determined as in
\eqref{flux_conserv}, i.e., we compute the fluxes  at a level $l+1$ on the
same edge, and we set the ingoing flux on the corresponding edge at level
$l$ equal to the sum of the outgoing fluxes on the sons cells of level
$l+1$ (for the same edge). To always have at hand the computed
fluxes as in \eqref{flux_conserv}, we need to perform the locally varying
time stepping recursively from fine to coarse levels.

\subsection{A Runge-Kutta-Fehlberg method}\label{sec:RKF}
In order to upgrade the FV scheme described in Section~\ref{sec:FV}
to at least second order so that the second-order spatial  accuracy
is effective, we  utilize an RKF method \cite{fehlberg},  which,
apart from providing the necessary accuracy, also allows  an
adaptive  control of the time step. For our models, we consider  a
vector-valued RKF method, i.e., $\u=(v,u_\mathrm{e},w)$ and its
time-discretized form at step~$m$, denoted by~$\u^m$. For ease of  discussion,
we assume that the  problem is written as $\pt\u = \mathcal{A} (t, \u)$.

We use two Runge-Kutta methods, of orders $p=3$ and $p-1=2$
 \begin{align*}
 \hat{\u}^{m+1}=\u^m+ \hat{b}_1 \bar{\kappa}_1
  +\hat{b}_2\bar{\kappa}_2 +\hat{b}_3\bar{\kappa}_3,\quad
 \check{\u}^{m+1}=\u^m+ \check{b}_1\bar{\kappa}_1
 +  \check{b}_2\bar{\kappa}_2
 +  \check{b}_3\bar{\kappa}_3,
 \end{align*}
 where
\begin{align}\label{kappas}
\begin{split}
\bar{\kappa}_1&:= \Delta t \mathcal{A}(t^m, \u^m),\\
\bar{\kappa}_2&:= \Delta t \mathcal{A}(t^m +c_2\Delta t,
\u^m+a_{21}\bar{\kappa}_1),\\
\bar{\kappa}_3 &:= \Delta t \mathcal{A}(t^m + c_3\Delta
 t,\u^m+a_{31}\bar{\kappa}_1+a_{32}\bar{\kappa}_2),
 \end{split}
 \end{align}
 and the coefficients  corresponding to  the RK3(2) method are
 $c_2= a_{21}=1$,
 $c_3=\frac{1}{2}$, $a_{31}= a_{32}=\frac{1}{4}$,
 $\smash{\hat{b}_1=\hat{b}_2=\frac{1}{6}}$,
 $\smash{\hat{b}_3=\frac{2}{3}}$,
 $\smash{\check{b}_1= \check{b}_2=\frac{1}{2}}$, and
$\smash{\check{b}_3=0}$.
  These  values  yield  an optimal pair of
  embedded TVD-RK methods of orders two and three.
The truncation  error between the two approximations for
 $\u^{m+1}$ is estimated by
 \begin{equation}\label{delta_old}
 \bar{\delta}_{\mathrm{old}}:=\hat{\u}^{m+1}-\check{\u}^{m+1}=
 \sum_{i=1}^{p}(\hat{b}_i-\check{b}_i)\bar{\kappa}_i\sim (\Delta t)^p,
 \quad \delta_{\mathrm{old}}:=\|\bar{\delta}_{\mathrm{old}}\|_\infty.
 \end{equation}
 Then we can adjust the step size to achieve a prescribed
 accuracy $\delta_{\mathrm{desired}}$ in time. The new time step is
 determined by
$\smash{   \Delta t_{\mathrm{new}}=\Delta t_{\mathrm{old}} |
 \delta_{\mathrm{desired}} / \delta_{\mathrm{old}} |^{1/p}}$
  with $p=3$. To avoid excessively large time steps,  we define
a  limiter function
 $\smash{\mathcal{S}(t):=(\mathcal{S}_0-\mathcal{S}_{\min})\exp (-t / \Delta t
 ) +\mathcal{S}_{\min}}$,
 where we choose $\mathcal{S}_0=0.1$ and $\mathcal{S}_{\min}=0.01$. The new time step
$\Delta t_{\mathrm{new}}$ is then defined as
\begin{align}\label{Delta t_new}
\Delta t_{\mathrm{new}}=\begin{cases}
\Delta t_{\mathrm{old}}|\delta_{\mathrm{desired}} / \delta_{\mathrm{old}} |^{1/p}
& \textrm{if $|(\Delta t_{\mathrm{new}}-\Delta t_{\mathrm{old}})/
\Delta t_{\mathrm{old}} |\leq \frac{1}{2}
\mathcal{S}(t,\Delta t_{\mathrm{old}})$},\\
  \frac{1}{2}
\mathcal{S}(t,\Delta t_{\mathrm{old}})\Delta t_{\mathrm{old}}
+\Delta t_{\mathrm{old}} & \text{otherwise.} \end{cases}
\end{align}

Notice that $\Delta t_{\mathrm{new}}$ is the time step size  for
computing $\u^{m+2}$. More details on the RKF scheme and its
implementation can be found in \cite{brss,DRS06}.
\section{Numerical Examples}\label{sec:num_results}
\begin{table}[t]
\begin{center}
\begin{tabular}{lcccccc}
\hline
Time & $\mathcal{V}$  & $\eta$&   $L^1-$error
& $L^2-$error &$L^\infty-$error $\vphantom{\int^X}$   \\
\hline
$t=0.0\,\mathrm{ms}$
& & 170.22 &
$4.31\times10^{-4}$&$2.47\times10^{-4}$&$3.99\times10^{-4}$\\
$t=1.5\,\mathrm{ms}$
& 27.81 & 37.56& $4.97\times10^{-4}$&$1.96\times10^{-4}$&$4.63\times10^{-4}$
\\
$t=3.5\,\mathrm{ms}$
& 26.47 & 29.89& $5.23\times10^{-4}$&$4.05\times10^{-4}$&$4.82\times10^{-4}$
\\
$t=4.5\,\mathrm{ms}$
& 31.41 & 28.12& $7.48\times10^{-4}$&$4.29\times10^{-4}$&$5.31\times10^{-4}$
\\
$t=5.5\,\mathrm{ms}$
& 30.62 & 24.70& $1.04\times10^{-3}$&$6.20\times10^{-4}$&$6.79\times10^{-4}$
\\
\hline
\end{tabular}
\end{center}

\vspace*{2mm}

\caption{Example~1 (Monodomain model):  Corresponding simulated time,
CPU ratio~$\mathcal{V}$, compression rate~$\eta$  and normalized errors for $v$,
 using a MR method.}
\label{table:ex1}
\end{table}

\begin{figure}[t]
\begin{center}
\begin{tabular}{cc}
\includegraphics[width=0.48\textwidth]{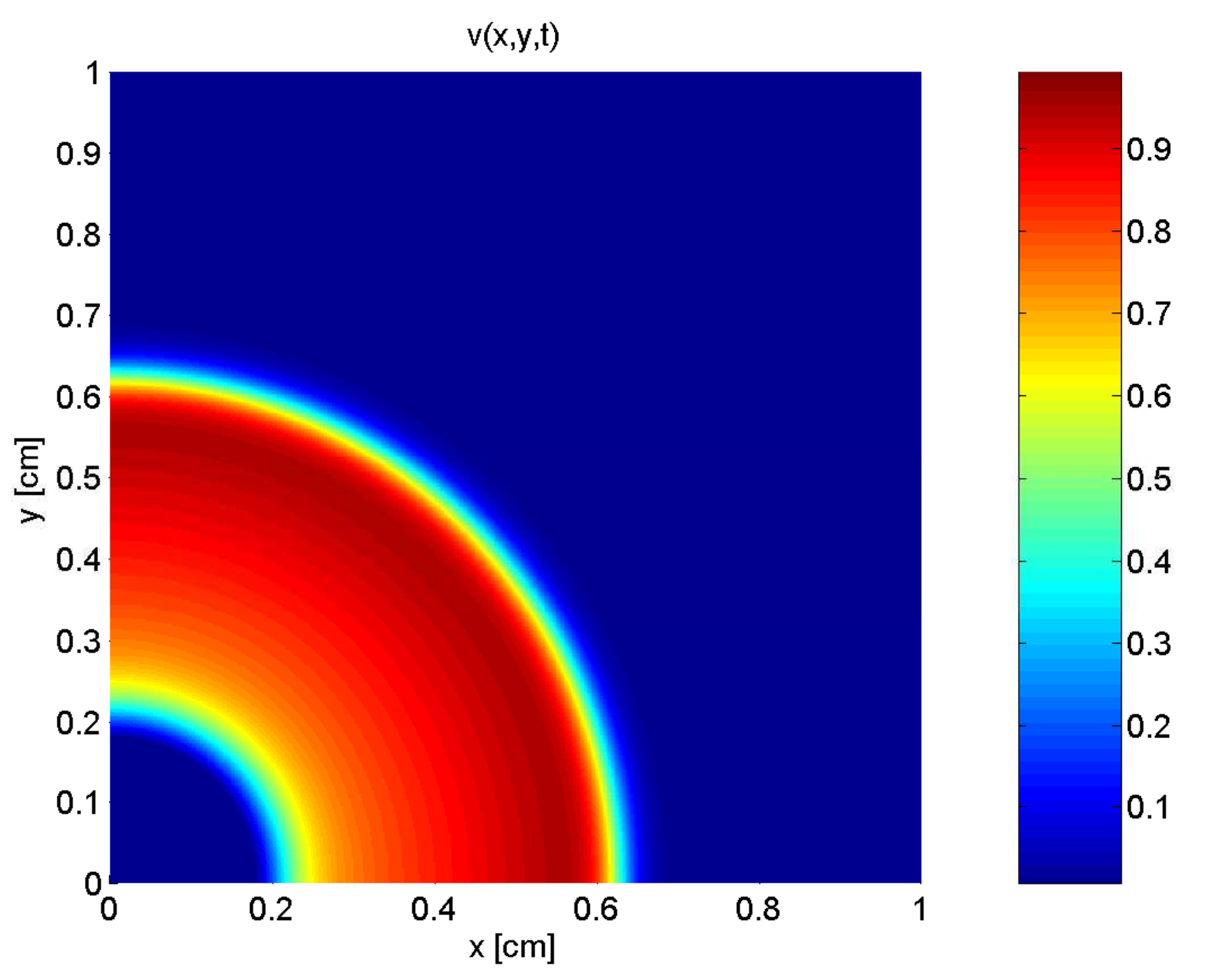}&
\includegraphics[width=0.36\textwidth]{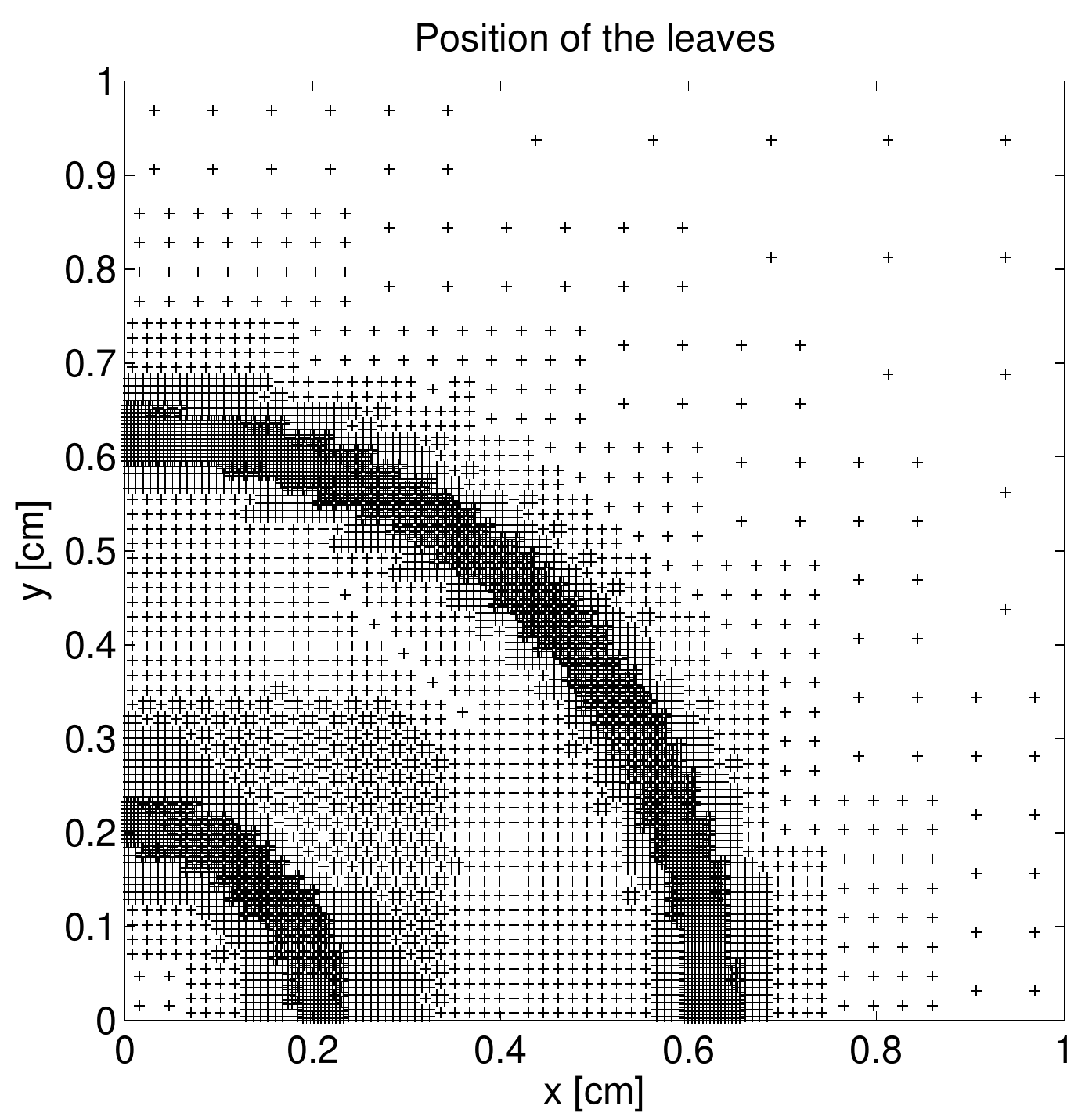}\\
\includegraphics[width=0.48\textwidth]{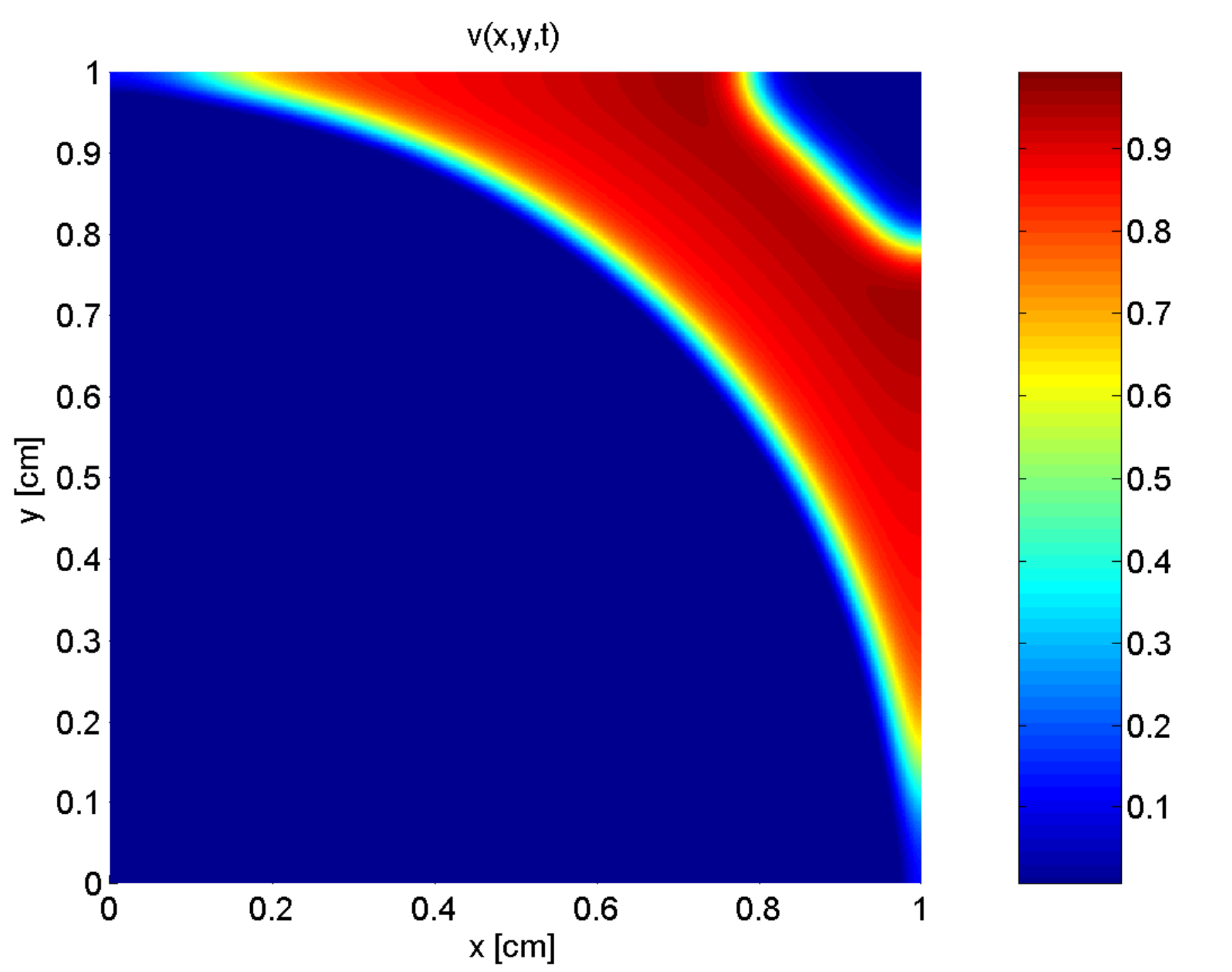} &
\includegraphics[width=0.36\textwidth]{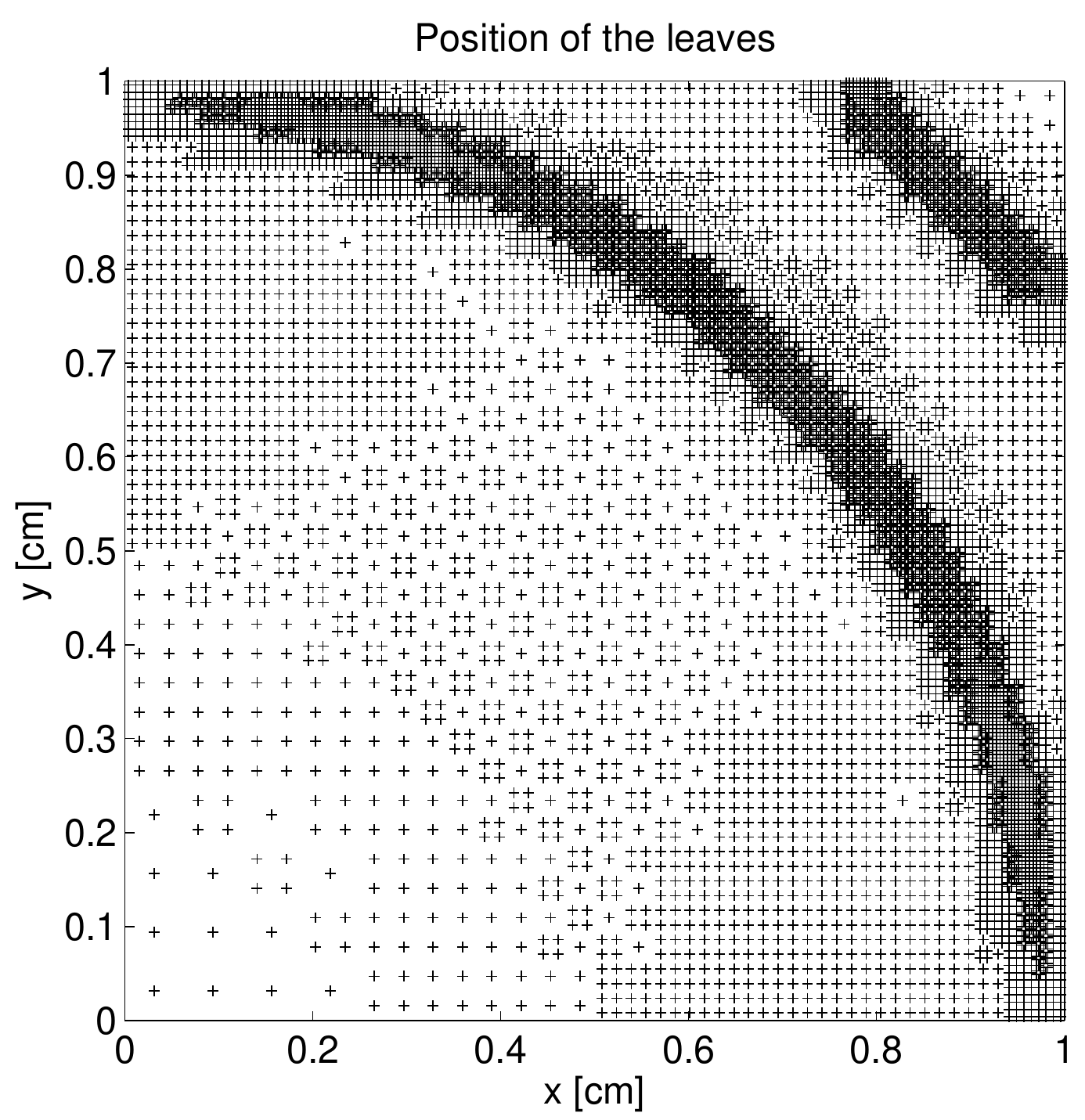}
\end{tabular}
\caption{Example~1 (monodomain model):
Numerical solution for $v$, measured in $[\mathrm{mV}]$ (left) and
leaves of the corresponding
tree at times (from top to bottom) 
$t=1.5\,\mathrm{ms}$, $t=3.5\,\mathrm{ms}$.} \label{fig:monodomain1}
\end{center}
\end{figure}

\begin{figure}[t]
\begin{center}
\begin{tabular}{cc}
\includegraphics[width=0.48\textwidth]{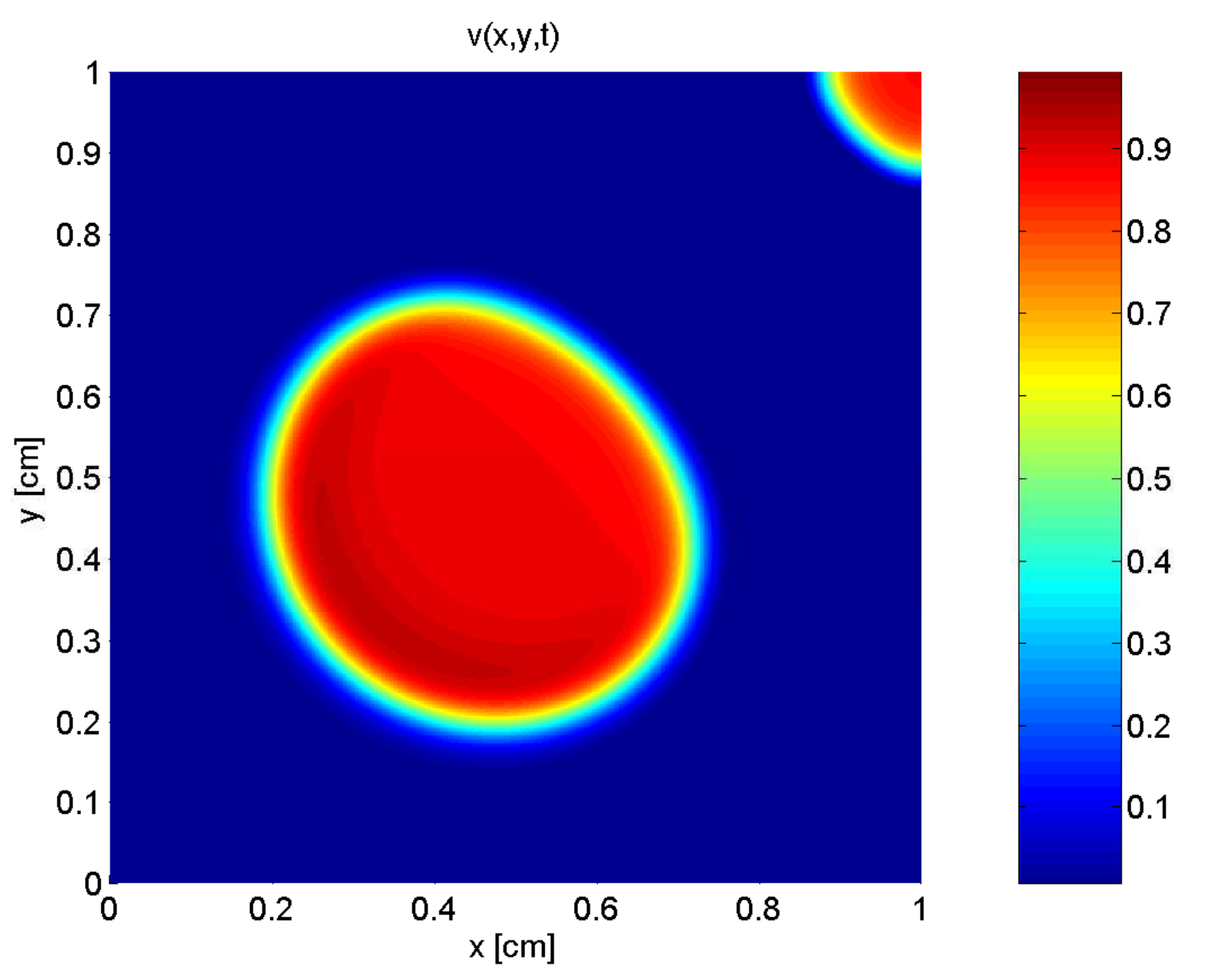}&
\includegraphics[width=0.36\textwidth]{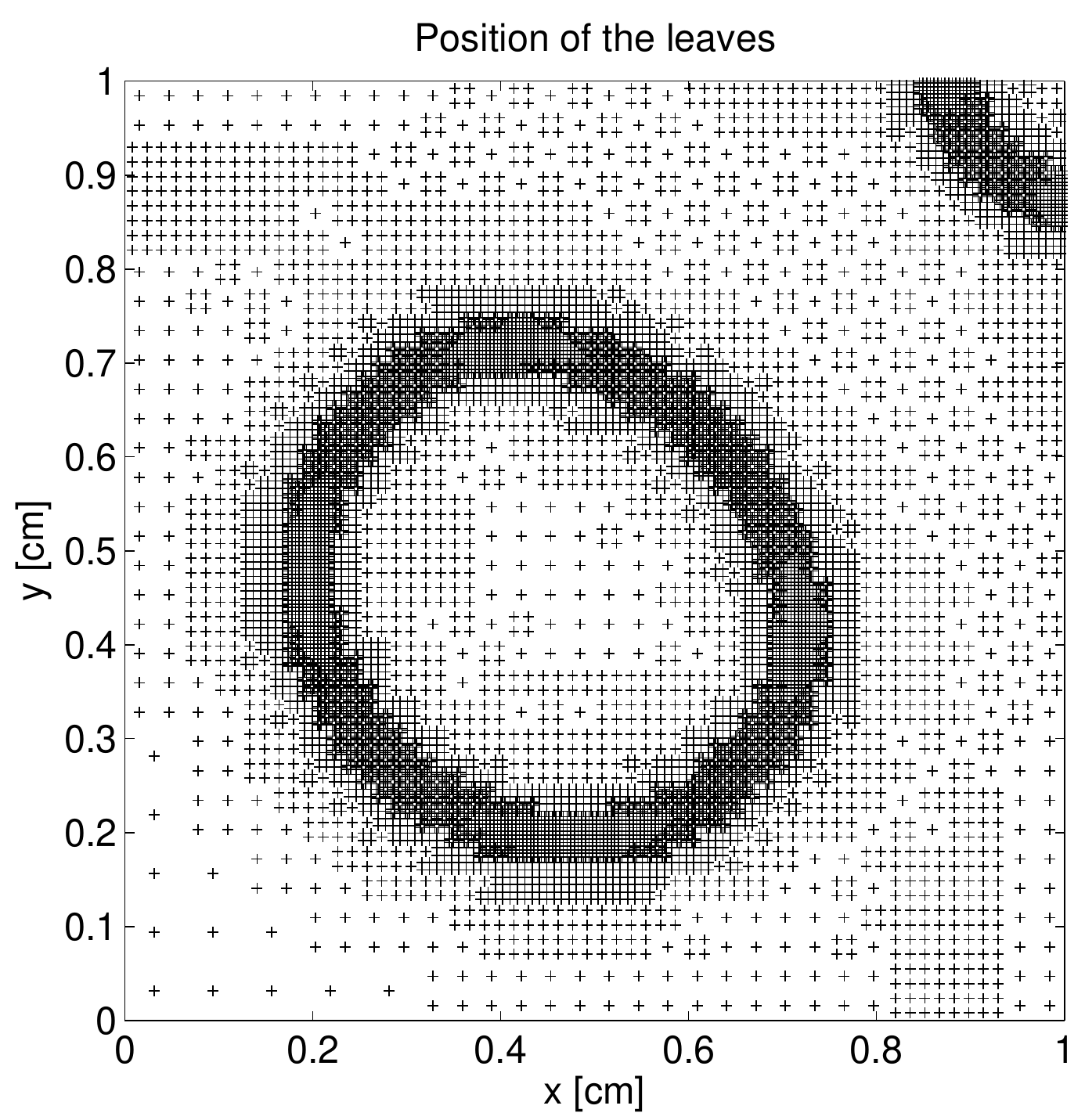}\\
\includegraphics[width=0.48\textwidth]{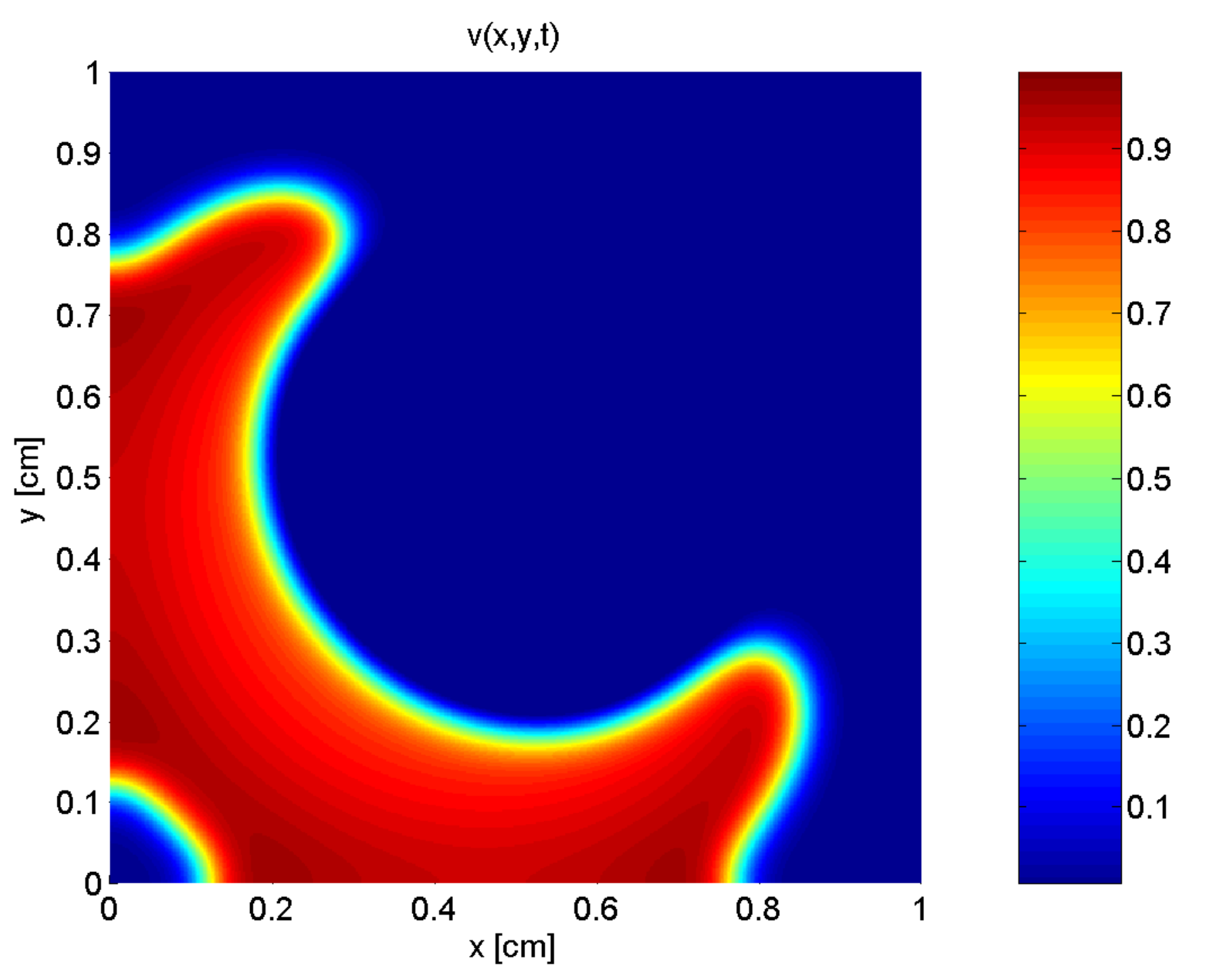}&
\includegraphics[width=0.36\textwidth]{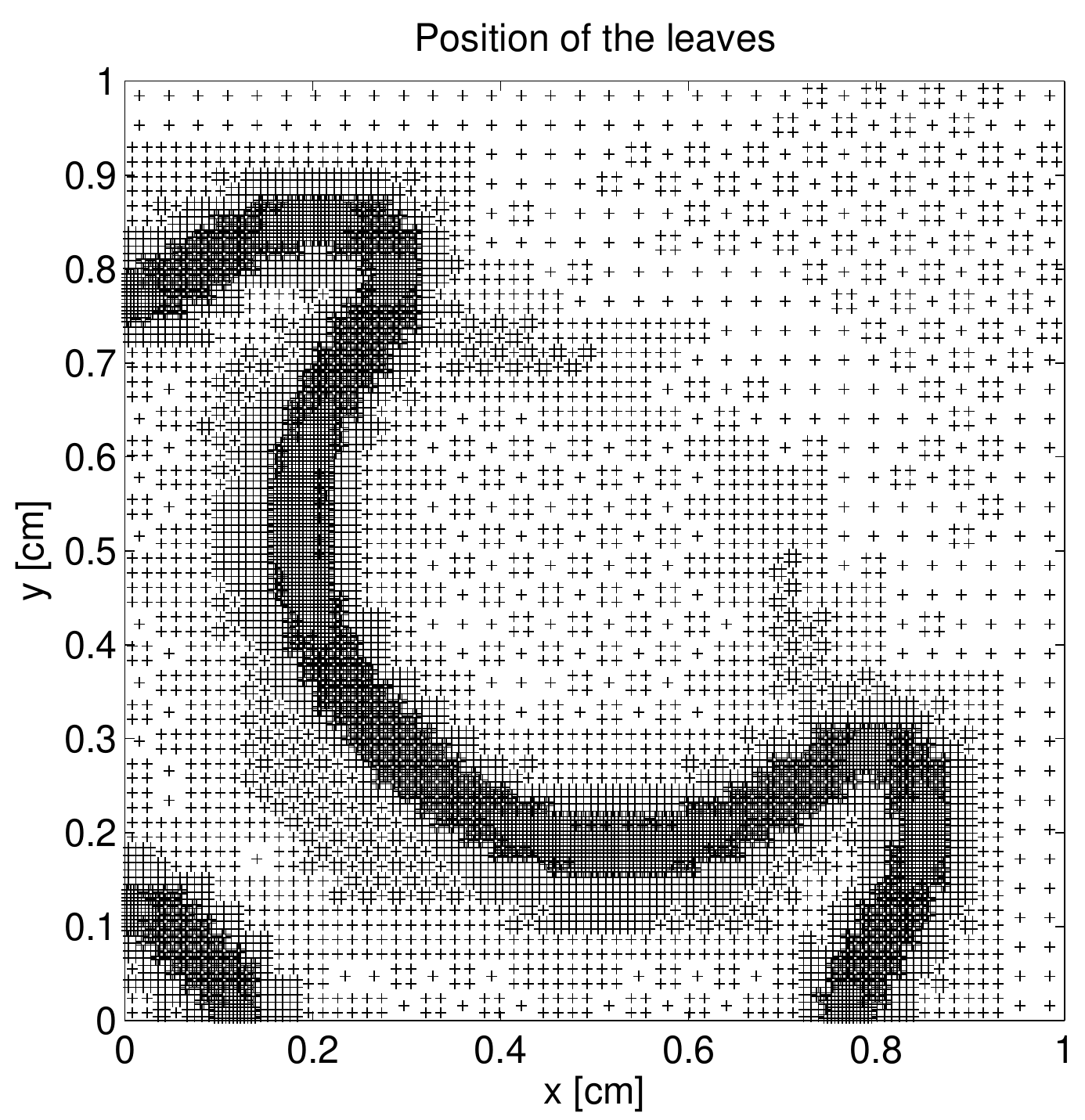}
\end{tabular}
\caption{Example~1 (monodomain model):
Numerical solution for $v$, measured in $[\mathrm{mV}]$ (left) and
leaves of the corresponding tree at times (from top to bottom)
$t=4.5\,\mathrm{ms}$ and $t=5.5\,\mathrm{ms}$.} \label{fig:monodomain2}
\end{center}
\end{figure}

 \begin{figure}[t]
 \begin{center}
 \begin{tabular}{cc}
 \includegraphics[width=0.48\textwidth]{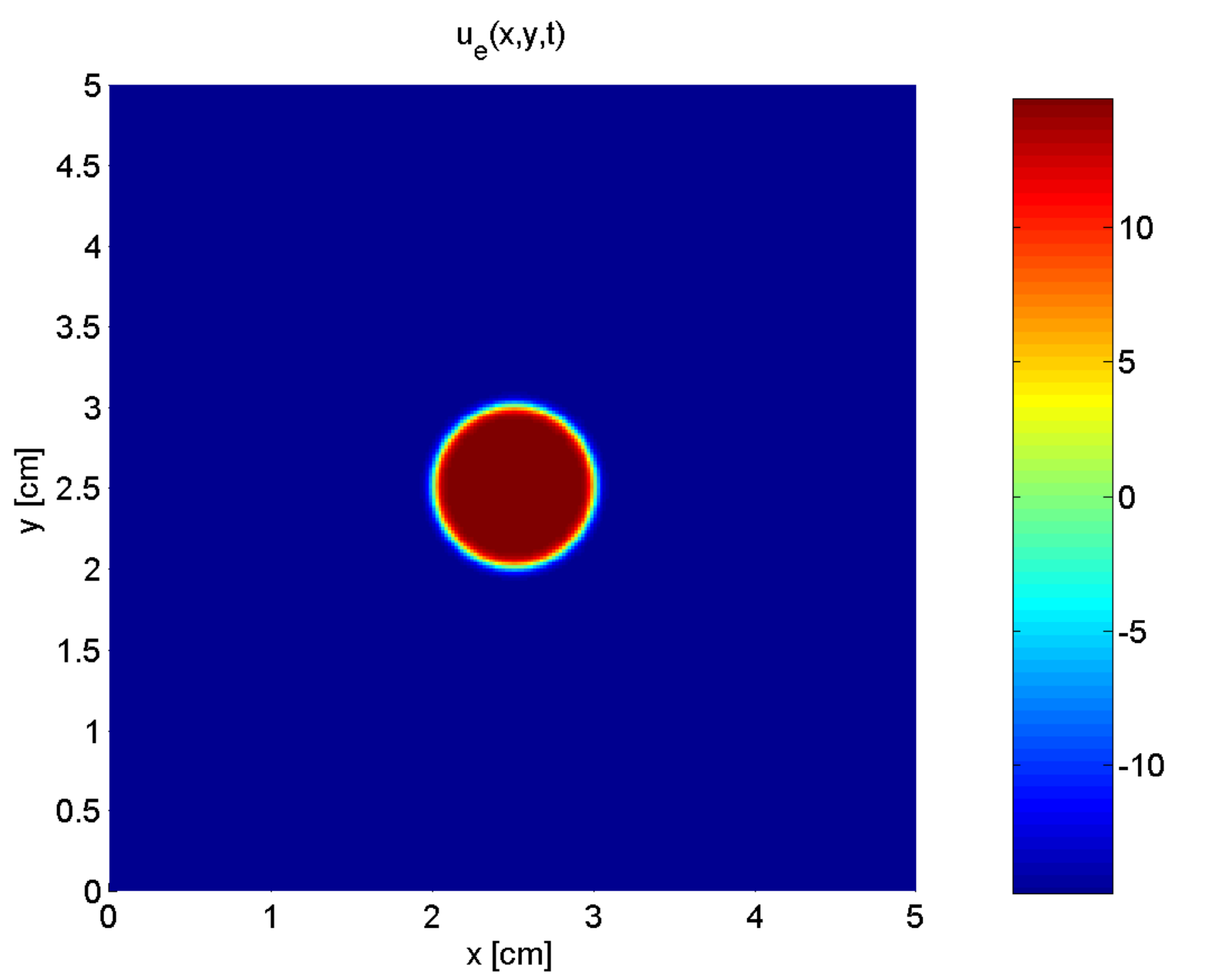}&
 \includegraphics[width=0.36\textwidth]{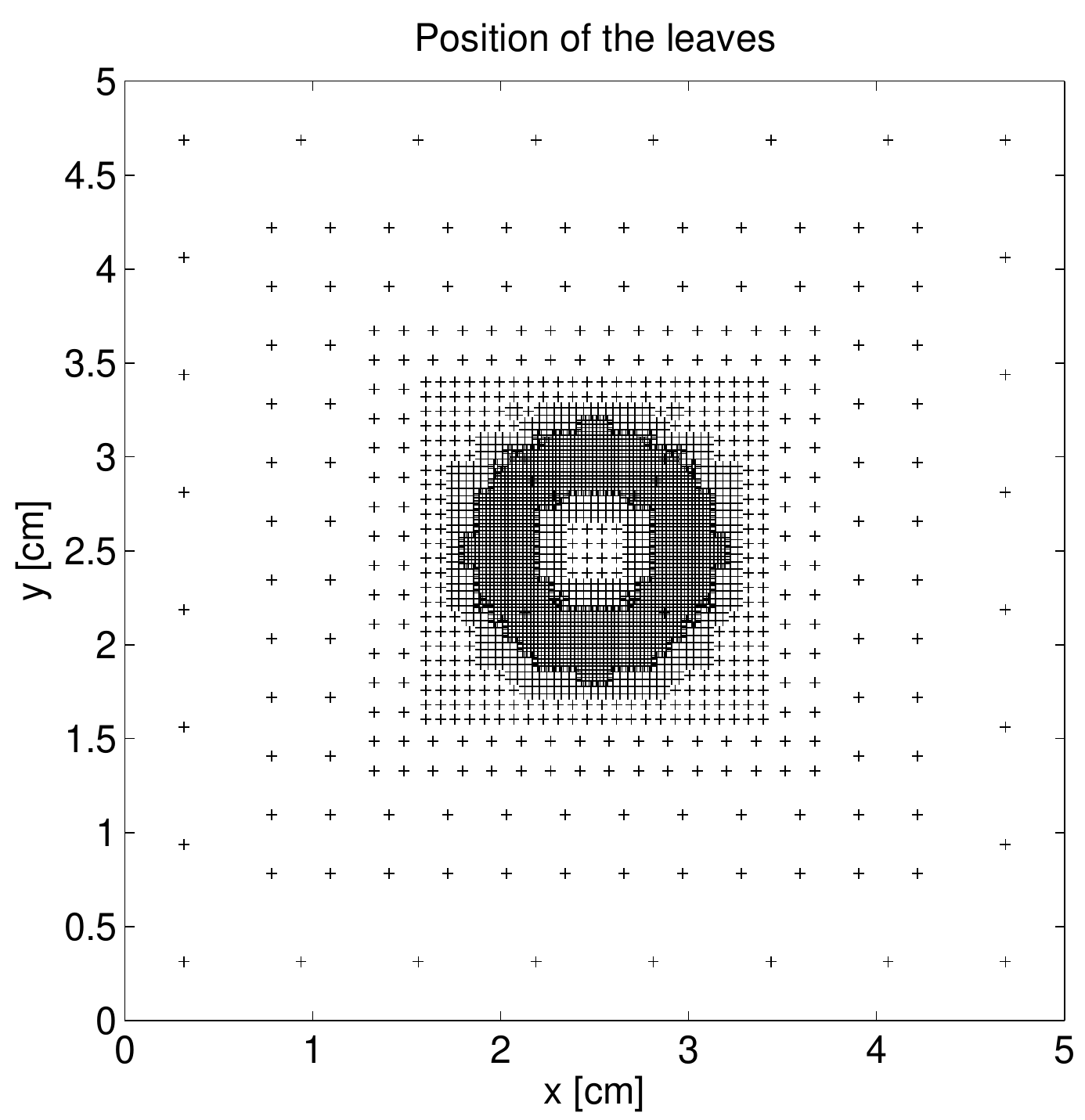}
 \end{tabular}
 \caption{Example~2 (bidomain model, one stimulus):
 Initial condition for the extracellular potential $u_\mathrm{e}$,
  and leaves of the corresponding tree data structure.} \label{fig:cond_ini}
 \end{center}
 \end{figure}

\begin{figure}[t]
\begin{center}
\begin{tabular}{cc}
$t=0.1\,\mathrm{ms}$&$t=0.5\,\mathrm{ms}$\\
\includegraphics[width=0.384\textwidth]{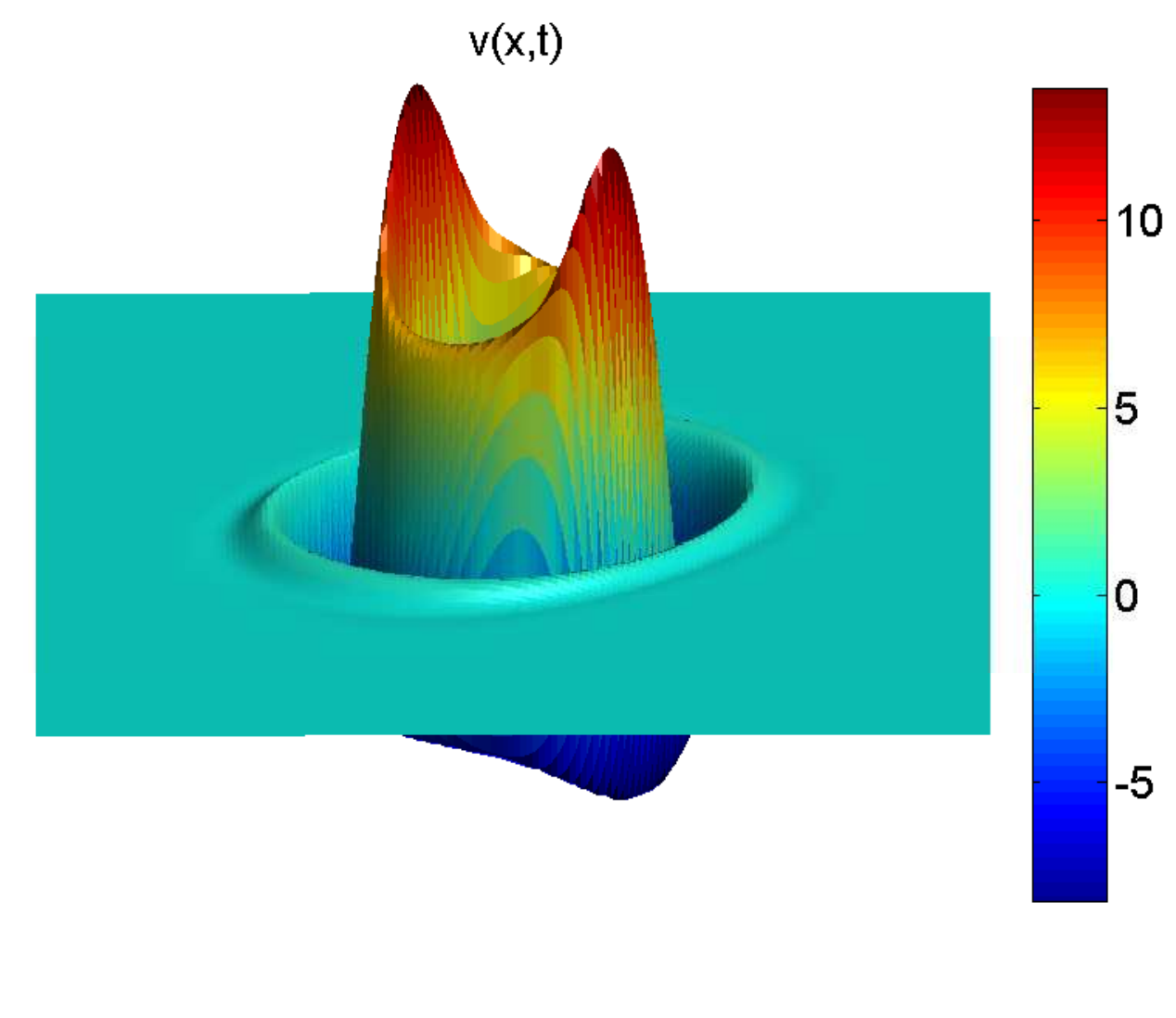}&
\includegraphics[width=0.384\textwidth]{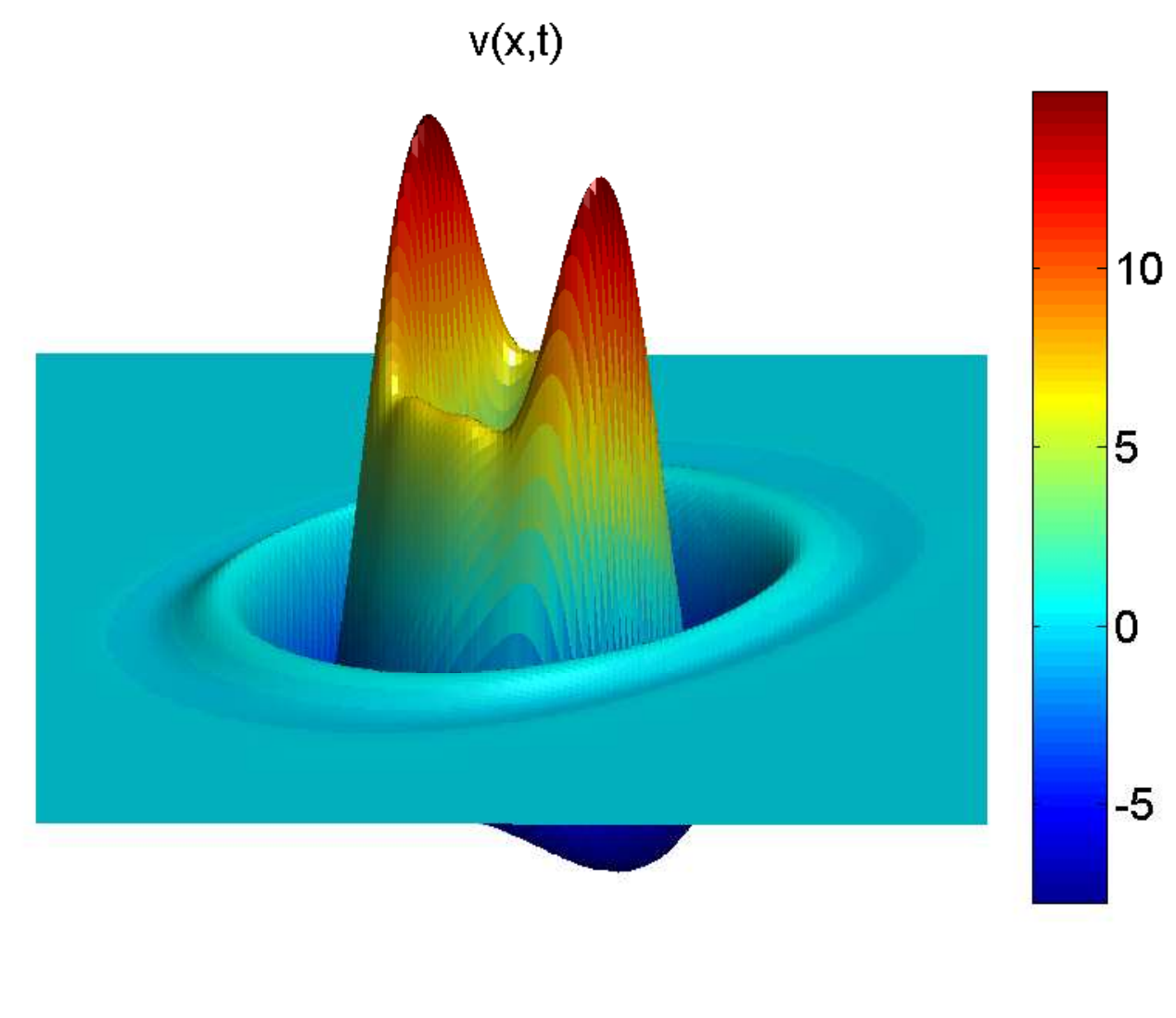}\\
\includegraphics[width=0.384\textwidth]{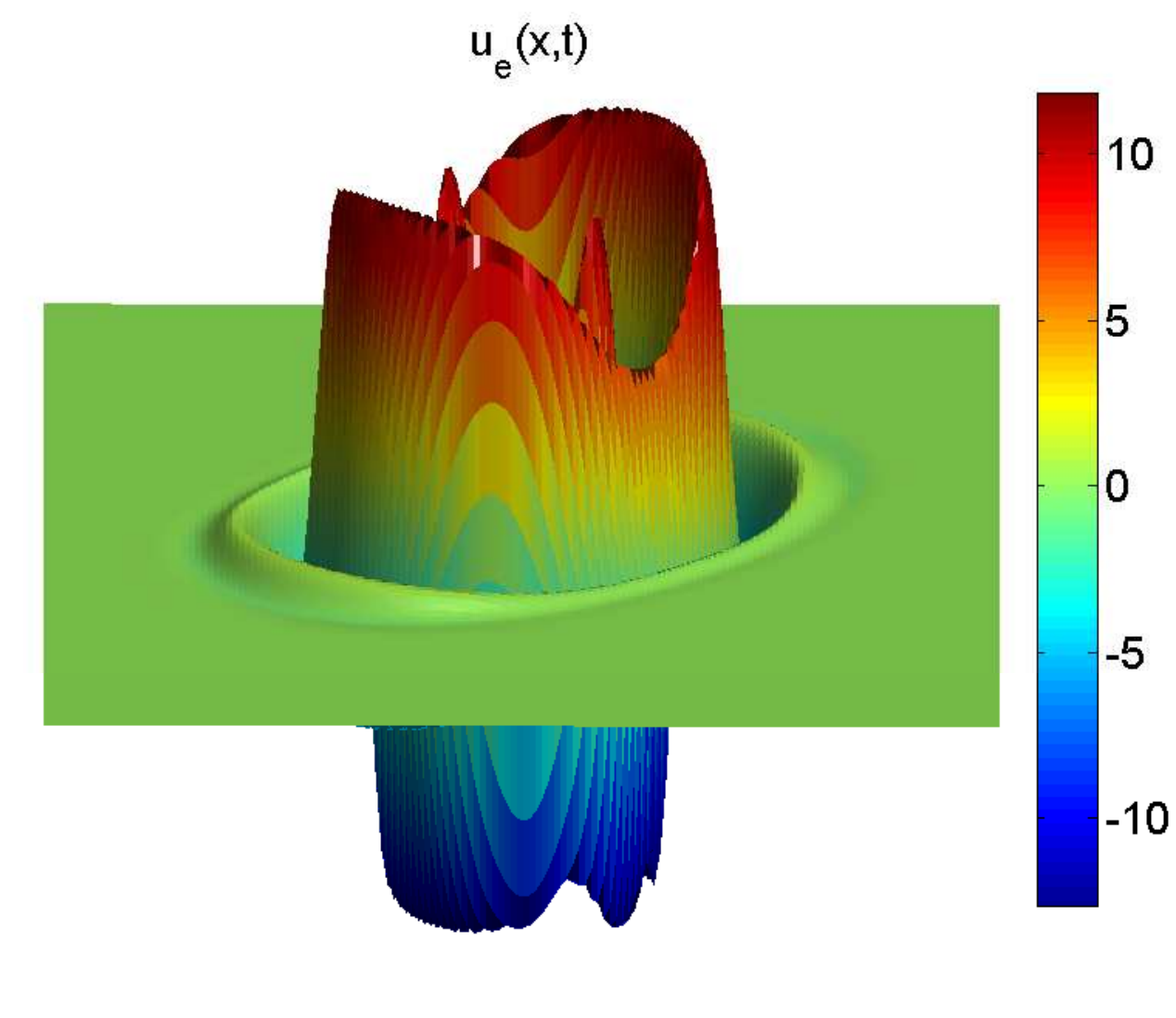}&
\includegraphics[width=0.384\textwidth]{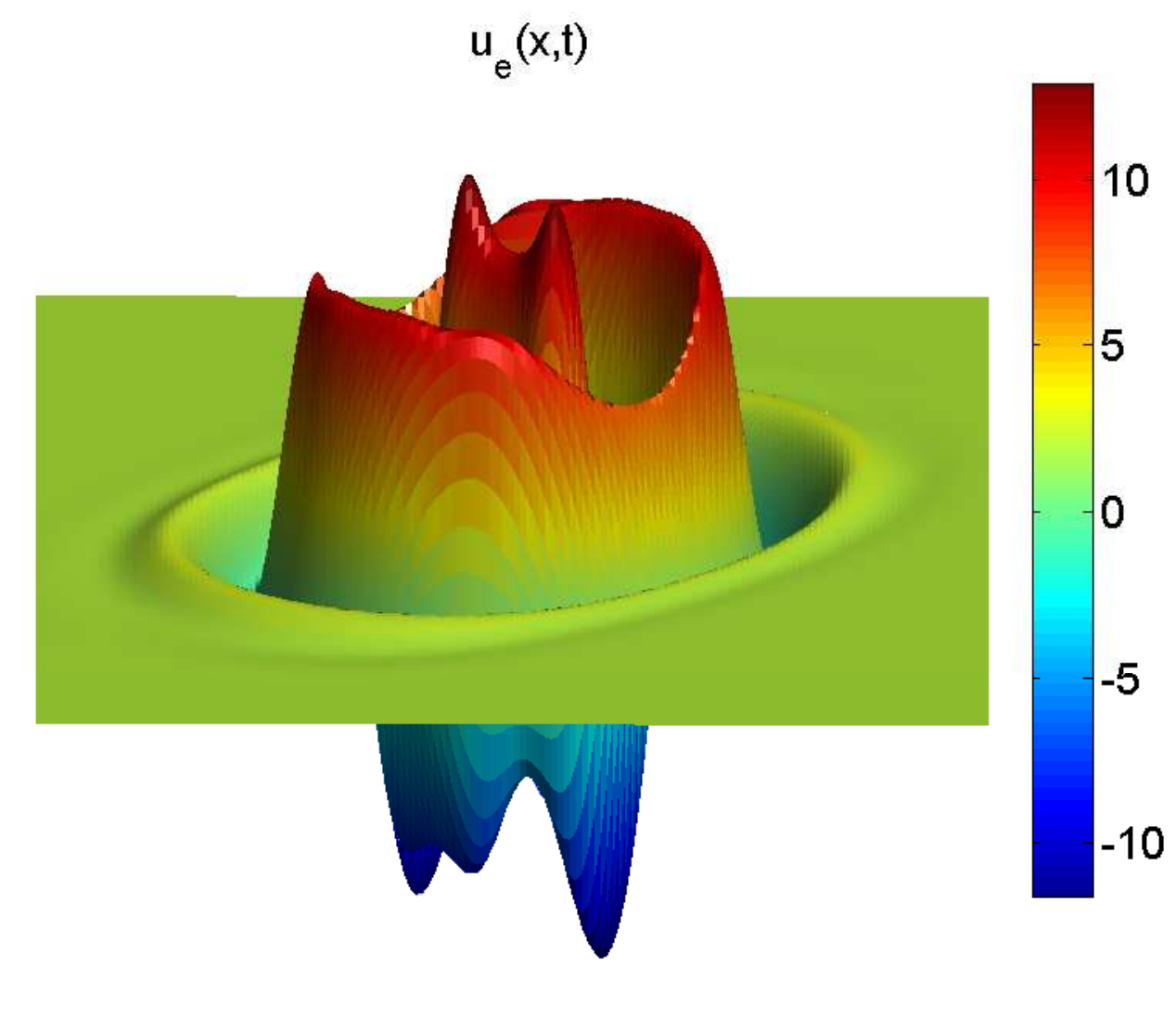}\\
\includegraphics[width=0.322667\textwidth]{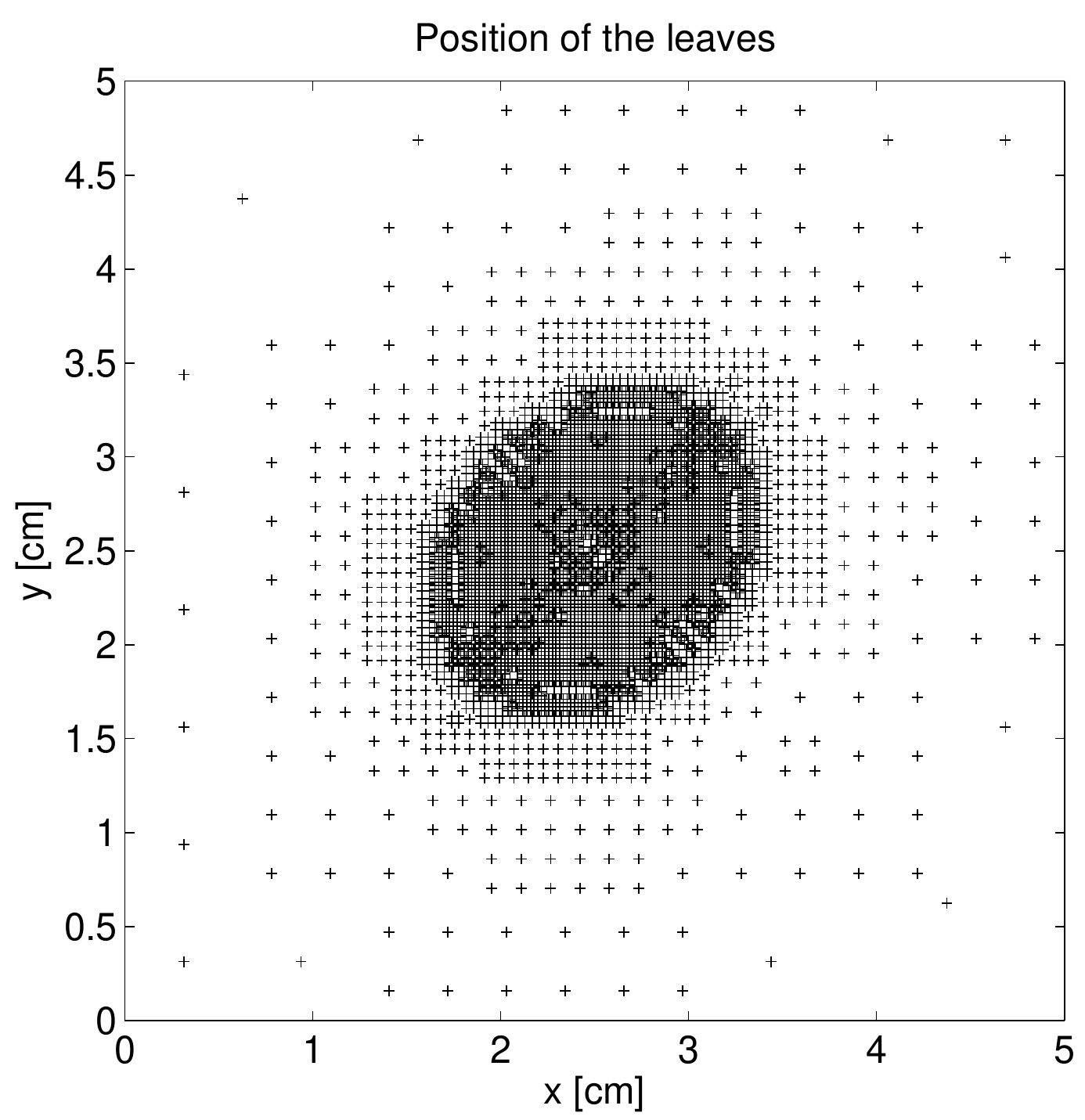}
\quad \qquad$\vphantom{X}$ $\vphantom{X}$&
\includegraphics[width=0.322667\textwidth]{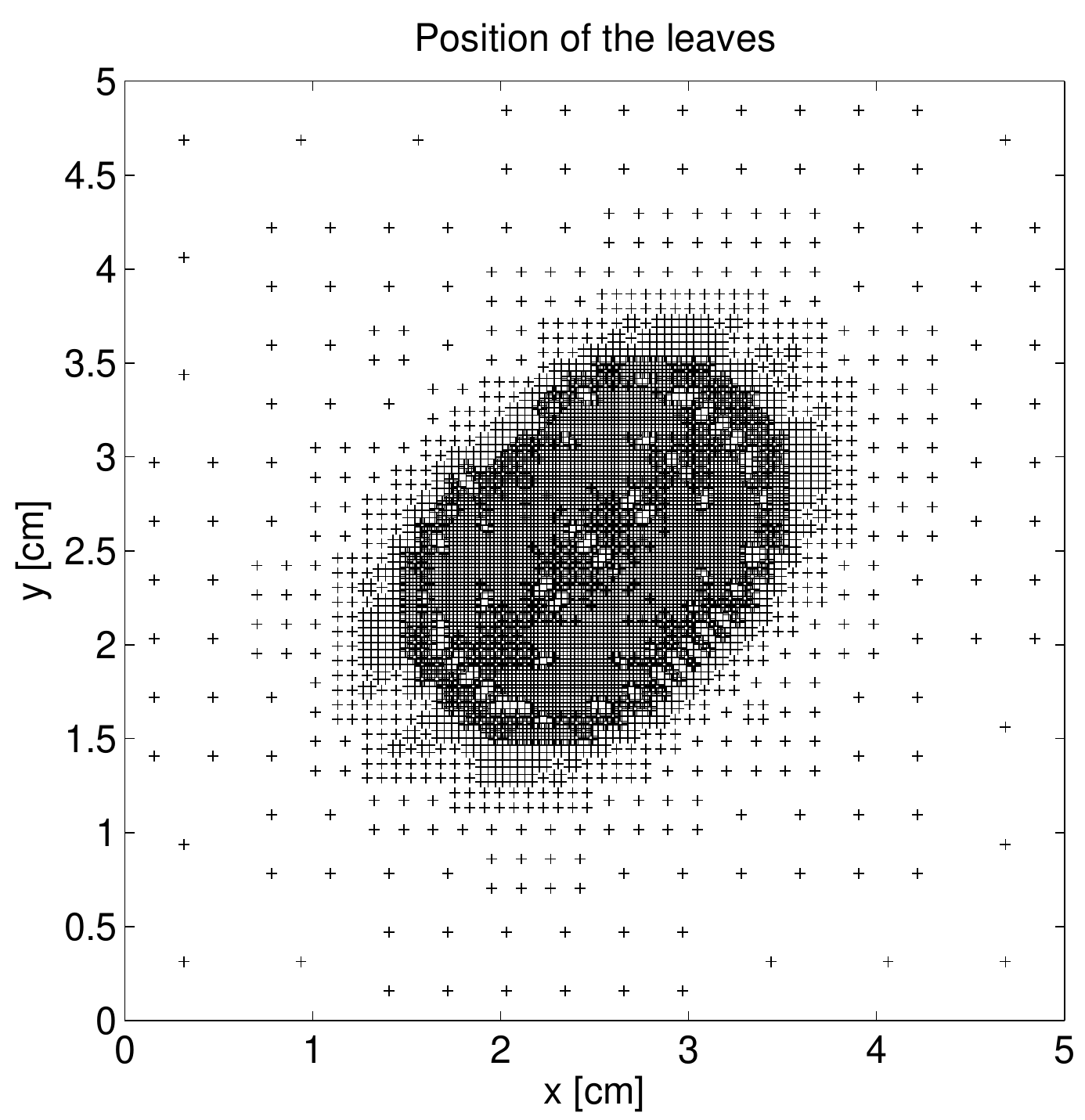}
\quad \qquad$\vphantom{X}$ $\vphantom{X}$
\end{tabular}
\caption{Example~2 (bidomain model, one stimulus):
Numerical solution for transmembrane potential $v$ and
extracellular potential $u_\mathrm{e}$ in $[\mathrm{mV}]$,
and leaves of the corresponding tree data structure at times
$t=0.1\,\mathrm{ms}$ and $t=0.5\,\mathrm{ms}$.} \label{fig:snapshots1}
\end{center}
\end{figure}

\begin{figure}[t]
\begin{center}
\begin{tabular}{cc}
$t=2.0\,\mathrm{ms}$&$t=3.5\,\mathrm{ms}$\\
\includegraphics[width=0.42667\textwidth]{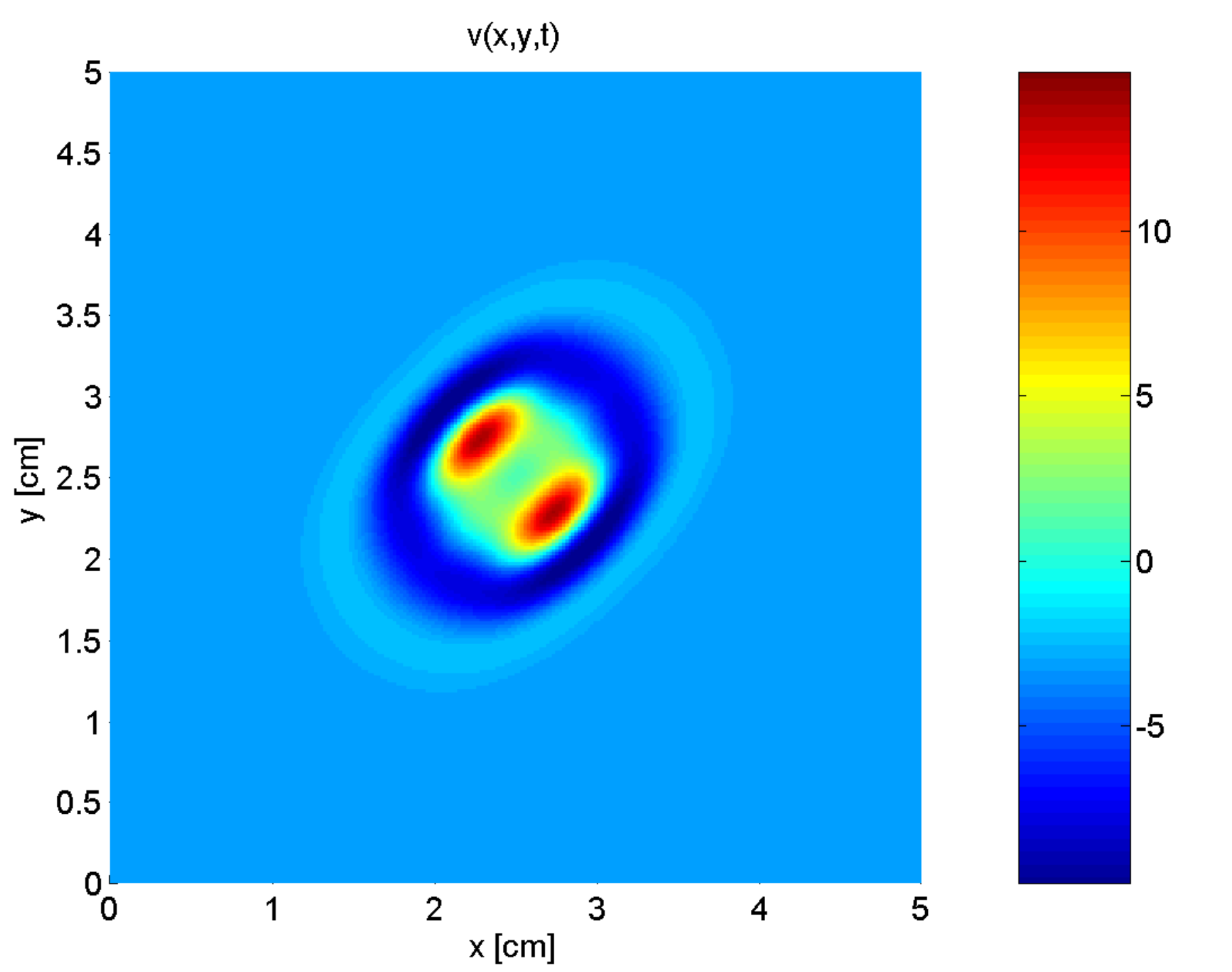}&
\includegraphics[width=0.42667\textwidth]{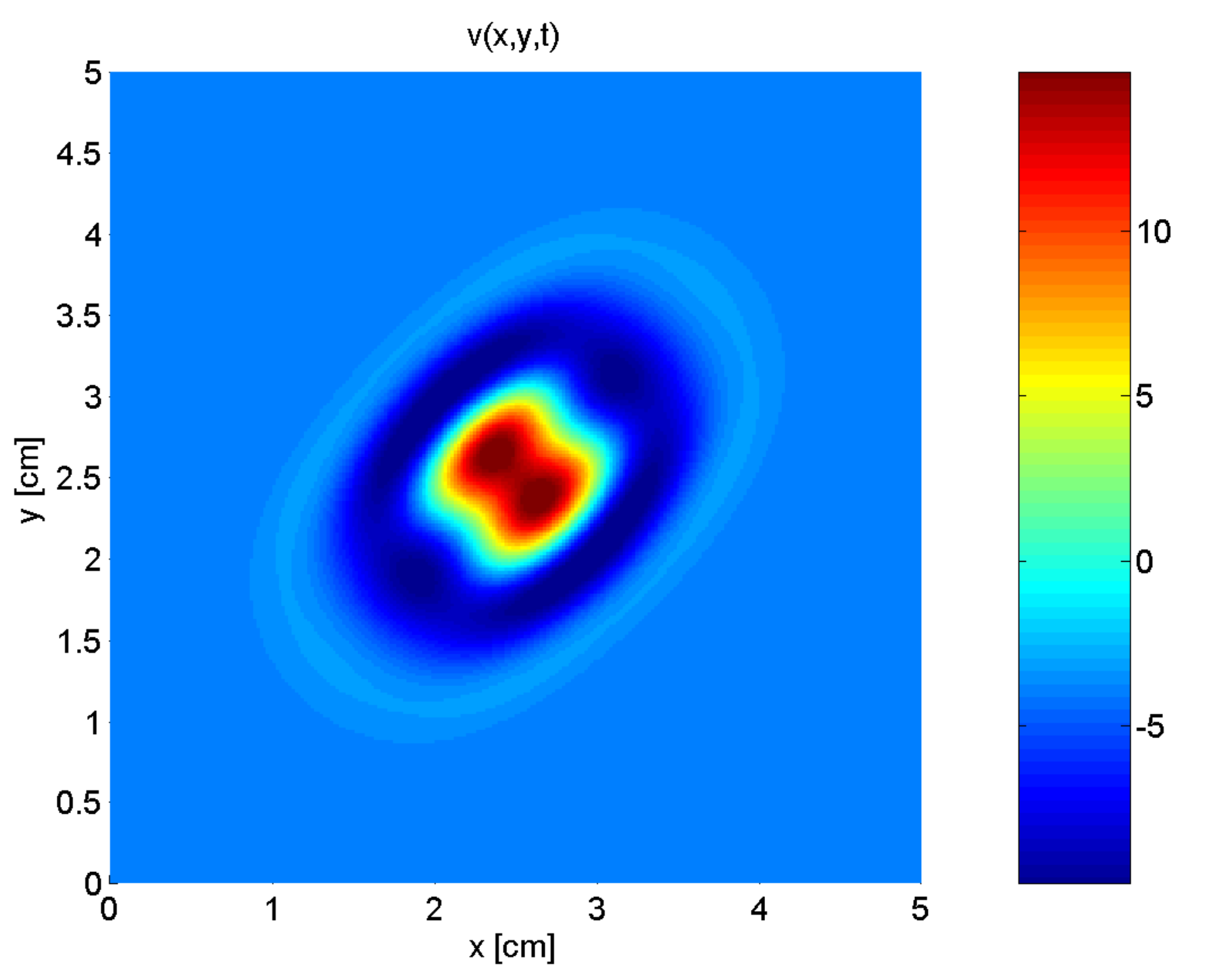}\\
\includegraphics[width=0.42667\textwidth]{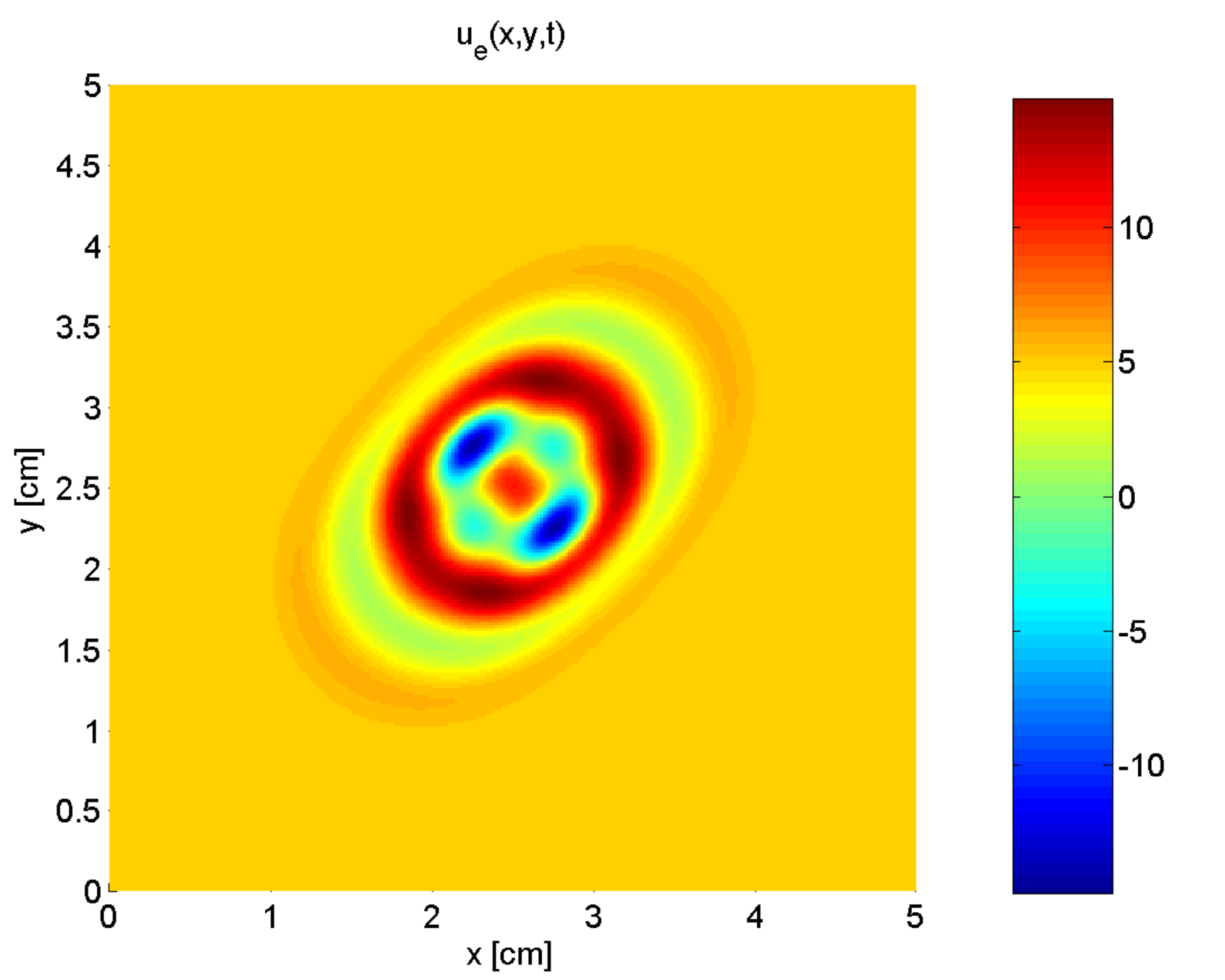}&
\includegraphics[width=0.42667\textwidth]{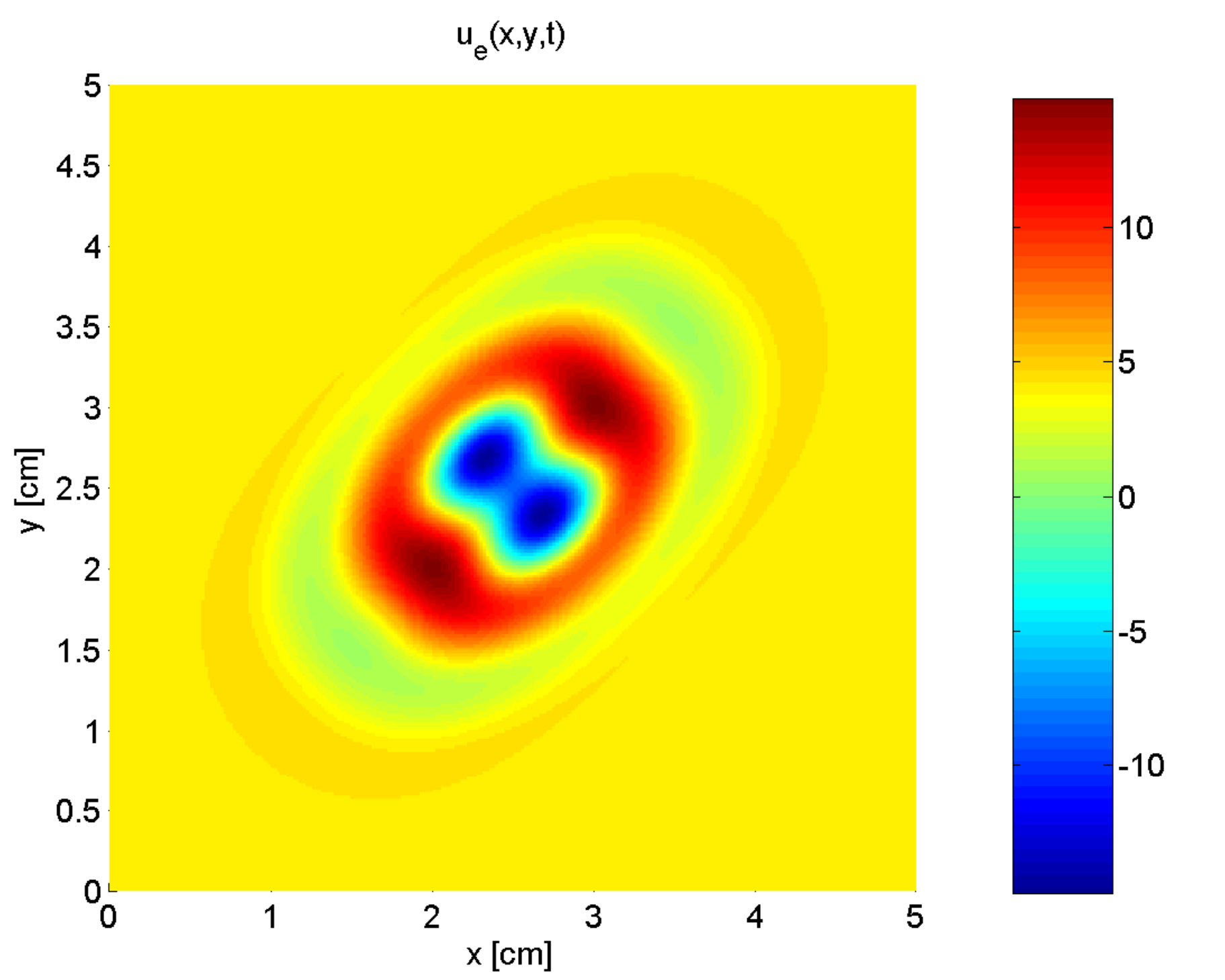}\\
\includegraphics[width=0.322667\textwidth]{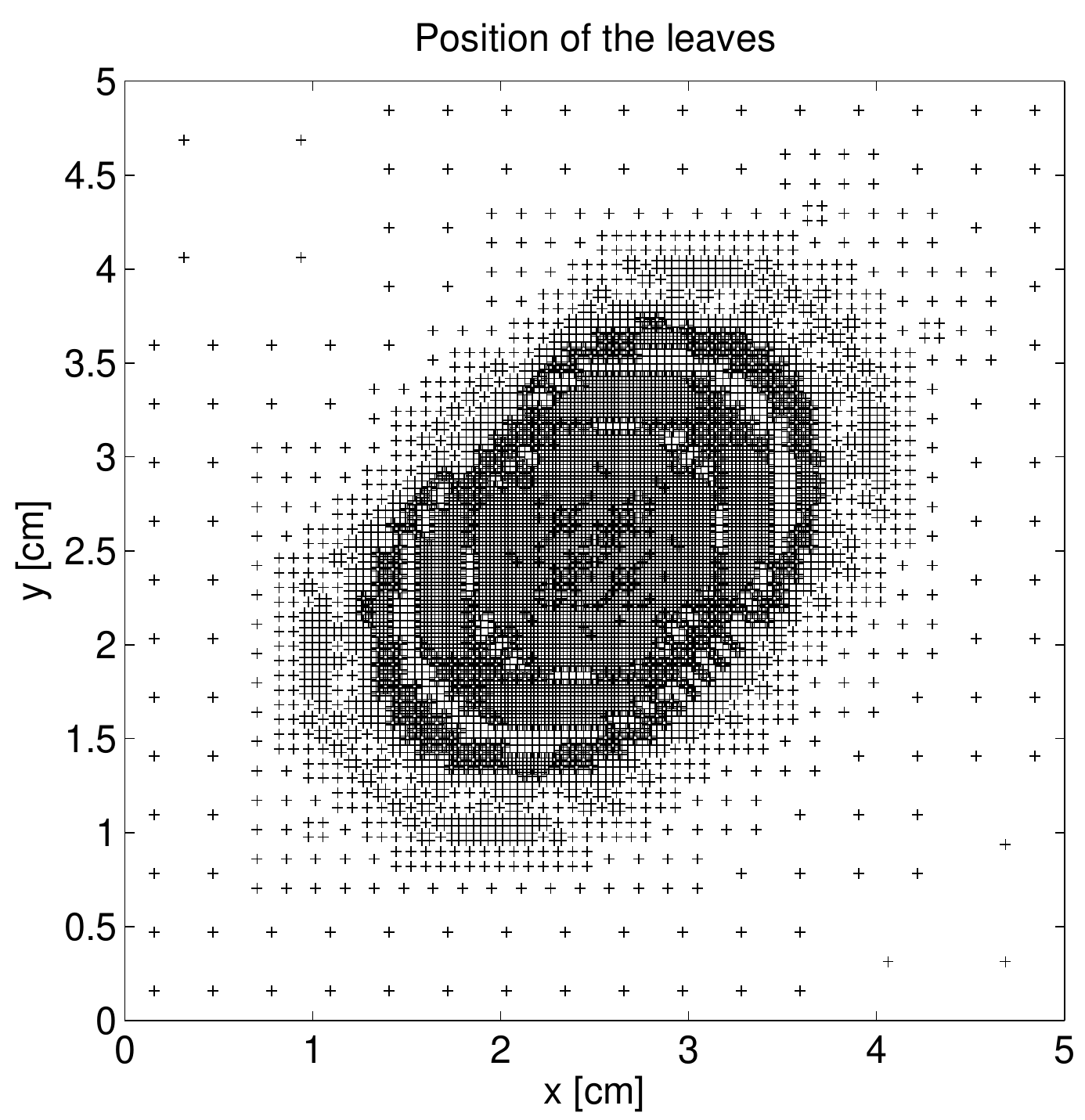}
\quad \qquad$\vphantom{X}$ $\vphantom{X}$ &
\includegraphics[width=0.322667\textwidth]{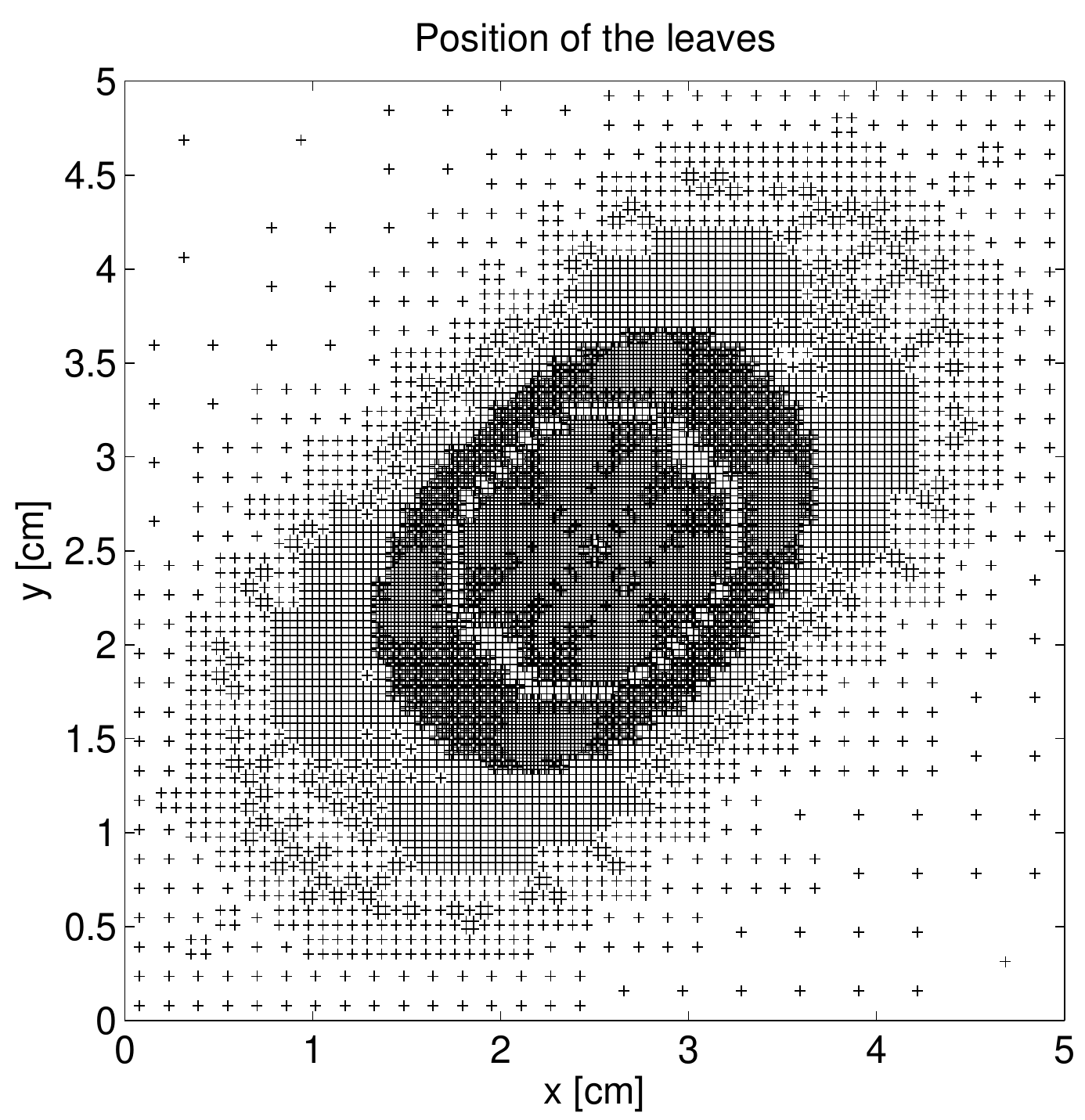}
\quad \qquad$\vphantom{X}$ $\vphantom{X}$
\end{tabular}
\caption{Example~2 (bidomain model, one stimulus):
Numerical solution for transmembrane potential $v$ and
extracellular potential $u_\mathrm{e}$ in $[\mathrm{mV}]$,
 and leaves of the corresponding
tree data structure at times $t=2.0\,\mathrm{ms}$ and
$t=3.5\,\mathrm{ms}$.} \label{fig:snapshots2}
\end{center}
\end{figure}
\begin{table}[t]
\begin{center}
\begin{tabular}{lccccccc}
\hline
Time $[\mathrm{ms}]$& $\mathcal{V}$  & $\eta$&  Potential & $L^1-$error
& $L^2-$error &$L^\infty-$error $\vphantom{\int^X}$   \\
\hline
$t=0.1 $ $\vphantom{\int^X}$
& 13.74 & 19.39& $v$  &$3.68\times10^{-4}$
&$8.79\times10^{-5}$&$6.51\times10^{-4}$ \\
&      &      & $u_\mathrm{e}$ &$2.01\times10^{-4}$
&$6.54\times10^{-5}$&$5.22\times10^{-4}$  \\
$t=0.5$
& 21.40 & 17.63 & $v$& $4.06\times10^{-4}$
&$9.26\times10^{-5}$&$6.83\times10^{-4}$ \\
&      &       & $u_\mathrm{e}$&$2.79\times10^{-4}$
&$8.72\times10^{-5}$&$5.49\times10^{-4}$ \\
$t=2.0$
& 25.23 & 17.74& $v$ & $4.37\times10^{-4}$
&$1.25\times10^{-4}$&$6.88\times10^{-4}$ \\
&      &     & $u_\mathrm{e}$ & $3.48\times10^{-4}$
&$9.44\times10^{-5}$&$6.11\times10^{-4}$ \\
$t=5.0$
& 26.09 & 16.35 & $v$& $5.29\times10^{-4}$
&$1.94\times10^{-4}$&$7.20\times10^{-4}$ \\
&       &     & $u_\mathrm{e}$& $4.15\times10^{-4}$
&$1.06\times10^{-4}$&$6.32\times10^{-4}$ \\
\hline
\end{tabular}
\end{center}

\vspace*{2mm}

\caption{Example~2 (bidomain model, one stimulus):  Corresponding simulated time,
CPU ratio~$\mathcal{V}$, compression rate~$\eta$  and normalized errors.}
\label{table:ex2}
\end{table}

\begin{figure}[t]
\begin{center}
\includegraphics[width=0.32\textwidth,height=0.4\textwidth]{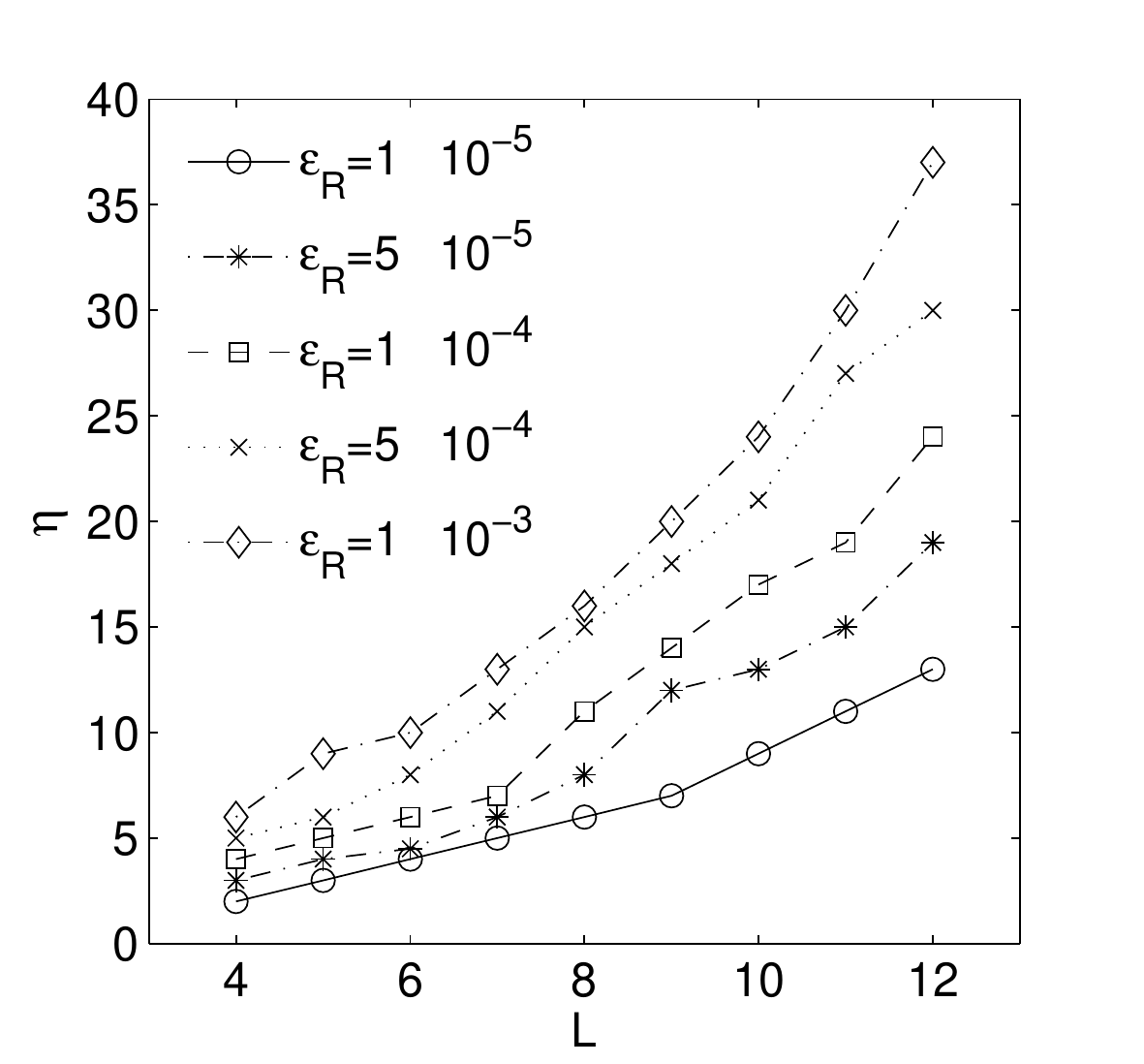}
\includegraphics[width=0.32\textwidth,height=0.4\textwidth]{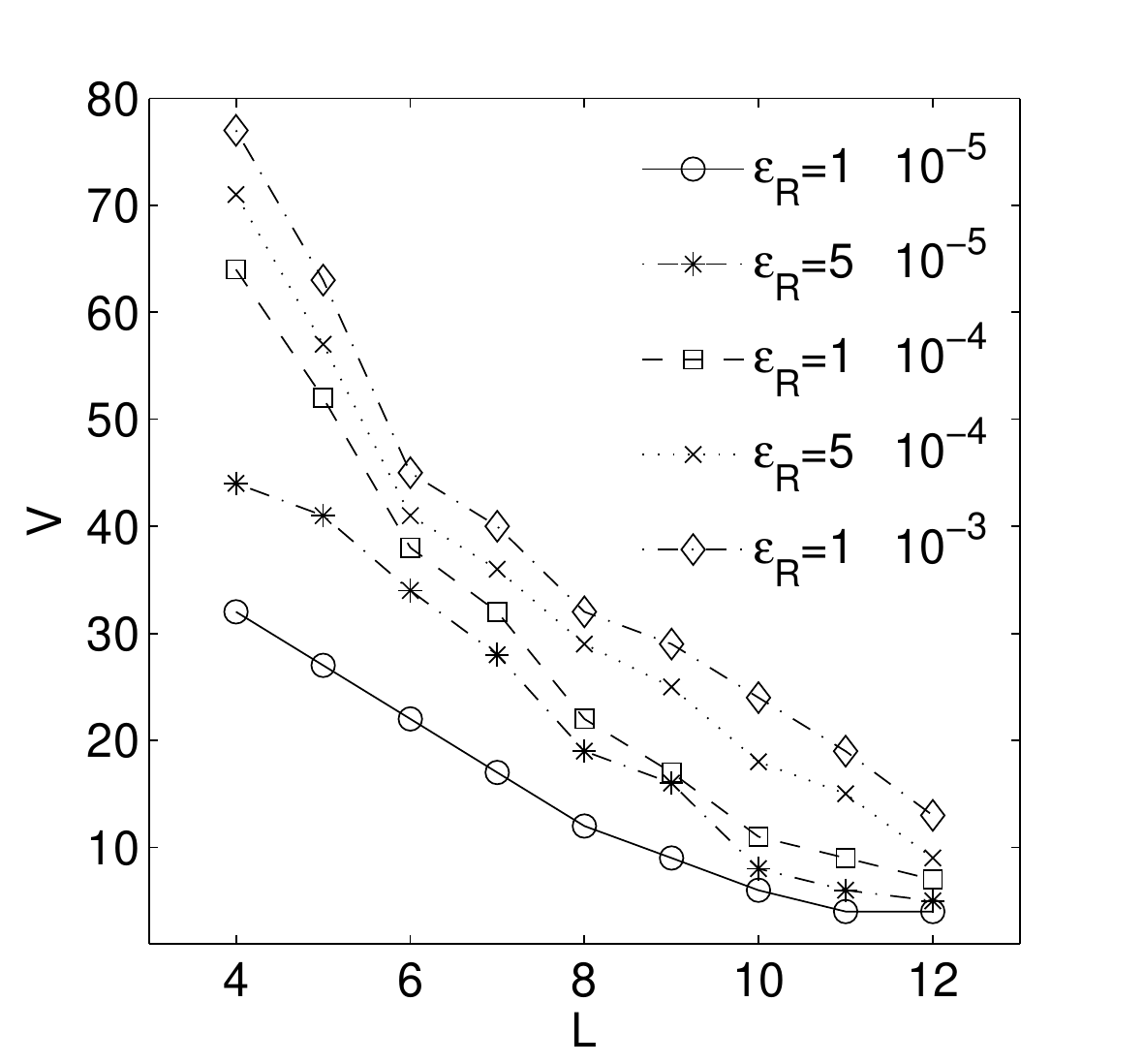}
\includegraphics[width=0.32\textwidth,height=0.4\textwidth]{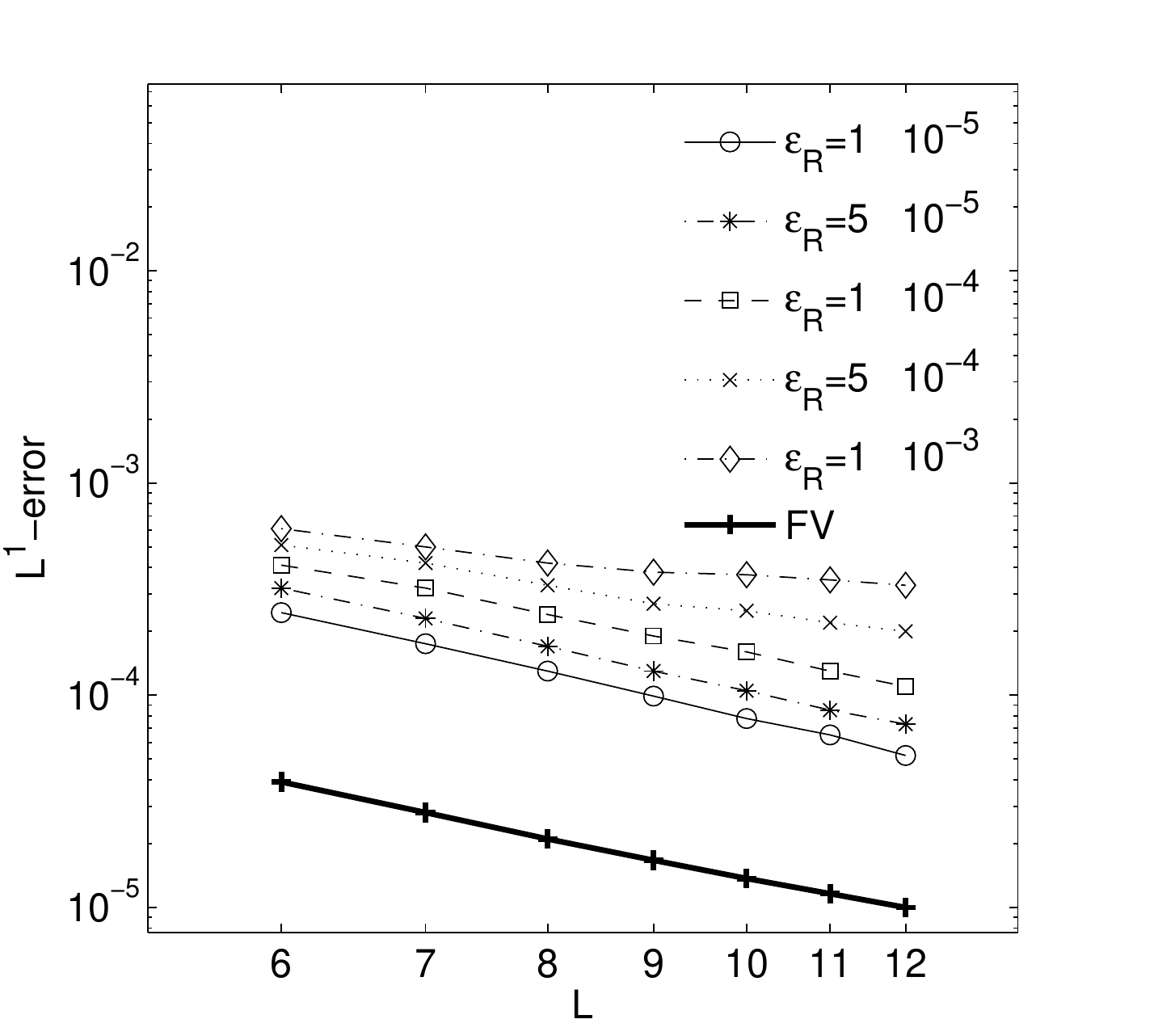}
\end{center}
\caption{Example~2 (bidomain model): data compression
rate~$\eta$ (left), speed-up factor~$\mathcal{V}$ (middle) and $L^1$-errors for
different scales~$L$ and values of~$\eps_{\mathrm{R}}$ (right). The simulated time is
$t=2.0\,\mathrm{ms}$.} \label{fig:factorC}
\end{figure}
\begin{figure}[t]
\begin{center}
\begin{tabular}{cc}
\includegraphics[width=0.42667\textwidth,height=0.3733\textwidth]{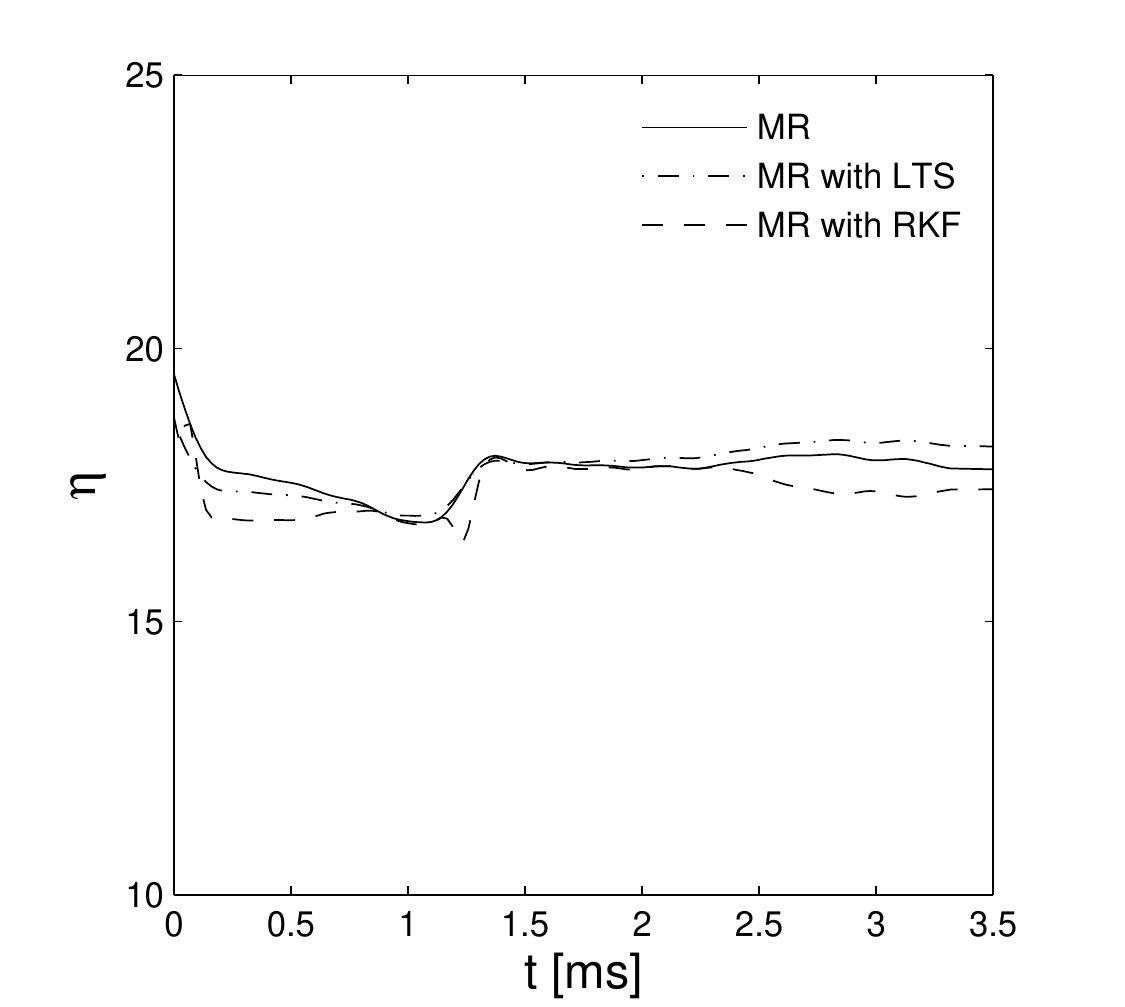}&
\includegraphics[width=0.42667\textwidth,height=0.3733\textwidth]{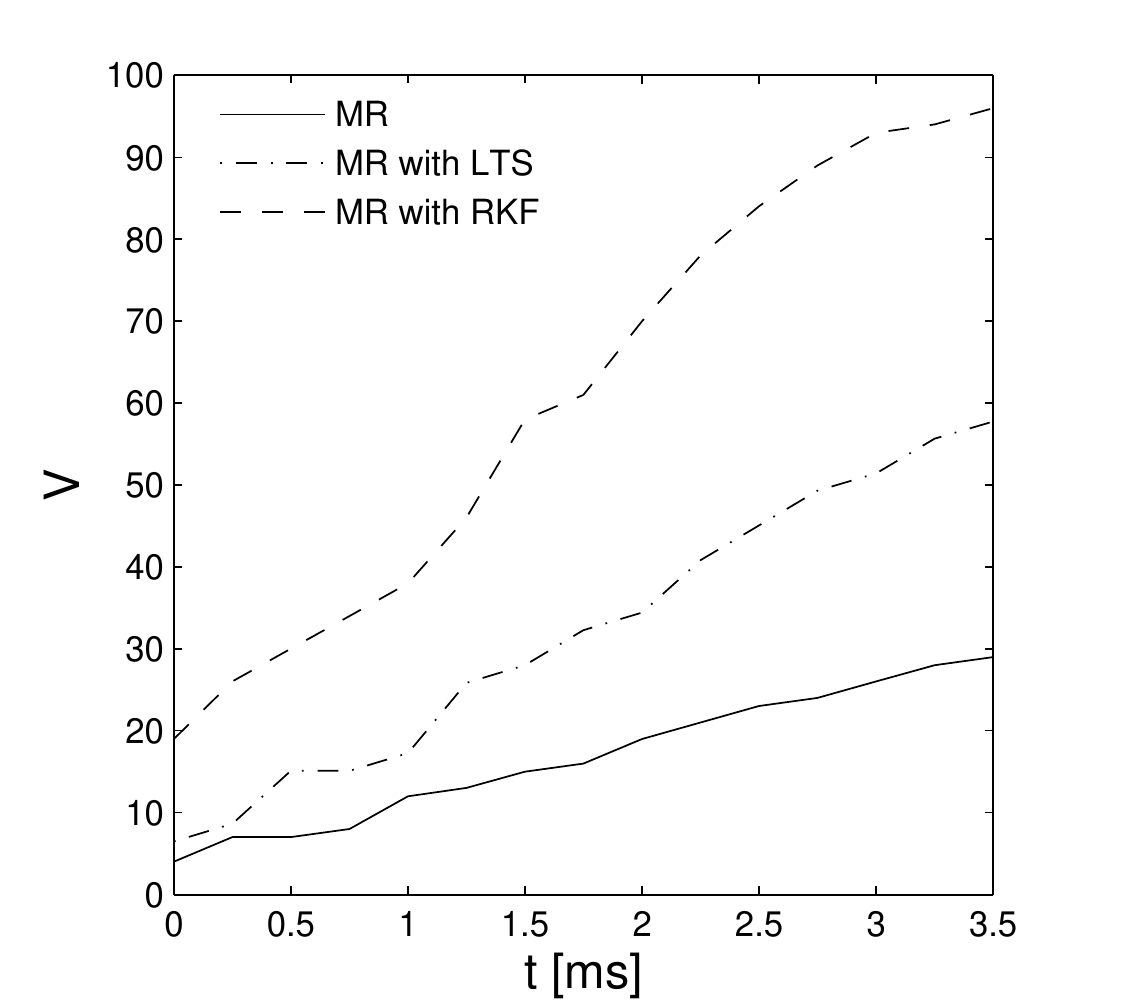}\\
\includegraphics[width=0.42667\textwidth,height=0.3733\textwidth]{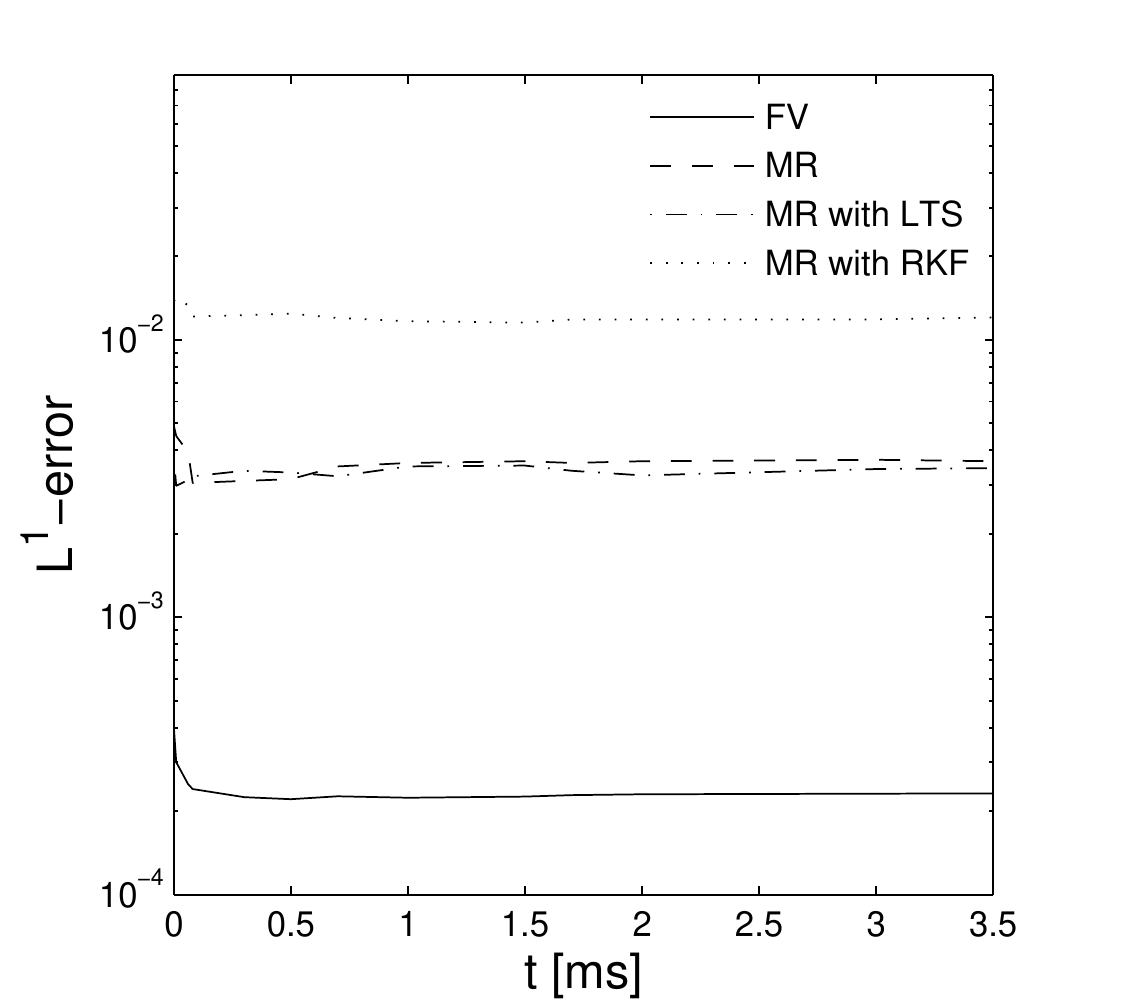}&
\includegraphics[width=0.42667\textwidth,height=0.3733\textwidth]{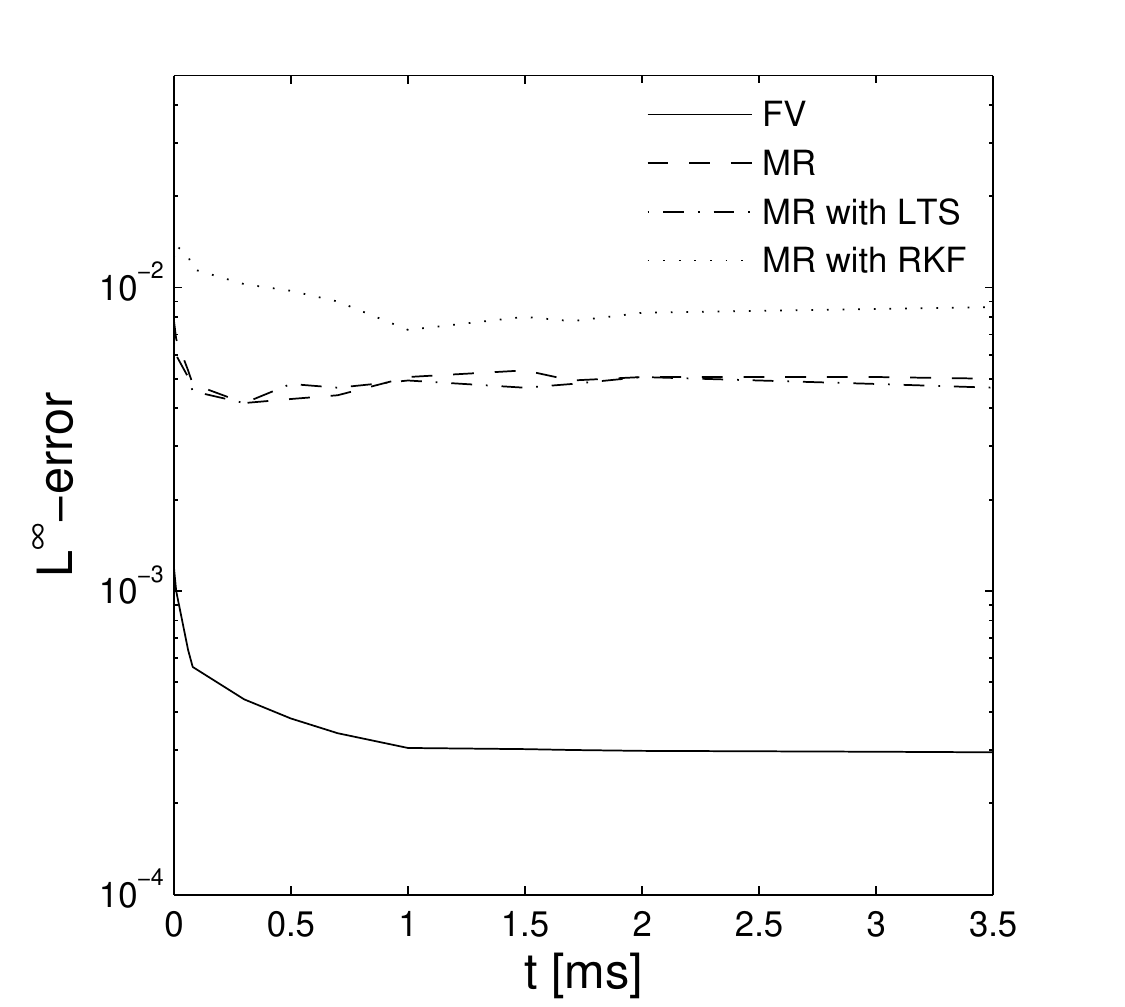}
\end{tabular}
\caption{Example~2 (bidomain model, one stimulus):  Time evolution for data
compression rate $\eta$, speed-up rate $\mathcal{V}$, and normalized errors for different
methods: MR scheme with  global time step, MR with locally
varying time stepping and MR with RKF time stepping.}
\label{fig:varios}
\end{center}
\end{figure}
 \begin{figure}[t]
 \begin{center}
 \begin{tabular}{cc}
 $t=0.1\,\mathrm{ms}$&$t=0.5\,\mathrm{ms}$\\
 \includegraphics[width=0.42667\textwidth]{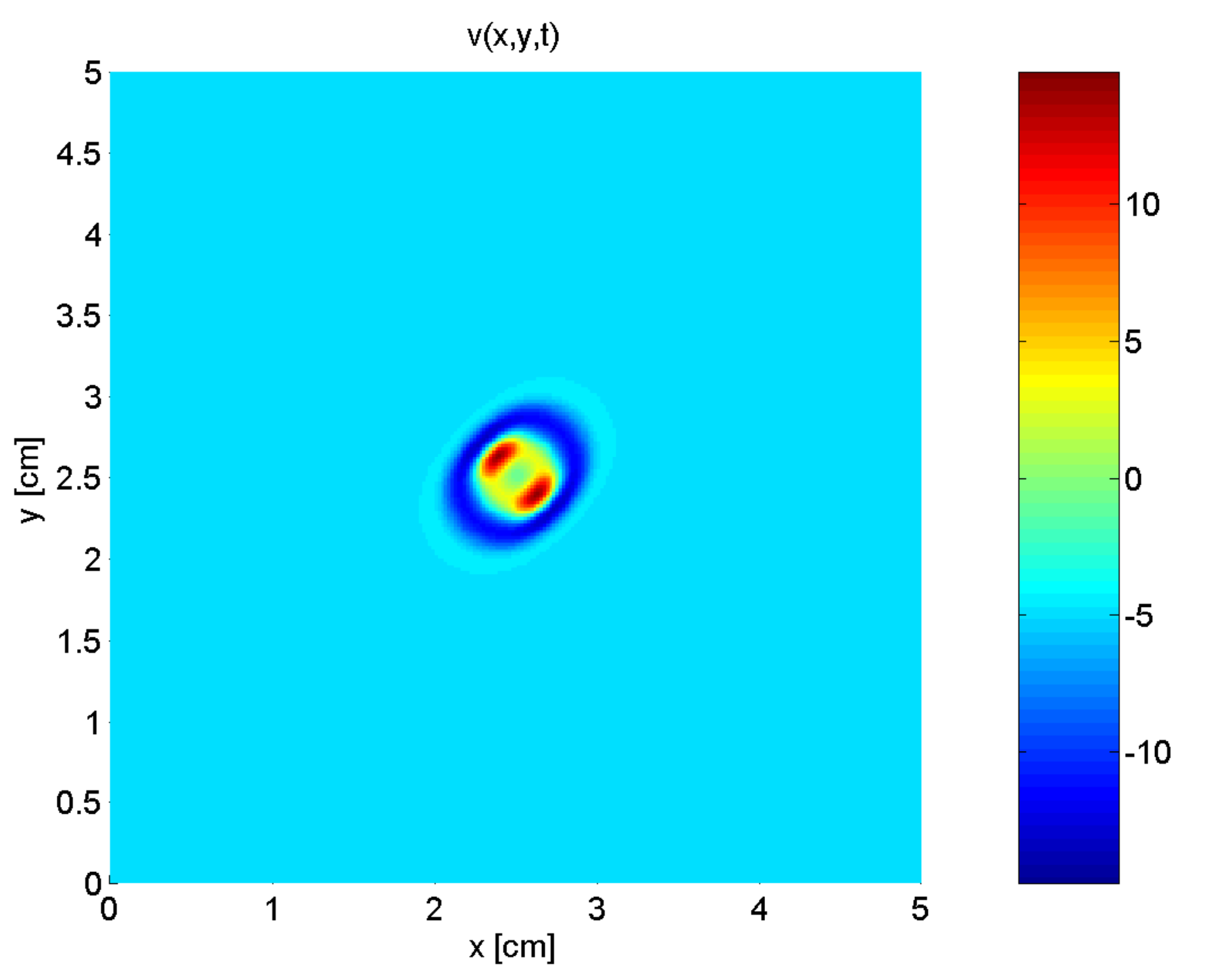}&
 \includegraphics[width=0.42667\textwidth]{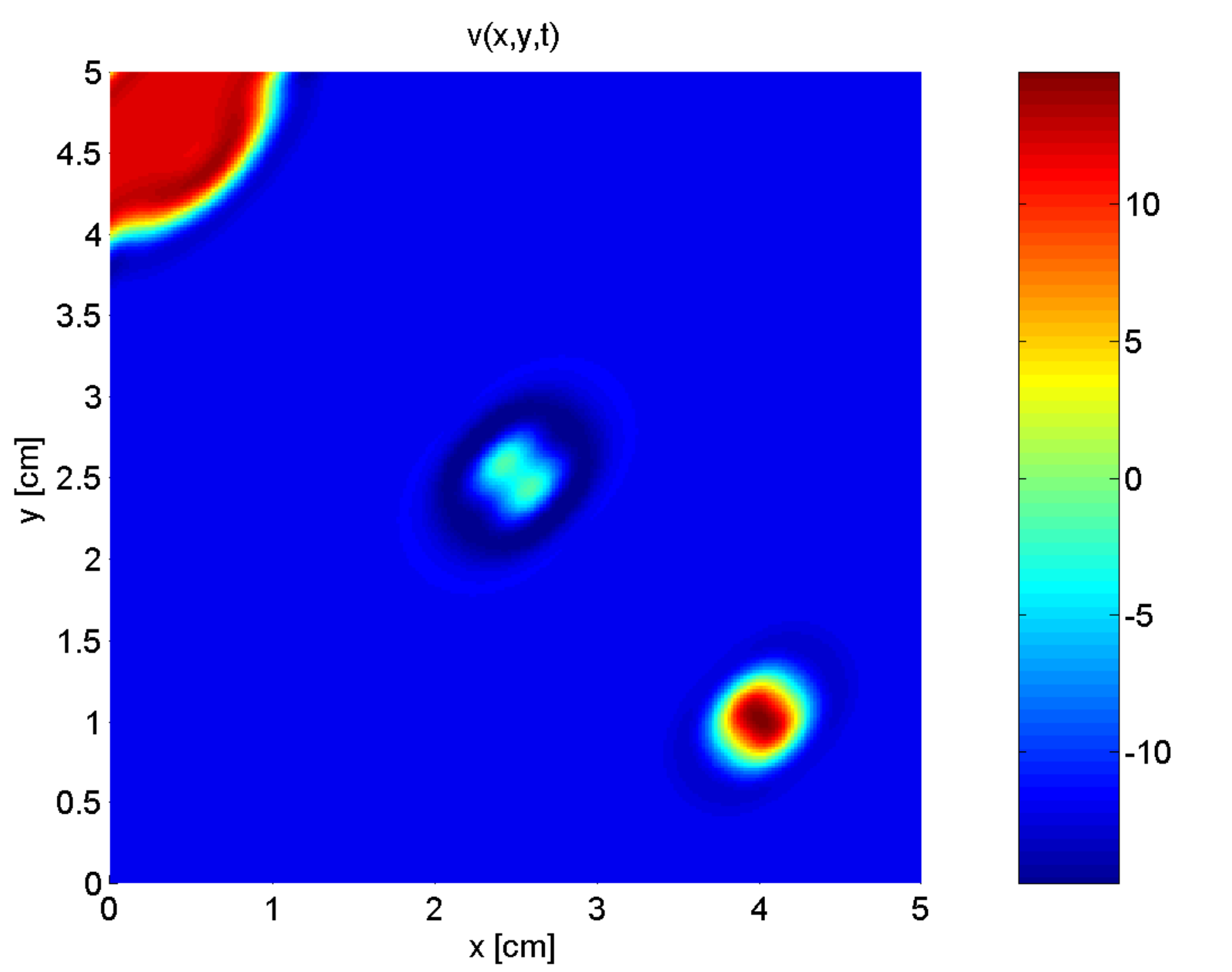}\\
 \includegraphics[width=0.42667\textwidth]{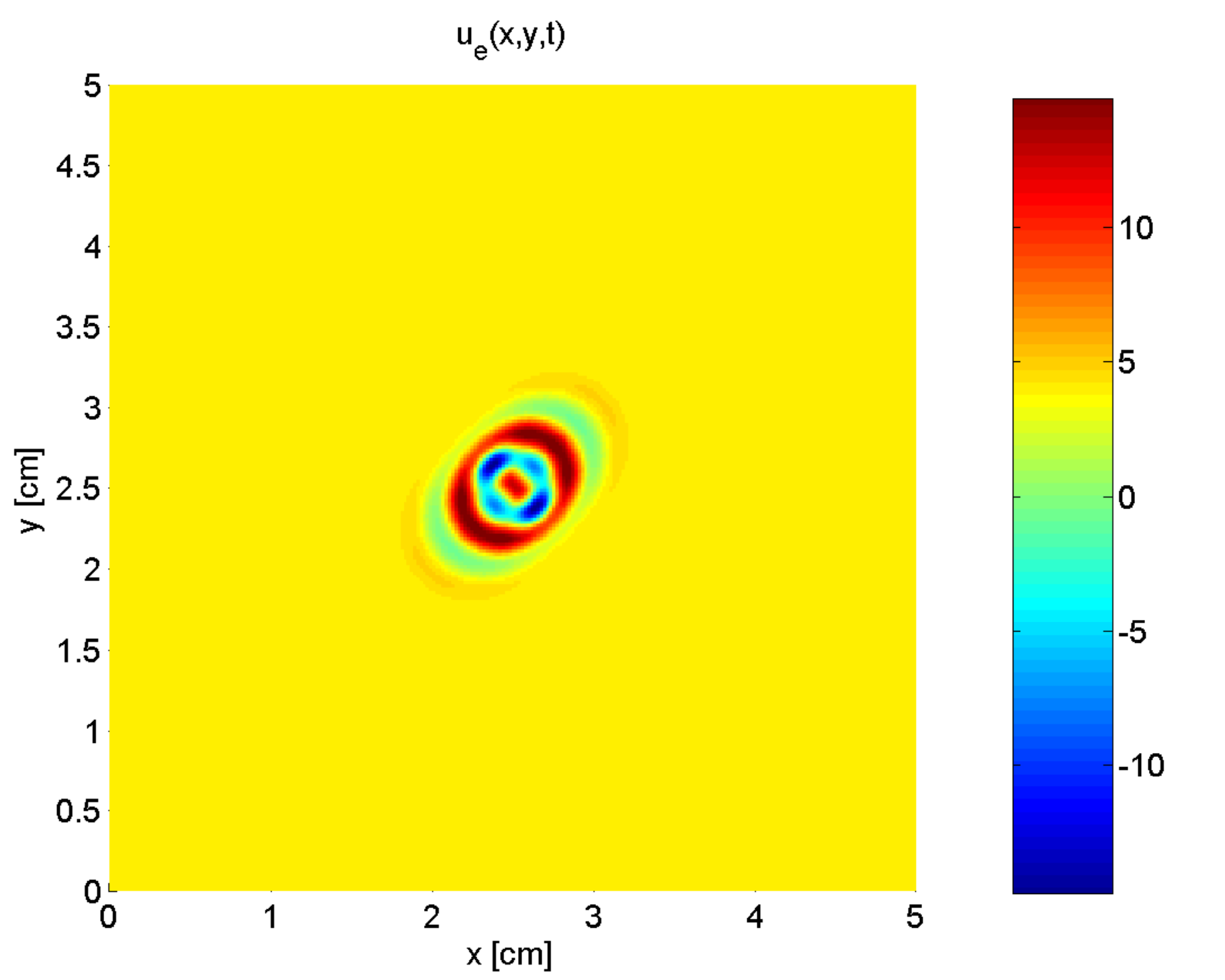}&
 \includegraphics[width=0.42667\textwidth]{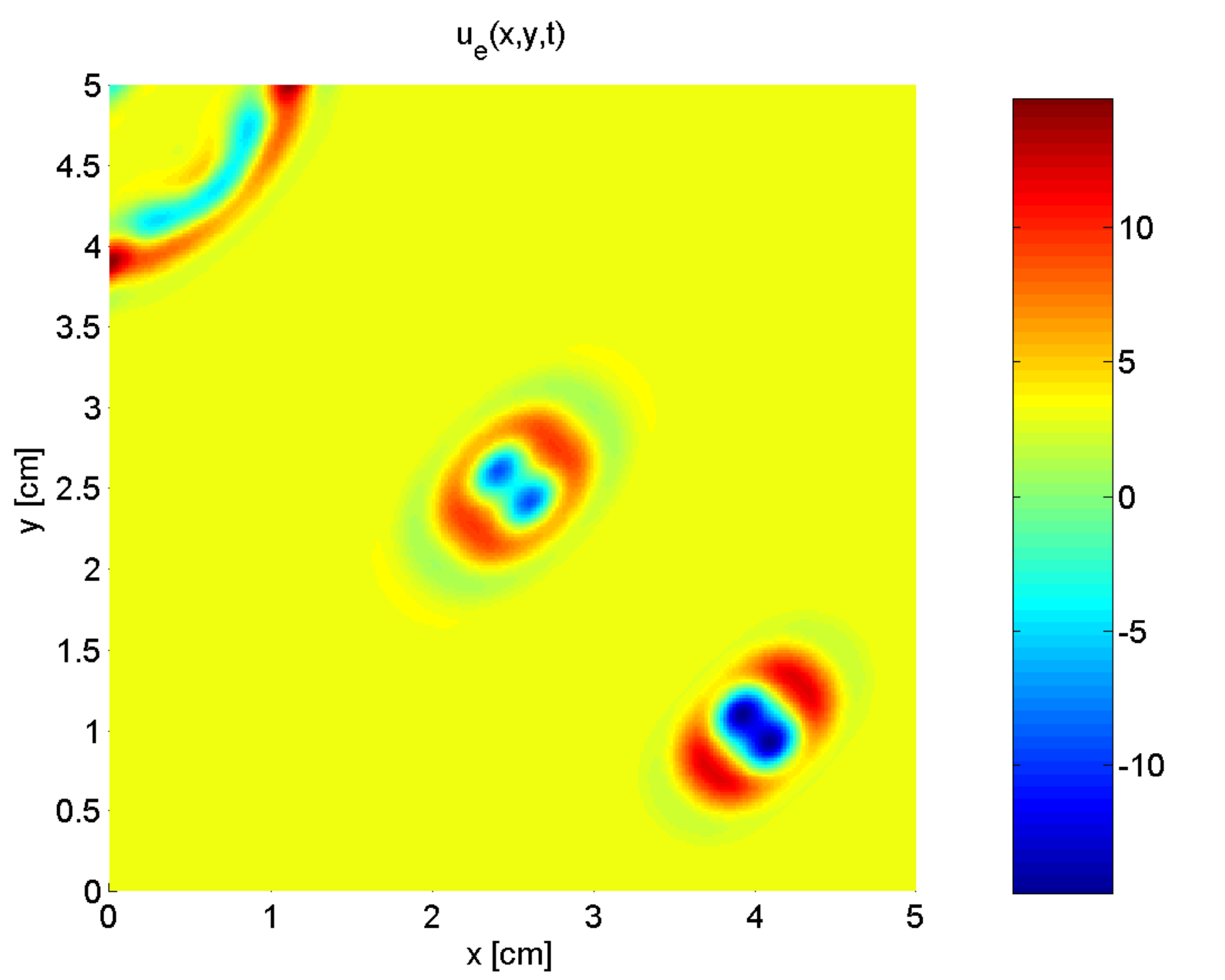}\\
 \includegraphics[width=0.322667\textwidth]{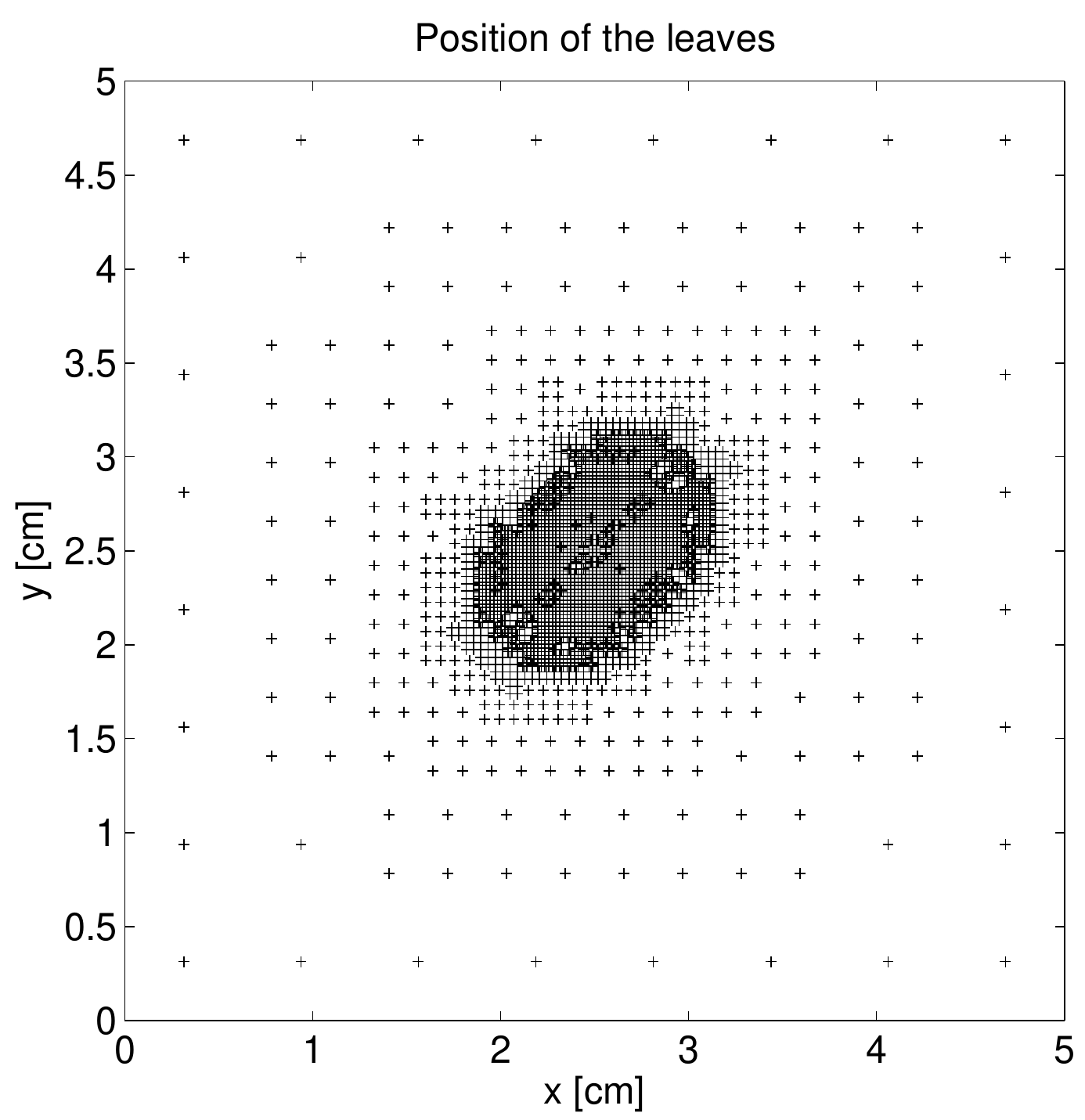}
\quad \qquad$\vphantom{X}$ $\vphantom{X}$ &
 \includegraphics[width=0.322667\textwidth]{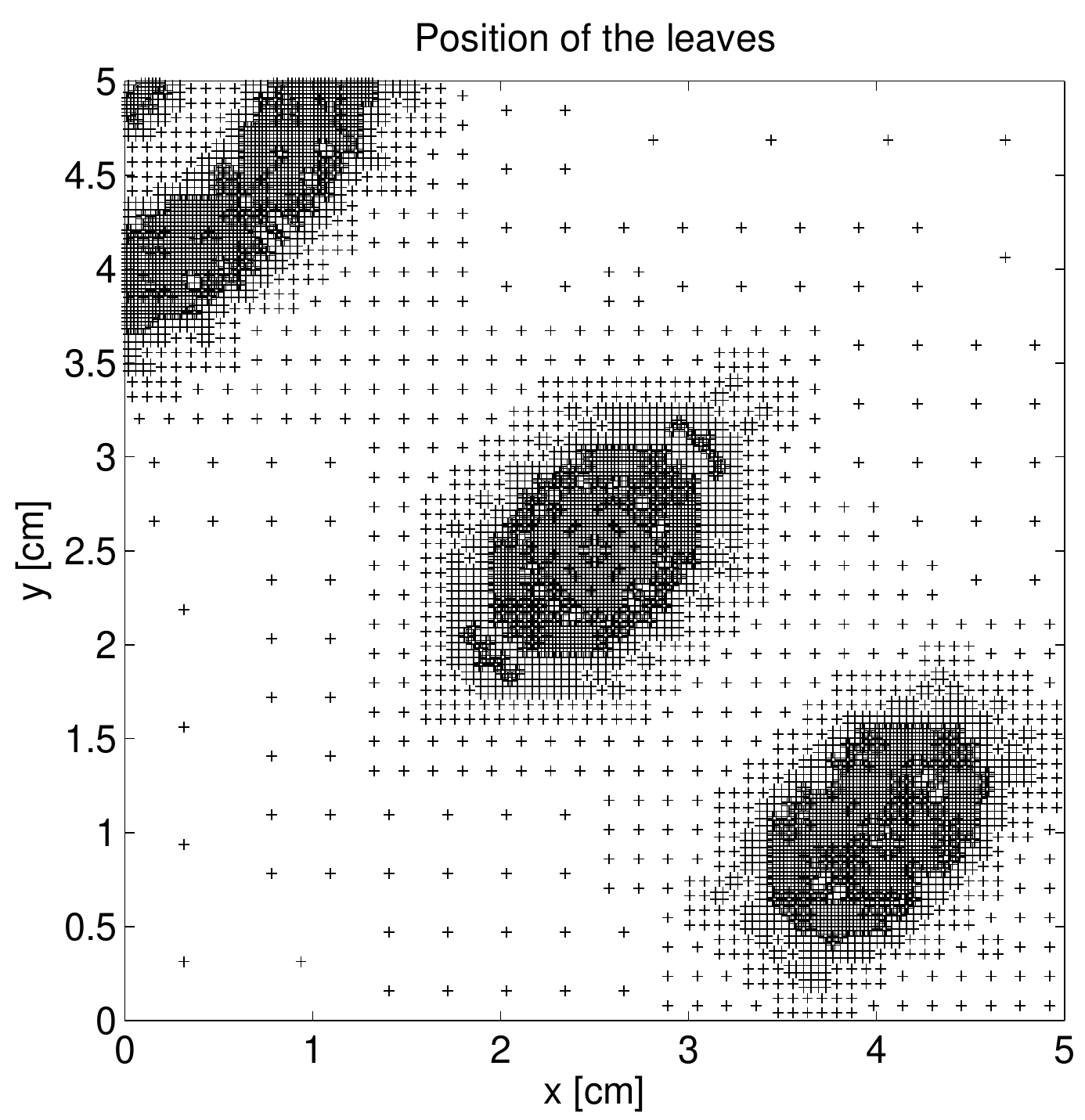}
\quad \qquad$\vphantom{X}$ $\vphantom{X}$
 \end{tabular}
 \caption{Example~3 (bidomain model, three stimuli):
 Numerical solution for transmembrane potential $v$ and
 extracellular potential $u_\mathrm{e}$ in $[\mathrm{mV}]$,
and leaves of the corresponding
 tree data structure at times $t=0.1\,\mathrm{ms}$ and
$t=0.5\,\mathrm{ms}$.} \label{fig:snapshots3}
 \end{center}
 \end{figure}

 \begin{figure}[t]
 \begin{center}
 \begin{tabular}{cc}
 $t=2.0\,\mathrm{ms}$&$t=5.0\,\mathrm{ms}$\\
 \includegraphics[width=0.42667\textwidth]{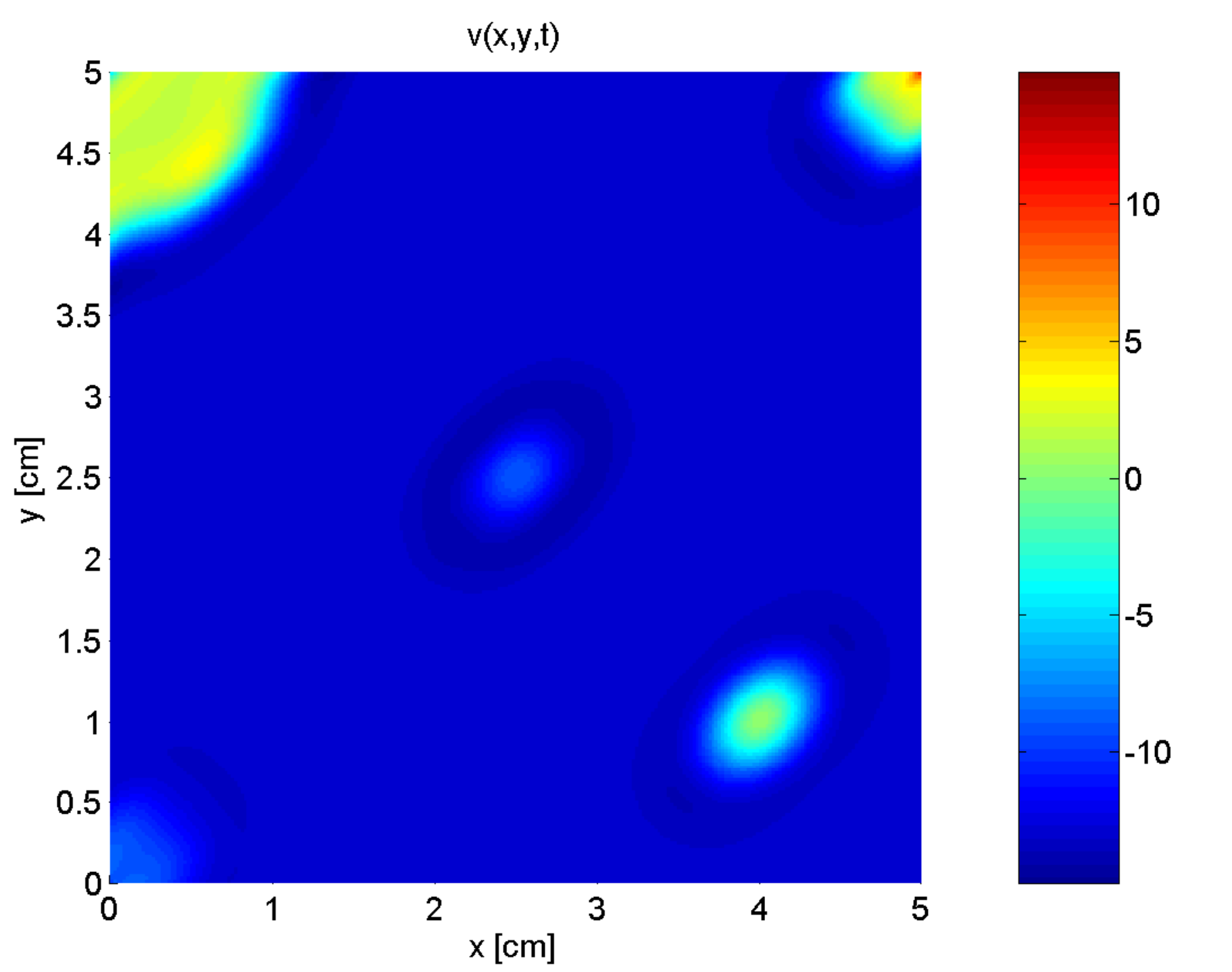}&
 \includegraphics[width=0.42667\textwidth]{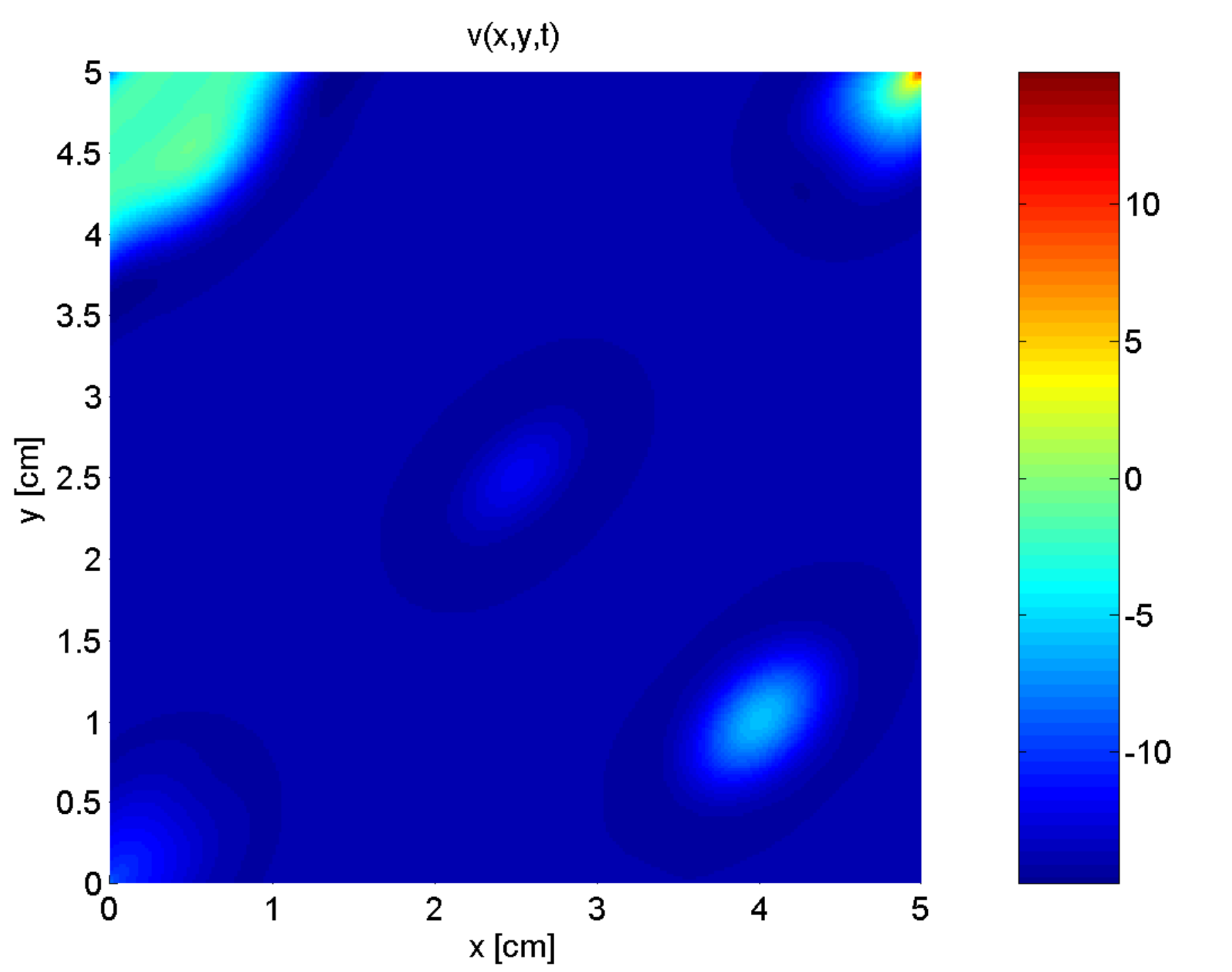}\\
 \includegraphics[width=0.42667\textwidth]{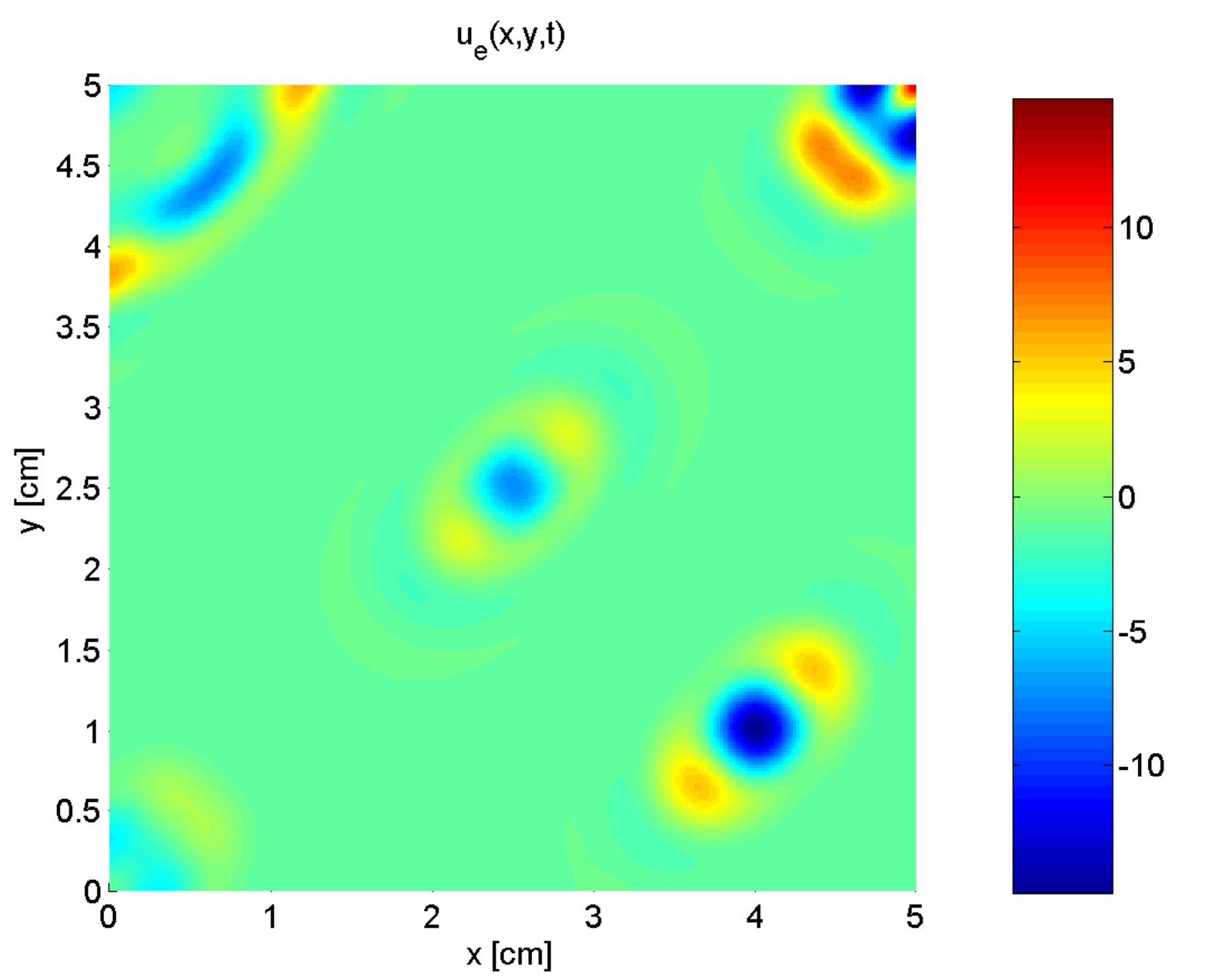}&
 \includegraphics[width=0.42667\textwidth]{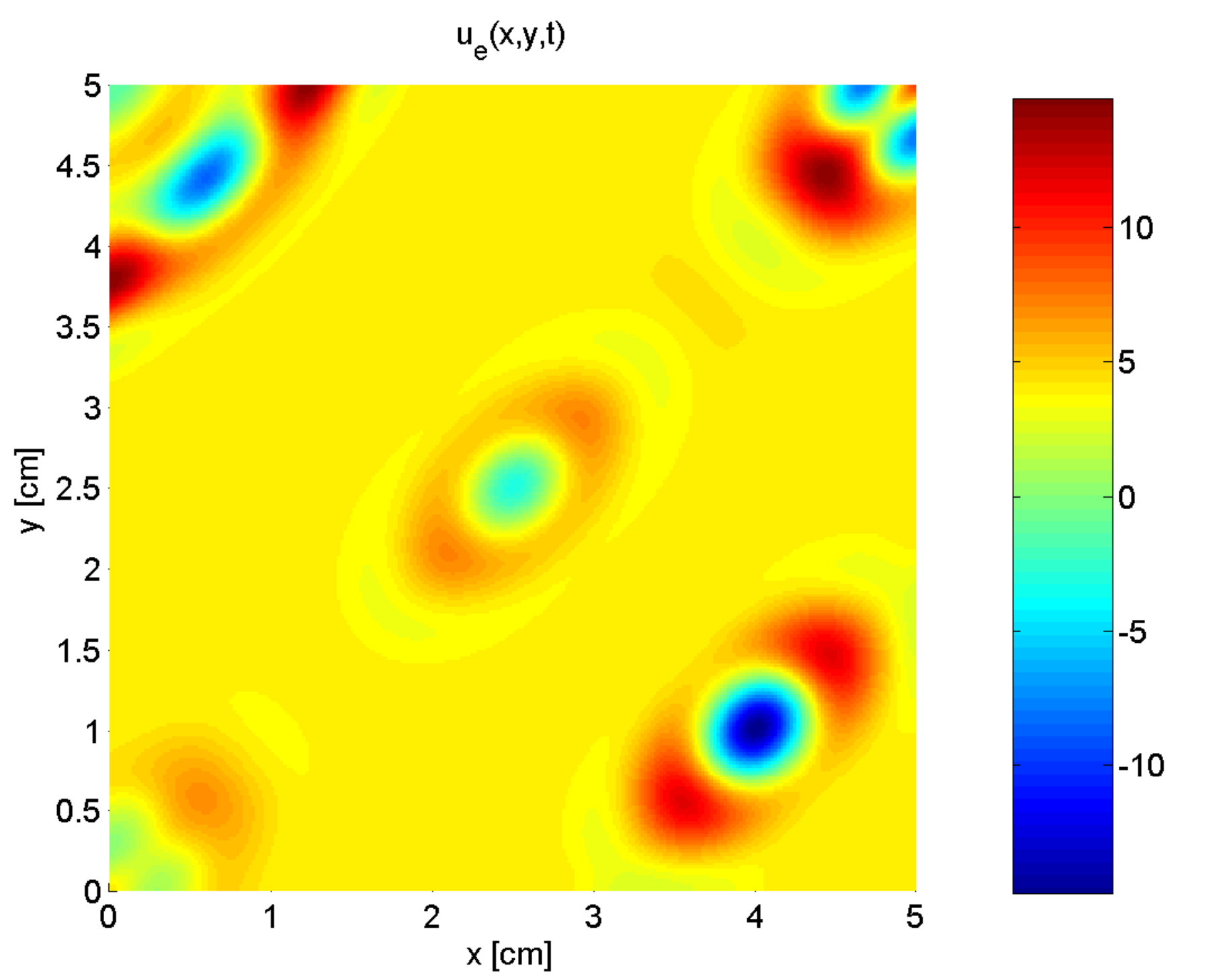}\\
 \includegraphics[width=0.322667\textwidth]{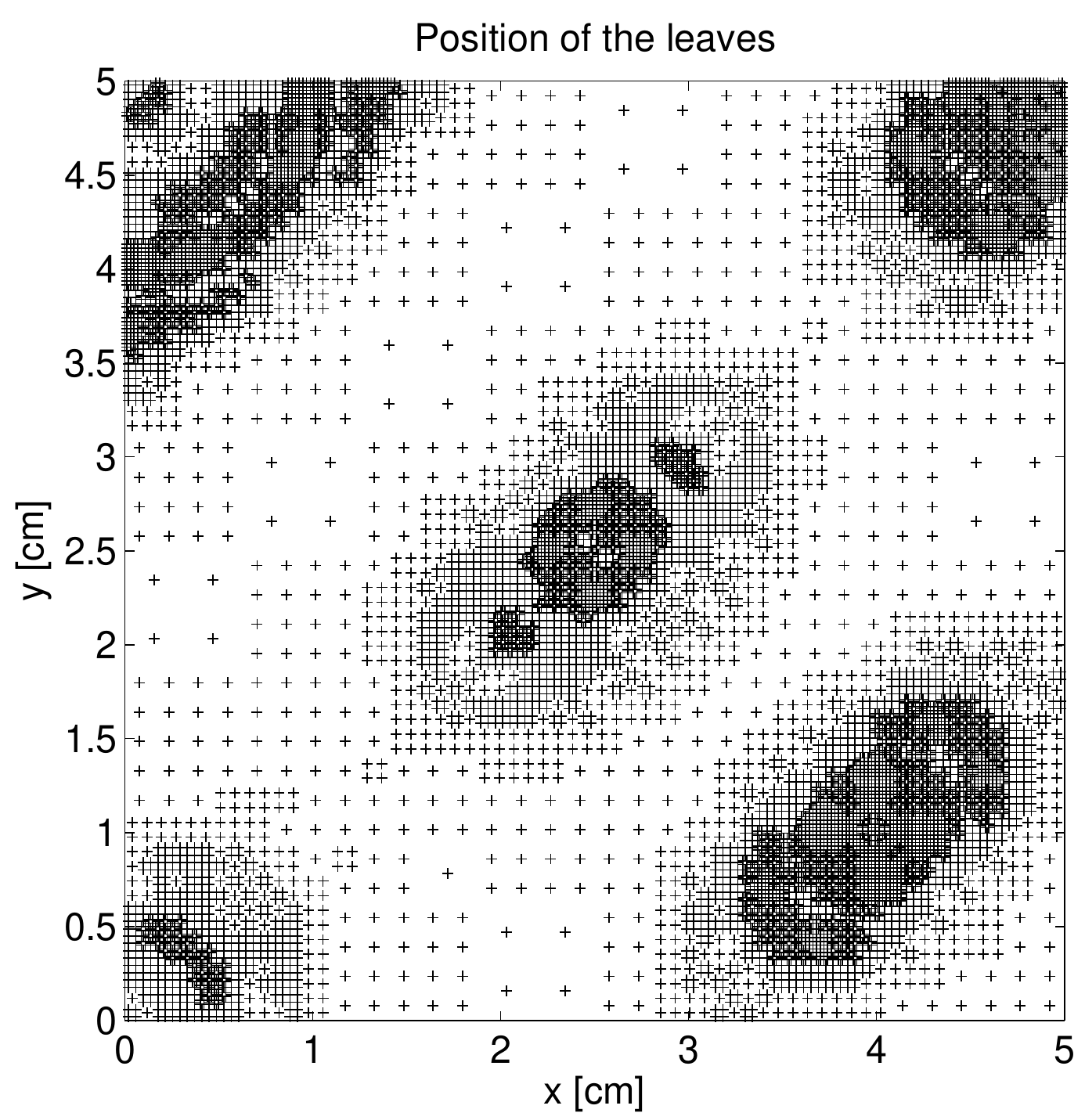}
\quad \qquad$\vphantom{X}$ $\vphantom{X}$ &
 \includegraphics[width=0.322667\textwidth]{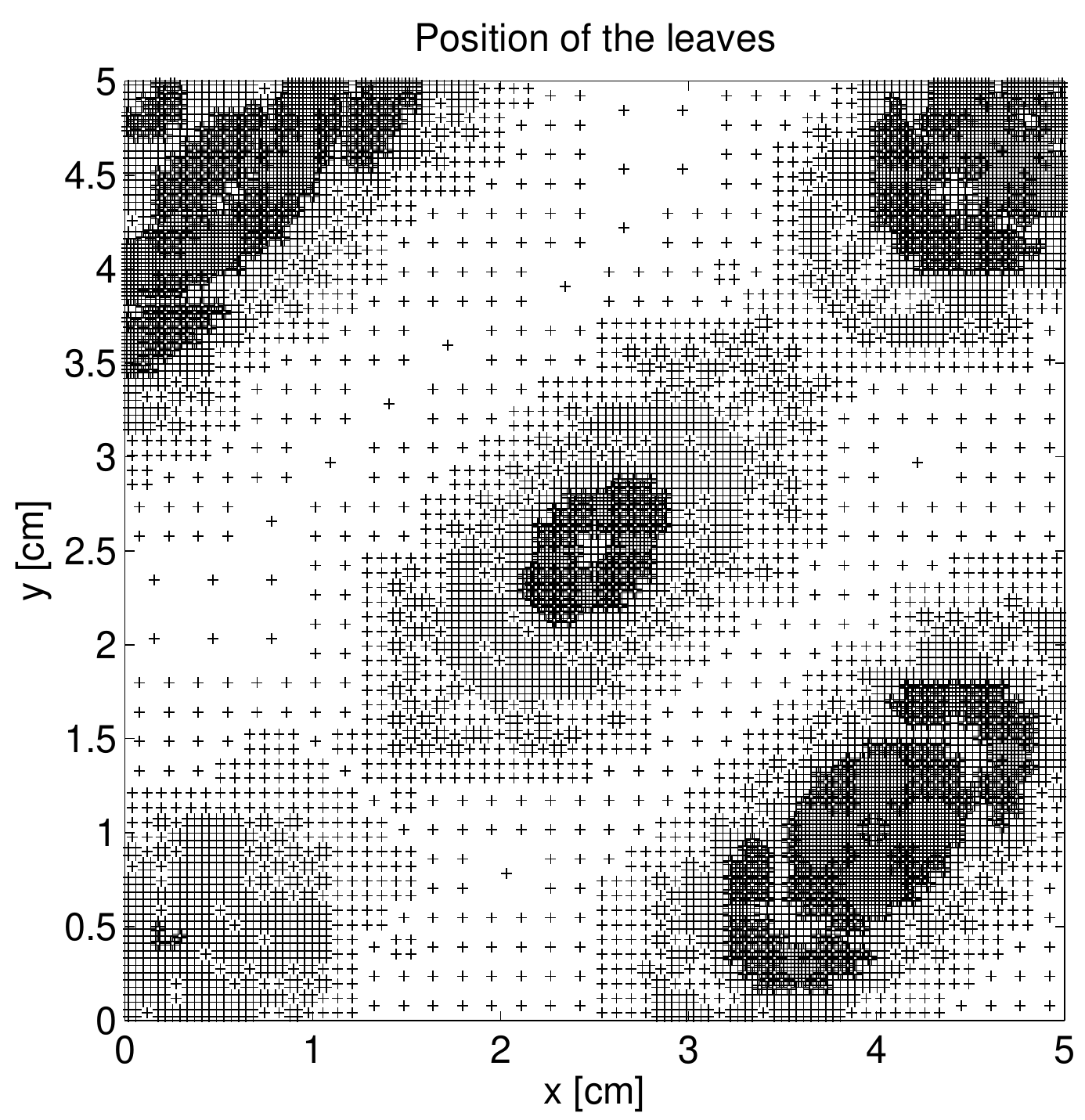}
\quad \qquad$\vphantom{X}$ $\vphantom{X}$
 \end{tabular}
 \caption{Example~3 (bidomain model, three stimuli):
 Numerical solution for transmembrane potential $v$ and
 extracellular potential $u_\mathrm{e}$ in $[\mathrm{mV}]$,
and leaves of the corresponding
 tree data structure at times $t=2.0\,\mathrm{ms}$ and
$t=5.0\,\mathrm{ms}$.} \label{fig:snapshots4}
 \end{center}
 \end{figure}

We will present three test cases showing the efficiency of the previously
described  methods in capturing the dynamical evolution of electro-physiological
waves for both the monodomain and bidomain models. Since we are dealing with
multicomponent solutions, we emphasize  that  a single
mesh is used to represent the vector of relevant variables.
In the bidomain model,
the anisotropies, mesh structures, and  the size of the problem cause the
sparse linear system corresponding to \eqref{S2-discr} to be ill-conditioned.
This system needs to be solved in each time step, which is
done by   the
Cholesky method. Before presenting the numerical results, we
 provide further  details on the implementation of the numerical
schemes.

\subsection{Implementation Issues}
 The following algorithm  shows how the solution
$\u^{n+1}=(v,u_\mathrm{e},w)^{n+1}$ is obtained in each time step
is obtained for each time step
\begin{alg}[General method]\label{alg:general} \hfill
\begin{enumerate}
\item Assume that  $u_\mathrm{i}^n$, $u_\mathrm{e}^n$, $v^n$ and $w^n$
 are known (at
time   $t^n$).
   \item Solve the ODE
   $$\pt w -H(v,w) =0,  \quad  x \in \Om,$$
 approximately
for $t^n<t\leq t^{n+1}$ with initial condition $w^n$ and data $v^n$, i.e.,
compute $w^{n+1}$ using  \eqref{S3-discr}.
\item Solve the parabolic PDE
\begin{align*}
\beta c_{\mathrm{m}}\pt v+\Div\bigl(\bM_\mathrm{e}(x)\Grad
u_\mathrm{e}\bigr)+\beta \Ion(v,w)& =
\Iap, \quad  x \in \Om, \\
\bigl(\bM_\mathrm{e}(x)\Grad u_\mathrm{e})\cdot\n &=0 \quad
 \text{\em on $\partial\Om$}
\end{align*}
approximately
 for $t^n<t\leq t^{n+1}$, with $v(t^n)=v^n$ and $w(t^n)=w^n$, i.e.,
calculate $v^{n+1}$ using \eqref{S1-discr}.
\item Solve the elliptic problem
\begin{align*}
\Div \bigl((\bM_\mathrm{i}(x)+\bM_\mathrm{e}(x))\Grad u_\mathrm{e}\bigr)+
\Div\bigl(\bM_\mathrm{i}(x)\Grad v\bigr)&=\Iap,\quad  x \in \Om\\
\bigl(\bM_j(x)\Grad u_j)\cdot\n &=0\quad \text{\em on $\partial\Om$, $j \in \{ \mathrm{e},\mathrm{i}\}$}
\end{align*}
approximately  for $t^n<t\leq t^{n+1}$ with
$v(t^n)=v^n$ and $u_\mathrm{e}(t^n)=u_\mathrm{e}^n$, i.e.,
determine $u_\mathrm{e}^{n+1}$ by solving the linear system
\eqref{S2-discr}.
\end{enumerate}
\end{alg}

This algorithm structure is usually preferred for systems involving parabolic
and elliptic equation, since  it explicitly isolates the
solution of the elliptic problem from the rest of the computations \cite{pennacc}.

For multicomponent solutions, there are  many  possible definitions for a
scalar detail $\smash{d_{(i,j),l}}$ that  is calculated from the details of the
components (see a brief discussion in \cite{brss}). To guarantee
that the  refinement and coarsening procedures are always on the safe side,
in the sense that we always prefer to keep a position with a detail triple
containing at least one component above the threshold \eqref{eq4.4},
we will use
$\smash{d^{\mathbf{u}}_{(i,j),l}=\min \{{d}^v_{(i,j),l},
d^{u_\mathrm{e}}_{(i,j),l},d^w_{(i,j),l} \}}$
 and  $\smash{d^{\mathbf{u}}_{(i,j),l}=\max \{d^v_{(i,j),l},
d^{u_\mathrm{e}}_{(i,j),l},d^w_{(i,j),l} \}}$
for the refinement and coarsening procedures, respectively.
In practice, the details introduced in
Section~\ref{subsec4.2} are computed  simply as the differences between the "exact" and the
predicted value:
\begin{align*}
 d^\mathbf{u}_{(i,j),l}:=\mathbf{u}_{(i,j),l}-\hat{\mathbf{u}}_{(i,j),l}.
\end{align*}

In \eqref{equ:epsref1}  $\alpha$ stands for the accuracy order of the
FV reference scheme, and numerical experience gives $\alpha=1.09$. To choose an
acceptable value for the factor $C$, a series
of computations (not completely shown here) with different tolerances are
needed in each case, prior to final computations. We basically choose
the largest available candidate value for $C$ such that  the same order of
accuracy as that of the reference FV scheme (same slope) is maintained
and the data compression rate~$\eta$ and the speed-up~$\mathcal{V}$
  are maximized (see Figure \ref{fig:factorC}).
In \cite{bbrs} we  give a detailed description of the multiresolution
algorithm and the LTS algorithm for systems of reaction-diffusion
equations.
For sake of completeness, we provide versions of  these algorithms
 adapted to the bidomain model of electrocardiology in the Appendix.


\subsection{Example 1} For this example, we consider the simple monodomain
model \eqref{monodomain} with homogeneous Neumann boundary conditions. The
ionic  current and membrane model is determined by  the FitzHugh-Nagumo
membrane kinetics \eqref{fhn}, with $a=0.16875$, $b=1.0$, $\lambda=-100$
and $\theta=0.25$. The computational domain is the square $\Om=[0,1\ \text{cm}]^2$,
and the remaining parameters are $c_{\mathrm{m}}=1.0\,\mathrm{mF/cm}^2$ and
$\beta=1.0\,\mathrm{cm}^{-1}$. The units for $v,w$ are $\mathrm{mV}$.
We consider in \eqref{monodomain}
 $\smash{(1+\lambda)^{-1} \bM_\mathrm{i}:=\diag (\gamma, \gamma)}$ with
 $\gamma=0.01$.
The respective initial data for~$v$ and~$w$ are
\begin{align*}
v^0(x,y)=\left(1-\frac{1}{1+ \exp( -50 (x^2+y^2)^{1/2} -0.1) } \right)
\,\mathrm{mV},\quad w_0=0\,\mathrm{mV}.
\end{align*}
After $4\,\mathrm{ms}$,
an instantaneous stimulus is applied in $(x_0,y_0)=(0.5\,\mathrm{cm},
 0.5\,\mathrm{cm})$ to the membrane potential $v$
 $$\frac{\lambda}{1+\lambda}\Iap:=\begin{cases}
 1\,\mathrm{mV} & \text{if $(x-x_0)^2+(y-y_0)^2<0.04\, \mathrm{cm}^2$,} \\
 0\,\mathrm{mV} & \text{otherwise.}
 \end{cases}.$$
In this example,  we use $L=10$ resolution levels, $\mathcal{N}=262144$ elements
in the finest level, a tolerance of
$\eps_{\mathrm{R}}=1\times10^{-3}$,
and we compute normalized errors by comparison with a reference
solution obtained with a fine mesh calculation with $\mathcal{N}=1024^2=1048576$
control volumes.  The time evolution is made using a first-order explicit Euler scheme.
Plots of the numerical solution with the corresponding adaptively refined meshes
at different times are shown in Figures~\ref{fig:monodomain1} and~\ref{fig:monodomain2}.

As can be seen from Table~\ref{table:ex1}, the normalized errors are controlled to be
of the same order of the reference tolerance $\eps_{\mathrm{R}}$. We also see that the
MR algorithm is efficient: we have high rates of memory compression
and speed-up.

\subsection{Example 2}  In Examples~2 and~3, we present computational
results for  the simulation of the bidomain equations. We
 consider  a computational
 domain $\Om=[0,5\,\mathrm{cm}]^2$,  and the parameters in  \eqref{Ion}
 and
 \eqref{S4}
(after \cite{entcheva,stm:00,ying,yrh:07}) are  given by the
 membrane capacitance $c_{\mathrm{m}}=1.0\,
\mathrm{mF/cm}^2$, the intracellular conductivity in the principal axis
         $\smash{\sigma_\mathrm{i}^{\mathrm{l}}}=6\,\Om^{-1}\mathrm{cm}^{-1}$, the
remaining intracellular  conductivity
 $\smash{\sigma_\mathrm{i}^{\mathrm{t}}}=0.6\,\Om^{-1}\mathrm{cm}^{-1}$
 (corresponding to  an anisotropy ratio of 10), the   extracellular conductivities
$\sigma_\mathrm{e}^{\mathrm{l}}=24\, \Om^{-1}\mathrm{cm}^{-1}$ and
 $\sigma_\mathrm{e}^{\mathrm{t}}=12\,\Om^{-1}\mathrm{cm}^{-1}$
 (corresponding to an anisotropy ratio of 2), the surface-to-volume ratio   $\beta=2000\,
 \mathrm{cm}^{-1}$, the surface resistivity $R_{\mathrm{m}}=2\times10^{4}\,\Om\,\mathrm{cm}^2$,
 $v_p=100\,\mathrm{mV}$, $\eta_1=0.005$, $\eta_2=0.1$, $\eta_3=1.5$,
 $\eta_4=7.5$,  and   $\eta_5=0.1$. The fibers   form an angle of
 $\pi/4$ with the $x$-axis.

In Example~2,  the initial datum is given by  a stimulus applied on the extracellular potential
$u_\mathrm{e}$ in the center of the domain, while  both~$v$ and the gating
 variable~$w$  are initially set to zero
(see Figure~\ref{fig:cond_ini}). The units for $v,u_e$ and $w$ are
$\mathrm{mV}$. In this example, the following MR setting is chosen. We utilize
  wavelets with $r=3$ vanishing moments, a maximal resolution level $L=9$, and
  therefore a finest mesh with $\mathcal{N}=65536$ elements. The  reference tolerance
  given by $\eps_{\mathrm{R}}=5.0\times 10^{-4}$.

We show in Figures~\ref{fig:snapshots1} and~\ref{fig:snapshots2} a sequence of
snapshots after an initial stimulus applied to the center of the domain,
corresponding to transmembrane potential $v$, extracellular potential
$u_\mathrm{e}$ and adaptive mesh.

Table~\ref{table:ex2}  illustrates  the efficiency and accuracy
of the base    MR   method by tabulating CPU ratio~$\mathcal{V}$,
compression rate~$\eta$  and normalized errors.   By using MR, we
obtain an average data compression rate of 17 and an increasing  speed-up rate
up to 26.09. Moreover, the errors in three different norms remain of the order
of $\eps_{\mathrm{R}}$. Here we have computed normalized errors using a  reference
FV solution on a grid with $\mathcal{N}=1024^2=1048576$ control volumes.

For the time integration using the  LTS method, we     choose the maximum
CFL number allowed by \eqref{cfl},   $\mathrm{CFL}_{l=0}=0.5$ for the coarsest
level and $\mathrm{CFL}_l=2^l\mathrm{CFL}_{l=0}$ for finer levels. For the RKF
computations, we use $\delta_\mathrm{desired}=1\times10^{-4}$,
$\mathcal{S}_0=0.1$, $\mathcal{S}_\mathrm{min}=0.01$,  and the initial CFL condition
$\mathrm{CFL}_{t=0}=0.5$.

We select this example for a detailed comparison of the
 performance of the  FV and MR
 methods with a global time step, the MR method with RKF
 adaptive global time stepping (MR-RKF), and the  MR method
 with local time stepping (MR-LTS). The evolution of the  speed-up
factor~$\mathcal{V}$ the   and
 data compression rate~$\eta$  for the MR versions and of the
normalized $L^1$ and $L^{\infty}$ errors for all these methods
 are   displayed in  Figure~\ref{fig:varios}. From these plots
it is observed that with RKF and LTS, the data compression rate is of the same
order during the time evolution, which means  that the adaptive meshes
 for both
 methods should be roughly the same. Also, a  substantial additional   gain is
obtained in speed-up rate when comparing with a MR calculation     using
global time stepping: The MR-LTS method gives us an additional speed-up factor
of about~2, while with the RKF alternative we obtain an additional speed-up of
about~4. This effect could be explained in part from the lack of
need of a synchronization procedure for the RKF computations,  and the fact that
the CFL condition \eqref{cfl} is not imposed during the time evolution with the
MR-RKF method, allowing larger time steps. (Although condition \eqref{cfl} guarantees numerical stability of the
solutions, in practice  this  is  observed to be a fairly conservative
 estimates, and moderately larger time steps may  be used.)  We can
also conclude that the errors of the MR-LTS  computations are  kept of the same
order that the errors obtained with a global time stepping, while the incurred
errors by using the MR-RKF method are larger during the whole time evolution.

\subsection{Example 3} For this example, we consider an initial
stimulus at the center of the domain, later at $t=0.2\,\mathrm{ms}$ we apply
another instantaneous stimulus to the northwest corner of the domain,
and
then  at $t=1.0\,\mathrm{ms}$ we apply a third stimulus of the same magnitude to the
northeast and southwest corners. The system is evolved and we show snapshots
of the numerical solution for $v,u_\mathrm{e}$ and the adaptive mesh. We use
the MR-RKF method with $\mathcal{N}=65536$, $\eps_{\mathrm{R}}=2.5\times10^{-3}$,
$\delta_\mathrm{desired}=1\times10^{-3}$ and the remaining parameters are
considered as in Example~2. As in Example~2, from Figures~\ref{fig:snapshots3}
and~\ref{fig:snapshots4} we clearly notice the anisotropic
orientation of the fibers.

\section{Conclusions} \label{sec:concl}
We address the application of a MR method for FV
schemes combined   with  LTS and RKF
adaptive time stepping for solving the bidomain equations. The
numerical experiments illustrate that these
 methods
are  efficient  and  accurate    enough  to
simulate the electrical activity in myocardial tissue with affordable
 effort. This is a real advantage in comparison with more
 involved methods that require large scale computations on clusters.
 We here contribute to the recent work done by several groups in testing
 whether the combination of MR, LTS and RKF strategies is indeed effective
 for a relevant class of problems.

From     a     numerical point of view, the
plateau-like structures, associated with very steep gradients,
 of typical solutions motivate the
use of a   locally  refined    adaptive mesh, since we require  high
resolution near these steep   gradients only. These
 areas of strong variation occupy   a very reduced
   part  of  the   entire  domain only, especially in the case of sharp fronts.
Consequently our gain will be less  significant in the presence of chaotic
electrical activity or when multiple waves interact in the considered tissue.

Based on our numerical examples, we conclude that using a LTS
strategy, we obtain a substantial gain in CPU time   speed-up for a factor of
about 2 for larger scales while the errors between the   MR-LTS solution and a
reference solution are of the same order as  those of
  the MR   solution. On the other
hand, using an MR-RKF strategy, we obtain  an  additional   speed-up factor of
about 4, but at the price of larger errors. However, in assessing our findings,
it is important to recognize limitations. The high rates of compression obtained
with our methods are problem-dependent and they may  depend
on the proper adjustment
  of parameters. We have only considered here very simple geometries,
because all computations are concentrated on adaptivity and performance.
Simulations on more complex and realistic geometries are part of possible future
work.

Finally, we remark that the FV method given in Section~\ref{sec:FV}
as well the MR framework detailed in Section~\ref{sec:MR} are both
straightforwardly extensible to the 3D case. A convergence analysis of
the  implicit version of the FV
 method presented in Section~\ref{sec:FV} is  being prepared
 \cite{bbrconvergence}.



\section*{Acknowledgements}

MB acknowledges support by Fondecyt project 1070682, RB acknowledges
support by Fondecyt project 1050728 and Fondap in Applied Mathematics,
 project 15000001 and RR acknowledges support by Conicyt Fellowship
and Mecesup project UCO0406.

\section*{Appendix}\label{appendix_1}
\appendix

For illustrative purposes, suppose that we are in the case of a cartesian mesh.
We give here an example of an interior
first-order flux calculation using LTS for the parabolic part of the bidomain
scheme (the equation for $v$),
needed to complete a full macro time step, by the following algorithm:

\begin{alg}[Locally varying intermediate time stepping] \hfill
\begin{enumerate}
   \item Grid adaptation (provided the former sets of leaves).
   \item {\bf do} $k=1,\ldots,2^L$ (for the local time steps
      $n+2^{-L},n+2\cdot 2^{-L}, n+ 3\cdot 2^{-L},\ldots,  n+1$)
   \begin{enumerate}
      \item Synchronization:
      \item[]{\bf do} $l=L,\ldots,1$
      \begin{itemize}
         \item[]{\bf do} $i=1, \dots, |\tilde{\Lambda}|_x(l)$, $j=1, \dots,
           |\tilde{\Lambda}|_y(l)$
         \begin{itemize}
            \item[] {\bf if} $1\leqslant l\leqslant \tilde{l}_{k-1}$
              {\bf then}
            \begin{itemize}
                \item[]  {\bf if} $V_{(i,j),l}$ is a virtual leaf {\bf then}
		  \begin{itemize}
		    \item[] ${\displaystyle
                  \bar{F}^{n+k2^{-L}}_{(i,j),l\to(i+1,j),l} \leftarrow
               \bar{F}^{n+(k-1)2^{-L}}_{(i,j),l\to(i+1,j),l}\vphantom{\int} }$
                 \item[]  Update reaction terms:

                 ${\displaystyle
                   I_{\mathrm{ion}(i,j),l}^{n+k2^{-L}} \leftarrow
                   I_{\mathrm{ion}(i,j),l}^{n+(k-1)2^{-L}}, \quad
                   I_{\mathrm{app}(i,j),l}^{n+k2^{-L}} \leftarrow
                   I_{\mathrm{app}(i,j),l}^{n+(k-1)2^{-L}}\vphantom{\int}}$
            \end{itemize}
		  \item[] {\bf endif}
            \end{itemize}
            \item[]{\bf else}
            \begin{itemize}
	        \item[] {\bf if} $V_{(i,j),l}$ is a leaf {\bf then}
		  \begin{itemize}
                \item[] ${\displaystyle
                   I_{\mathrm{ion}(i,j),l}^{n+k2^{-L}} \leftarrow \Ion
                   \bigl(v_{(i,j),l}^{n+k2^{-L}},w_{(i,j),l}^{n+k2^{-L}}
                   \bigr)}$,
                \item[] ${\displaystyle I_{\mathrm{app}(i,j),l}^{n+k2^{-L}} \leftarrow \Iap
                   \bigl(v_{(i,j),l}^{n+k2^{-L}},w_{(i,j),l}^{n+k2^{-L}}
                   \bigr)\vphantom{\int}}$
                \item[] {\bf if}
                     $V_{(i+1,j),l}$ is a leaf {\bf then}
                    \item[] $\quad {\displaystyle  \bar{F}_{(i,j),l\to(i+1,j),l}
                        \leftarrow -\frac{d^*_{\mathrm{e}}}{h(l)}\frac{|\sigma_{V
                        _{(i,j),l},V_{(i+1,j),l}}|}{d(V_{(i,j),l},V_{(i+1,j),l})} (u_{\mathrm{e},(i+1,j),l}
                        -u_{\mathrm{e},(i,j),l})\vphantom{\int}}$
                    \item[] $\quad {\displaystyle  \bar{F}_{(i,j),l\to(i,j+1),l}
                        \leftarrow -\frac{d^*_{\mathrm{e}}}{h(l)}\frac{|\sigma_{V
                        _{(i,j),l},V_{(i,j+1),l}}|}{d(V_{(i,j),l},V_{(i,j+1),l})}
                        (u_{\mathrm{e},(i,j+1),l}-u_{\mathrm{e},(i,j),l})\vphantom{\int}}$
                \item[]{\bf endif}
		\item[] {\bf if}
                    $V_{(2i+2,2j),l+1}$, $V_{(2i+2,2j+1),l+1}$ are leaves (interface
                   edges) {\bf then}
                     \item[] $\quad {\displaystyle\vphantom{\int}
                          \bar{F}_{(i,j),l\to(i+1,j),l}\leftarrow
                          \bar{F}_{2(i+1,j),l+1\to(2i+1,2j),l+1}
                         +\bar{F}_{(2i+2,2j+1),l+1\to(2i+1,2j+1),l+1}
                          \vphantom{\int}}$
                     \item[] $\quad {\displaystyle\vphantom{\int}
                          \bar{F}_{(i,j),l\to(i,j+1),l}\leftarrow
                          \bar{F}_{2(i,j+1),l+1\to(2i,2j+1),l+1}
                         +\bar{F}_{(2i+1,2j+2),l+1\to(2i+1,2j+1),l+1}
                          \vphantom{\int}}$
                \item[] {\bf endif}
            \end{itemize}
  \item[] {\bf endif}
            \end{itemize}
            \item[] {\bf endif}
         \end{itemize}
         \item[]{\bf enddo}
      \end{itemize}
      \item[]{\bf enddo}
   \item Time evolution:
      \item[]{\bf do} $l=1,\ldots,L$, $i=1, \dots, |\tilde{\Lambda}|_x(l)$,
      $j=1, \dots, |\tilde{\Lambda}|_y(l)$
      \begin{itemize}
         \item[] {\bf if} $1\leqslant l\leqslant \tilde{l}_{k-1}$ {\bf then}
             there is no evolution:
         \begin{itemize}
             \item[] ${\displaystyle v_{(i,j),l}^{n+(k+1)2^{-L}}
                \leftarrow v_{(i,j),l}^{n+k2^{-L}}}$
         \end{itemize}
         \item[]{\bf else}
         \begin{itemize}
            \item[] Marching formula only for the leaves $V_{(i,j),l}$:
            \item[] ${\displaystyle v_{(i,j),l}^{n+(k+1)2^{-L}}\leftarrow \frac{1}{c_m}v_{(i,j),l}^{n+k2^{-L}}+
                  \frac{\Delta t_l}{c_m}I_{\mathrm{app}(i,j),l}^{n+k 2^{-L}}-
                  \frac{\Delta t_l}{c_m}I_{\mathrm{ion}(i,j),l}^{n+k 2^{-L}}}$
            \item[] \qquad \qquad ${\displaystyle - \frac{\Delta t_l}{c_m}\sum_{(n,m)\in \mathcal{S}((i,j),l)}\frac{d^*_{\mathrm{e}}}{h(l)}\frac{|\sigma_{V
                        _{(i,j),l},V_{(n,m),l}}|}{d(V_{(i,j),l},V_{(n,m),l})}
                        (u_{\mathrm{e},(n,m),l}-u_{\mathrm{e},(i,j),l})\vphantom{\int}}$
           \end{itemize}
         \item[] {\bf endif}
      \end{itemize}
      \item[]{\bf enddo}
      \item Partial grid adaptation each odd intermediate time step:
      \item[] {\bf do} $l=L,\ldots,\tilde{l}_k+1$
	\begin{itemize}
     \item[] Projection from the leaves.
       \end{itemize}
              {\bf enddo}
      \item[] {\bf do} $l=\tilde{l}_k,\ldots,L$ 	\begin{itemize}
     \item[] Thresholding, prediction, and
        addition of the safety zone.
	  \end{itemize}
	  {\bf enddo}
   \end{enumerate}
   \item[]{\bf enddo}
\end{enumerate}
\end{alg}
Here, $\tilde{l}_k$ denotes the coarsest level containing leaves in the intermediate step $k$,
$h(l)$ is the mesh size on level $l$, and $|\tilde{\Lambda}|_z(l)$ is the size of the set formed
by leaves and virtual leaves per
resolution level $l$ in the direction $z$. The marching formula  corresponds to
\eqref{S1-discr}, for the intermediate time steps $k=1,\ldots,2^L$, for the leaf $V_{(i,j),l}$.

Now we give a brief description of the general multiresolution procedure.
\begin{alg}[Multiresolution procedure] \hfill
\begin{enumerate}
\item Initialization of parameters.
\item Creation of the initial tree:
  \begin{enumerate}
      \item Create the root and compute its cell average value.
      \item Split the cell, compute the cell average values in the sons and
compute details.
      \item Apply thresholding for the splitting of the
	new sons.
      \item Repeat this until all sons have details below the required
tolerance $\varepsilon_l$.
  \end{enumerate}
\item {\bf do}  $n=1, \dots,  total\_ time\_ steps$
\begin{enumerate}
\item Determination of the leaves and virtual leaves sets.
\item Time evolution with global time step: \label{xxx}
Compute the discretized space operator $\mathcal{A}$  for all the leaves.
\item Updating the tree structure:
  \begin{itemize}
      \item Recalculate the values on the nodes by
projection from the leaves. Compute the details for all positions $(\cdot,\cdot,l)$
for $l\geqslant \tilde{l}_k$. If the detail
 in a node and in its brothers is smaller than
the prescribed tolerance, then the cell and its brothers are \emph{deletable}.
      \item If some node and all its sons are deletable, and the sons
are leaves without virtual sons, then delete sons. If this node has no
sons and it is not deletable and it is not at level $l=L$, then create sons.
      \item Update the values in the new sons by prediction from the
former leaves.
\end{itemize}
\end{enumerate}
\item[] {\bf enddo}
\item Output: Save meshes, leaves and cell averages.
\end{enumerate}
\end{alg}
Here  $total\_ time\_ steps$ is the total time steps needed to reach
$T_{\text{final}}$ using $\Delta t$ as the maximum time step allowed by the
CFL condition using the finest space step.

When using a RKF strategy for the time evolution, replace step (\ref{xxx}) by
the new step
\begin{enumerate}
\item[(3)]
\begin{itemize}
\item \emph{Compute the discretized space operator $\mathcal{A}$ for all the
leaves as in \eqref{kappas}}.
\item \emph{Compute the difference between the two solutions obtained as in
\eqref{delta_old}}.
\item \emph{Apply the limiter for the time step variation and compute the new
time step by \eqref{Delta t_new}}.
\end{itemize}
\end{enumerate}

When using a LTS strategy, replace step (\ref{xxx}) by the new step
\begin{enumerate}
\item[(3)]  {\bf do}  $n=1, \dots,  total\_ time\_ steps$
\begin{enumerate}
\item \emph{Determination of the leaves and virtual leaves sets}.
\item \emph{Time evolution with local time stepping:
 Compute the discretized space operator $\mathcal{A}$ for all the leaves and
virtual leaves}
\item {\bf do} $k=1, \dots , 2^L$ \emph{($k$ counts intermediate time
  steps)}
\begin{itemize}
\item {\em Compute the
intermediate time steps depending on the position of the
 leaf as explained in Section~\ref{sec:LTS}.}
\item {\bf if} $k$ is odd {\bf then}  \emph{update the tree
  structure:}
\begin{itemize}
      \item \emph{Recalculate the values on the nodes and the virtual nodes by
projection from the leaves. Compute the details in the whole tree. If the detail
in a node is smaller than the prescribed
tolerance, then the cell and its brothers are \emph{deletable}}.
      \item \emph{If some node and all its sons are deletable, and the sons
are leaves without virtual sons, then delete sons. If this node has no
sons and it is not deletable and it is not at level $l=L$, then create sons.}
      \item \emph{Update the
values in the new sons by prediction from the
former leaves.}
\end{itemize}
\item[] {\bf endif}
\end{itemize}
\item[] {\bf enddo}
\item[] \emph{(Now, after $2^L$ intermediate
  steps,  all the elements are synchronized.)}
\end{enumerate}
\item[] {\bf enddo}
\end{enumerate}
Here  $total\_ time\_ steps$ is the total time steps needed to reach
$T_{\text{final}}$, with $\Delta t_0$ as the maximum time step allowed by the
CFL condition using the coarsest space step.

\end{document}